\documentclass[a4paper,12pt]{book}
\usepackage{a4wide}
\usepackage{amsmath}
\usepackage{amssymb}
\usepackage{amsthm}
\usepackage{latexsym}
\usepackage{graphicx}
\usepackage[english]{babel}
\usepackage{makeidx}
                     
\newtheorem{dfn} [subsection]{Definition}
\newtheorem{obs} [subsection]{Remark}
\newtheorem{rem} [subsection]{Remark}
\newtheorem{exm} [subsection]{Example}
\newtheorem{ex}  [subsection]{Example}

\newtheorem{prop}[subsection]{Proposition}

\newtheorem{teor}[subsection]{Theorem}
\newtheorem{teo}[subsection]{Theorem}
\newtheorem{lema}[subsection]{Lemma}
\newtheorem{cor} [subsection]{Corollary}

\newcommand{\Nn}{\mathbb N^{n}}

\newcommand{\de}{\mathbf{d}}
\newcommand{\me}{\mathbf{m}}
\newcommand{\en}{\mathbf{n}}
\newcommand{\Gin}{Gin}
\newcommand{\Dei}{$\mathcal D$-fixed ideal }

\newcommand{\Deisv}{$\mathcal D$-fixed ideals, }
\newcommand{\Di}{$\mathcal D$-fixed ideal }

\newcommand{\di}{$\de$-fixed ideal }
\newcommand{\dis}{$\de$-fixed ideals }
\newcommand{\leqd}{\leq_{\mathbf{d}}}
\begin{document}

\selectlanguage{english}
\frenchspacing

\huge
\begin{center}
\textbf{Contributions in combinatorics in commutative algebra}
\vspace{20 pt}

Ph.D. Thesis
\large
\vspace{20 pt}

Mircea Cimpoea\c s

\vspace{20 pt}

Adviser: Professor dr. Dorin Popescu

\vspace{200 pt}

University of Bucharest,\\
Faculty of Mathematics and Informatics \\

\vspace{20 pt}
October 2007
\end{center}

\normalsize





\newpage

.

\newpage
\section*{Abstract}

In the first chapter we present new results related on monomial ideals of Borel type. Also, we introduce a new
class of monomial ideals, called $\de$-fixed ideals, which generalize the class of $p$-Borel ideals and we extend
several results to this new class. In section $1.2$ we extend a result of Eisenbud-Reeves-Totaro in the frame of ideals of Borel type. This allows us to obtain some nice consequences related to the regularity of the Borel type ideals.
In section $1.3$ we introduce a new class of ideals, called strong Borel type ideals, and we compute the Mumford-Castelnouvo regularity for a special case of strong Borel type ideals. In the sections $1.4$, $1.5$ and $1.6$ we show how some results for $p$-Borel ideals can be transfered to $\de$-fixed ideals. In particular, we give the form of a principal $\de$-fixed ideal and we compute the socle of factors of these ideals, using methods similar as in \cite{hpv}, see section $1.5$. This allows us to give a generalization of Pardue's formula, i.e. a formula of the regularity for a principal $\de$-fixed ideal, see section $1.6$. In the last section of the first chapter, we describe the $\de$-fixed ideals generated by powers of variables.

In the second chapter, we compute the generic initial ideal, with repect to the reverse lexicographic order, of an ideal which define a complete intersection of embedding dimension three with strong Lefschetz property and we show that it is an almost reverse lexicographic ideal, see sections $2.2$ and $2.3$. This enable us to give a proof for Moreno's conjecture in the case $n=3$ and characteristic zero, see section $2.1$. In section $2.4$ we prove that the $d$-component of the generic initial ideal, with respect to the reverse lexicographic order, of an ideal generated by a regular sequence of homogeneous polynomials of degree $d$ is revlex, in a particular, but important, case. Using this property, in section $2.5$ we compute the generic initial ideal for several complete intersections with strong Lefschetz property.

\section*{Acknowledgements.}

The author is deeply indebted to his adviser Professor Dorin Popescu, for all of his help. The author owes a sincerely thank to Assistant Professor Alin \c Stefan for discussions on the section $1.7$. Also, he owes a special thank to Dr. Marius Vladoiu for his help on the sections $2.2$ and $2.3$. My thanks goes also to the School of Mathematical Sciences, GC University, Lahore, Pakistan, especially to Dr. A. D. Raza Choudary, for supporting and facilitating an important part of the work on my thesis.

\newpage
\huge \noindent \textbf{Contents.} \vspace{30 pt} \normalsize

\contentsline {chapter}{\numberline {}Abstract.}{3}
\contentsline {chapter}{\numberline {}Acknowledgements.}{3}
\contentsline {chapter}{\numberline {}Contents.}{4}
\contentsline {chapter}{\numberline {}Introduction.}{5}
\contentsline {chapter}{\numberline {1}Ideals of Borel type and $\mathbf {d}$-fixed ideals.}{9}
\contentsline {section}{\numberline {1.1}Ideals of Borel type.}{9}
\contentsline {section}{\numberline {1.2}Stable properties of Borel type ideals.}{11}
\contentsline {section}{\numberline {1.3}Monomial ideals of strong Borel type.}{14}
\contentsline {section}{\numberline {1.4}$\mathbf {d}$-fixed ideals.}{16}
\contentsline {section}{\numberline {1.5}Socle of factors by principal $\mathbf {d}$-fixed ideals.}{22}
\contentsline {section}{\numberline {1.6}A generalization of Pardue's formula.}{29}
\contentsline {section}{\numberline {1.7}$\mathbf {d}$-fixed ideals generated by powers of variables.}{31}
\contentsline {chapter}{\numberline {2}Generic initial ideal for complete intersections.}{35}
\contentsline {section}{\numberline {2.1}Main results.}{35}
\contentsline {section}{\numberline {2.2}Case $d_1+d_2\leq d_3+1$.}{37}
\contentsline {section}{\numberline {2.3}Case $d_1+d_2 > d_3+1$.}{45}
\contentsline {section}{\numberline {2.4}Generic initial ideal for $(n,d)$-complete intersections.}{64}
\contentsline {section}{\numberline {2.5}Several examples of computations of the Gin.}{69}
\contentsline {paragraph}{The case $n=4$, $d=2$.}{69}
\contentsline {paragraph}{The case $n=5$, $d=2$.}{70}
\contentsline {paragraph}{The case $n=4$, $d=3$.}{71}
\contentsline {chapter}{\numberline {}Bibliography.}{73}

\newpage
\chapter*{Introduction}

In the first chapter, we prove a stable property for monomial ideals of Borel type and give some nice consequences. Also, we discuss issues related to $\de$-fixed ideals, a new class of ideals which generalize
the class of $p$-Borel ideals. In order to explain the context, we need some preparations.

Let $K$ be an infinite field, and let $S=K[x_1,...,x_n],n\geq 2$ the polynomial ring over $K$.
Bayer and Stillman \cite{BS} note that a Borel fixed ideal $I\subset S$ satisfies the following property:
\[(*)\;\;\;\;(I:x_j^\infty)=(I:(x_1,\ldots,x_j)^\infty)\;for\; all\; j=1,\ldots,n.\] Herzog, Popescu and Vladoiu \cite{hpv} say that a monomial ideal $I$ is of \emph{Borel type}, if it satisfy $(*)$. We mention that this concept appears also in \cite{CS} as the so called \emph{weakly stable ideal}. Herzog, Popescu and Vladoiu proved in \cite{hpv} that $I$ is of Borel type, if and only if for any monomial $u\in I$ and for any $1\leq j<i \leq n$ and $q>0$ with $x_i^{q}|u$, there exists an integer $t>0$ such that $x_j^{t}u/x_i^{q}\in I$. In the first section, we present some
facts on ideals of Borel type, following \cite{hpv}. In the second section we prove that if $I$ is an ideal of Borel type, then $I_{\geq e} = $ the ideal generated by the monomials of degree $\geq e$ from $I$, is stable whenever $e\geq reg(I)$ (Theorem $1.2.10$). This allows us to give a generalization of a result of Eisenbud-Reeves-Totaro (Corollary $1.2.11$). Also, we prove that the regularity of a product of ideals of Borel type is bounded by the
sum of the regularity of those ideals (Theorem $1.2.15$). In the third section, we introduce a new class of ideals, called ideals of strong Borel type and we compute the regularity of a principal strong Borel type ideal in a special case ($1.3.6$).

A $p$-Borel ideal is a monomial ideal which satisfies certain combinatorial condition, where $p>0$ is a prime number.
It is well known that any positive integer $a$ has an unique $p$-adic decomposition $a=\sum_{i\geq 0}a_{i}p^{i}$.
 If $a,b$ are two positive integers, we write $a\leq_{p} b$ iff $a_{i}\leq b_{i}$ for any $i$,
 where $a=\sum_{i\geq 0}a_{i}p^{i}$ and $b=\sum_{i\geq 0}b_{i}p^{i}$. We say that a monomial ideal
 $I\subset S=K[x_{1},\ldots,x_{n}]$ is $p$-Borel if for any monomial $u\in I$ and for any indices $j<i$,
 if $t \leq_{p} \nu_{i}(u)$ then $x_{j}^{t}u/x_{i}^{t} \in I$, where $\nu_{i}(u) = max\{k:\; x_{i}^{k}|u\}$.

This definition suggests a natural generalization. The idea is to consider a strictly increasing sequence of positive integers $\de: 1=d_{0}|d_{1}|\cdots|d_{s}$, which we called a $\de$-sequence. Lemma $1.4.1$ states that for any positive integer $a$, there exists an unique decomposition $a=\sum_{i= 0}^{s}a_{i}d_{i}$ with $0\leq a_i < d_{i+1}/d_i$
for any $i$. If $a,b$ are two positive integers, we write $a\leq_{\de} b$ iff $a_{i}\leq b_{i}$ for any $i$, where $a=\sum_{i= 0}a_{i}d_{i}$ and $b=\sum_{i= 0}b_{i}d_{i}$. We say that a monomial ideal $I\subset S$ is $\de$-fixed if for any monomial $u\in I$ and for any indices $j<i$, if $t \leq_{\de} \nu_{i}(u)$ then $x_{j}^{t}u/x_{i}^{t} \in I$. Obviously, the $p$-Borel ideals are a special case of $\de$-fixed ideals for $\de: 1|p|p^{2}|\cdots$.

A principal $\de$-fixed ideal, is the smallest $\de$-fixed ideal which contains a given monomial. $1.4.8$ and $1.4.11$ give the explicit form of a principal $\de$-fixed ideal. In section $1.5$ we compute the socle of factors for a principal $\de$-fixed ideal ($1.5.1$ and $1.5.4$). The proofs are similar as in \cite{hpv} but we consider that is necessary to present them in this context.

In the section $1.6$ we give a formula (Theorem $1.6.1$) for the regularity of a principal
$\de$-fixed ideal $I$, which generalize the Pardue's formula for the regularity of principal $p$-Borel ideals, proved by Aramova-Herzog \cite{ah} and Herzog-Popescu \cite{hp}. Using a theorem of Popescu \cite{p} we compute the extremal Betti numbers of $S/I$ ($1.6.3$). S.Ahmad and I.Anwar proved in \cite{saf} than any ideal of Borel type has the regularity bounded by the number $q(I)=m(I)(deg(I)-1)+1$, where $m(I)=max\{i:\;x_i|u$ for some $u\in G(I)\}$ , $deg(I)=max\{deg(u):\;u\in G(I)\}$ and $G(I)$ is the set of minimal monomial generators of $I$.
As a consequence of Pardue's formula, we obtain another proof of this result in the particular case of $\de$-fixed ideals (Corollary $1.6.2$). Also, we introduce a new class of monomial ideals, called \Deisv which are sums of various \dis. Since any $\de$-fixed ideal is in particular, an ideal of Borel type, it follows that Theorem $1.2.10$ can be applied for them, so we get Corollary $1.6.6$. This result was first obtain by Herzog-Popescu \cite{hp} in the special case of a principal $p$-Borel ideal. In the last section we give the explicit form of a $\de$-fixed ideal generated by powers of variables (Proposition $1.7.2$) and make some remarks on its regularity.

In the second chapter, we discuss issues related to the generic initial ideal of an ideal generated by a regular sequence of homogeneous polynomials.

Let $S=K[x_{1},x_{2},x_{3}]$ be the polynomial ring over a field $K$ of characteristic zero. Let $f_{1}$, $f_{2}$, $f_{3}$ be a regular sequence of homogeneous polynomials of degrees $d_{1},d_{2}$ and $d_{3}$ respectively. We consider the ideal $I=(f_{1},f_{2},f_{3})\subset S$. Obviously, $S/I$ is a complete intersection Artinian $K$-algebra. One can easily check that the Hilbert series of $S/I$ depends only on the numbers $d_{1},d_{2}$ and $d_{3}$. More precisely,
\[ H(S/I,t) = (1+t+\cdots +t^{d_{1}-1})(1+t+\cdots +t^{d_{2}-1})(1+t+\cdots +t^{d_{3}-1}).\]
\cite[Lemma 2.9]{hpv} gives an explicit form of $H(S/I,t)$.

Let $S=K[x_1,\ldots,x_n]$ and let $I=(f_1,\ldots,f_n)\subset S$ be an ideal generated by a regular sequence
of homogeneous polynomials. We say that a homogeneous polynomial $f$ of degree $d$ is \emph{semiregular} for $S/I$ if the maps $(S/I)_{t}\stackrel{\cdot f}{\longrightarrow} (S/I)_{t+d}$ are either injective, either surjective for all $t\geq 0$. We say that $S/I$ has the \emph{weak Lefschetz property} (WLP) if there exists a linear form $\ell\in S$, semiregular on $S/I$, in which case we say that $\ell$ is a weak Lefschetz element for $S/I$. A theorem of Harima-Migliore-Nagel-Watanabe (see \cite{hmnw}) states that $S/I$ has (WLP) for $n=3$. We say that $S/I$ has the \emph{strong Lefschetz property} $(SLP)$ if there exists a linear form $\ell\in S$ such that $\ell^{b}$ is semiregular on $S/I$ for all integer $b\geq 1$. In this case, we say that $\ell$ is a strong Lefschetz element for $S/I$. Of course, $(SLP)\Rightarrow (WLP)$ but the converse is not true in general. In the case $n=3$, it is not known if $S/I$ has (SLP) for any regular sequence of homogeneous polynomials $f_{1}$, $f_{2}$, $f_{3}$. However, this is known true for certain cases, for example, when $f_{1},f_{2},f_{3}$ is generic, see \cite{Par} or when $f_{2}\in K[x_{2},x_{3}]$ and $f_{3}\in K[x_{3}]$, see \cite{HW} and \cite{HP}.

We say that a property $(P)$ holds for a \emph{generic} sequence of homogeneous polynomials 
$f_1,f_2,\ldots, f_n \in S=K[x_1,x_2,\ldots,x_n]$ of given degrees $d_1,d_2,\ldots,d_n$ if there exists a nonempty open Zariski subset $U\subset S_{d_1}\times S_{d_2}\times \cdots \times S_{d_n}$ such that for every $(f_{1},f_{2},\ldots,f_{n})\in U$ the property (P) holds. For example, a generic sequence of homogeneous polynomials $f_1,f_2,\ldots,f_n \in S$ is regular.

Now, we present some conjectures and the relations between them (see \cite{Par}). 

\noindent
\textbf{Conjecture A}.(Fr\"oberg) If $f_1,f_2,\ldots,f_r\in S=K[x_1,\ldots,x_n]$ is a generic sequence of homogeneous polynomials of given degrees $d_1,d_2,\ldots,d_r$ and $I=(f_1,f_2,\ldots,f_r)$ then the Hilbert series of $S/I$ is
\[ H(S/I) = \left| \frac{\prod_{i=1}^{r} (1-t^{d_i})}{(1-t)^{n}} \right|,\]
where $|\sum_{j\geq 0}a_t t^{j}| = \sum_{j\geq 0}b_t t^{j}$, with $b_{j}=a_{j}$ if $a_{i}>0$ for all $i\leq j$ and $b_{j}=0$ otherwise. \vspace{5 pt}

\noindent
\textbf{Conjecture B}. If $f_1,f_2,\ldots,f_n\in S=K[x_1,\ldots,x_n]$ is a generic sequence of homogeneous polynomials of given degrees $d_1,d_2,\ldots,d_n$ and $I=(f_1,\ldots,f_n)$ then $x_n,x_{n-1},\ldots,x_1$ is a semi-regular sequence on $A=S/I$, i.e. $x_i$ is semiregular on $A/(x_n,\ldots,x_{i+1})$ for all $1\leq i\leq n$. \vspace{5 pt}

\noindent
\textbf{Conjecture C}. If $f_1,f_2,\ldots,f_n\in S=K[x_1,\ldots,x_n]$ is a generic sequence of homogeneous polynomials of given degrees $d_1,d_2,\ldots,d_n$, $I=(f_1,\ldots,f_n)$ and $J$ is the initial ideal of $I$ with respect to the revlex order, then $x_n,x_{n-1},\ldots,x_1$ is a semi-regular sequence on $A=S/(f_1,\ldots,f_n)$. \vspace{5 pt}

\noindent
\textbf{Conjecture D}.(Moreno) If $f_1,f_2,\ldots,f_n\in S=K[x_1,\ldots,x_n]$ is a generic sequence of homogeneous polynomials of given degrees $d_1,d_2,\ldots,d_n$, $I=(f_1,\ldots,f_n)$ and $J$ is the initial ideal of $I$ with respect to the revlex order, then $J$ is an almost revlex ideal, i.e. if $u\in J$ is a minimal generator of $J$ then every monomial of the same degree which preceeds $u$ must be in $J$ as well. \vspace{5 pt}
\vspace{5 pt}

Pardue proved in \cite{Par} that if conjecture $A$ is true for some positive integer $n$ then the conjecture $B$ is true for the same $n$. Also, conjecture $C$ is true for $n$ if and only if $B$ is true for $n$ and if conjecture $B$ is true for some $r$ then $A$ is true for $n<r$ and exactly for that $r$. Also, if conjecture $D$ is true for some $n$ then $B$, and thus $C$, are true for the same $n$. Fr\"oberg \cite{F} and Anick \cite{A} proved that $A$ is true for $n\leq 3$ and so $B$ and $C$ are true for $n\leq 3$. Moreno \cite{Mor} remarked that $D$ is true for $n=2$. Note that Conjecture $A$ for $n=3$ does not imply the Moreno's conjecture $D$ for $n=3$.

Let $I\subset S=K[x_1,\ldots,x_n]$ be a graded ideal, i.e. an ideal generated by homogeneous polynomials. We choose a monomial order "$\leq$" on the set of monomials of $S$. 
If $\alpha=(\alpha_{ij})$ is a $n\times n$ invertible matrix with entries in $K$ and $f\in I$ is a polynomial, we denote by $\alpha f$ the polynomial obtained from $f$ by the changing of variables, $x_i \mapsto \sum_{j=1}^{n}\alpha_{ij}x_j$ for all $i=1,\ldots,n$. We denote $\alpha I = (\alpha f|\;f\in I)$. Galligo and Bayer-Stillman proved that there exists a nonempty open Zariski subset $U\in GL_n(K)$, such that for any $\alpha,\alpha'\in U$, $in_{\leq}(\alpha I) = in_{\leq}(\alpha' I)$. For an 
$\alpha\in U$, we denote $in_{\leq}(\alpha I):=gin_{\leq}(I)$ and we called it the generic initial ideal of $I$,
with respect to "$\leq$". For an introduction on generic initial ideal, see \cite[\S 15.9]{E}. The generic initial ideal is Borel fixed. In the case of characteristic zero, that means that it is strongly stable and in the case of positive characteristic $p$ that means it is $p$-Borel. This remark shows a connection between the two parts of my thesis.

Now, let $I=(f_{1},f_{2},f_{3})\subset S=K[x_1,x_2,x_3]$ an ideal generated by a regular sequence of homogeneous
polynomials $f_{1},f_{2},f_{3}$ of degrees $d_1$,$d_2$ and $d_3$, respectively. Let $J=Gin(I)$ be the generic initial ideal of $I$, with respect to the reverse lexicographic order. Our aim is to compute $J$ for all regular sequences $f_1,f_2,f_3$ of homogeneous polynomials of given degrees $d_1,d_2,d_3$ such that $S/I$ has (SLP). We will do this in the sections $2.2$ and $2.3$. These computations shown us in particular, that $J$ depends only on the numbers $d_1,d_2,d_3$ (this has been proved also by Popescu and Vladoiu in \cite{PV}) and more important, that $J$ is an almost reverse lexicographic ideal (Theorem $2.1.1$). As a consequence, conjecture Moreno ($D$) is true for $n=3$ and $char(K)=0$ (Theorem $2.1.2$). 

Now, let $K$ be an algebraically closed field of characteristic zero. Let $S=K[x_{1},\ldots,x_{n}]$ be the polynomial ring in $n$ variables over $K$. Let $n,d\geq 2$ be two integers. We consider \linebreak $I=(f_{1},\ldots,f_{n})\subset S$ an ideal generated by a regular sequence $f_{1},\ldots,f_{n}\in S$ of homogeneous polynomials of degree $d$. We say that $A=S/I$ is a \emph{$(n,d)$-complete intersection}. Let $J=Gin(I)$ be the generic initial of $I$, with respect to the reverse lexicographic (revlex) order. 
With the above notations, Conca and Sidman \cite{CS} proved that $J_d =$ the set of monomials of degree $d$ of $J$, is revlex if $f_{1},\ldots,f_{n}$ is a generic regular sequence, (see \cite[Theorem 1.2]{CS}). 

In the section $2.4$, we prove that $J_{d}$ is a revlex set in another case, namely, when $f_{i}\in K[x_{i},\ldots,x_{n}]$ for all $1\leq i\leq n$. It is likely to be true that $J_{d}$ is revlex for any $(n,d)$-complete intersection, but we do not have the means to prove this assertion. As Example $2.4.10$ shows, the hypotheses $char(K)=0$ and $f_1,\ldots,f_n$ is a regular sequence are essential.

In the section $2.5$, we compute the generic initial ideal for some particular cases of $(n,d)$-complete intersections: $(n=4,d=2)$, $(n=5,d=2)$ and $(n=4,d=3)$. In order to do this, we suppose in addition that $S/I$ has 
(SLP). Note that this property holds for generic complete intersection (see \cite{Par}) and for the case when 
$f_{i}\in k[x_{i},\ldots,x_{n}]$. Also, it was conjectured that (SLP) holds for any standard complete intersection. By a theorem of Wiebe \cite{W}, $S/I$ has (WLP) (respectively (SLP)) if and only if $x_n$ is a weak (respectively strong) Lefschetz element for $S/J$. This result is very important for our computations.

In the writing on this thesis, we used new results from our articles and preprints. In the section $1.2$ we followed \cite{mir4} and \cite{mir7}. In the sections $1.4$, $1.5$ and $1.6$ we followed \cite{mir} and \cite{mir2}. In the sections $1.3$ and $1.7$ we used \cite{mir6}. In the sections $2.1$, $2.2$ and $2.3$ we followed \cite{mir3} and in the sections $2.4$ and $2.5$ we followed \cite{mir5}, respectively.

\chapter{Ideals of Borel type and $\mathbf {d}$-fixed ideals.}
\section{Ideals of Borel type.}

Let $K$ be a field and $S=K[x_1,\ldots,x_n]$ the ring of polynomials over $K$.
Herzog, Popescu and Vl\u adoiu introduced in \cite{hpv} the following definition.

\begin{dfn}
A monomial ideal $I\subset k[x_{1},\ldots,x_{n}]$ is said to be of Borel type if
\[ (I:x_{j}^{\infty}) = (I:(x_{1},\ldots,x_{j})^{\infty}),\;for\;any\;j=1,\ldots,n.\]
\end{dfn}

We have the following equivalent characterization of ideals of Borel type.

\begin{prop}\cite[Proposition 2.2]{hpv}
Let $I\subset S$ be a monomial ideal. The following conditions are equivalent:

(a) $I$ is an ideal of Borel type.

(b) For any $1\leq j<i\leq n$, we have $(I:x_i^{\infty}) \subset (I:x_j^{\infty})$;

(c) Let $u\in I$ be a monomial and suppose that $x_i^q|u$ for some $q>0$. Then for any $j<i$ there exists
    an integer $t$ such that $x_j^tu/x_i^{q}\in I$;

(d) Let $u\in I$ be a monomial; then for any $1\leq j<i\leq n$, there exists an integer $t>0$ such that
    $x_j^{t}u/x_i^{\nu_i(u)} \in I$.

Moreover, it is easy to see that the conditions (c) and (d) are satisfied if and only if they are satisfied for
all $u\in G(I)$.
\end{prop}

\begin{proof}
$(a)\Rightarrow (b)$ is trivial. For the converse, we use induction on $1\leq j\leq n$, the assertion 
being obvious for $j=1$. Suppose $j<n$ and $(I:x_j^{\infty}) = (I:(x_1,\ldots,x_{j})^{\infty})$. Since by (b),
$(I:x_{j+1}^{\infty})\subset (I:x_j^{\infty})$ it follows that $(I:x_{j+1}^{\infty})\subset (I:(x_1,\ldots,x_{j})^{\infty})$ and thus $(I:x_{j+1}^{\infty})\subset (I:(x_1,\ldots,x_{j+1})^{\infty})$. Since
the converse inclusion it is obvious, we get the required conclusion.

$(c)\Rightarrow (d)$ is trivial. For the converse, let $u\in I$ be a monomial such that $x_i^q|u$ for some $q>0$
and let $j<i$. By (d) there exists $t$ such that $x_j^{t}u/x_i^{\nu_i(u)} \in I$. Therefore $x_j^{t}u/x_i^q = x_i^{\nu_i(u)-q}x_j^{t}u/x_i^{\nu_i(u)} \in I$.

$(b)\Rightarrow (c)$: Let $u\in I$ be a monomial such that $x_i^q|u$ for some $q>0$ and let $j<i$. Then $u=x_i^qv$
with $v\in (I:x_i^{\infty})$. Therefore, there exists $t$ such that $x_j^tu/x_i^q = x_j^tv\in I$.

$(c)\Rightarrow (b)$: Let $u\in (I:x_i^{\infty})$ be a monomial. Then $x_i^qu\in I$ for some $q>0$ and so $(c)$ implies that $x_j^tu\in I$ for some $t$, that is, $u\in (I:x_j^{\infty})$.
\end{proof}

\begin{prop}
(a) If $I,J\subset S$ are two ideals of Borel type then $I+J$, $I\cap J$ and $I\cdot J$ are also ideals of Borel type.

(b) If $I\subset S$ is an ideal of Borel type and $J\subset S$ is an arbitrary monomial ideal, then $(I:J)$ is an ideal
    of Borel type.
\end{prop}

\begin{proof}
(a) Since a monomial of $I+J$ is either in $I$, either in $J$ it follows immediately that $I+J$ is of Borel type, using
    the characterization $(d)$ from the previous proposition. A similar argument holds for $I\cap J$. Now, let $u\in
    I\cdot J$ be a monomial. It follows that $u=v\cdot w$, where $v\in I$ and $w\in J$ are monomials. Let $1\leq i\leq
    n$ such that $x_i|u$ and let $1\leq j<i$. Since $I$ is of Borel type, then there exists some $t_1\geq 0$
    such that $x_j^{t_1}\cdot v/x_i^{\nu_i(v)}\in I$. Analogously, there exists some $t_2\geq 0$ such that
    $x_j^{t_2}\cdot w/x_i^{\nu_i(w)}\in J$. It follows that $x_j^{t_1+t_2}u/x_i^{\nu_i(u)}\in I\cdot J$, therefore
    $I\cdot J$ is of Borel type.

(b) Suppose $J=(v_1,\ldots,v_m)$, where $v_i\in S$ are monomials. Since $(I:J)=\cap_{k=1}^{m} (I:v_k)$ and the
    intersection of Borel type ideals is still of Borel type, we can assume $m=1$. Denote $v_1:=v$. Let $u\in (I:v)$
    be a monomial. We have $u\cdot v\in I$. Let $1\leq i\leq n$ such that $x_i|u$ and let $1\leq j<i$. Since $I$ is of
    Borel type, there exists some $t\geq 0$ such that $x_j^t u\cdot v/x_i^{\nu_i(uv)}\in I$. In particular, multiplying
    by $x_i^{\nu_i(v)}$, it follows that $v\cdot(x_j^t u/x_i^{\nu_i(u)})\in I$ and thus $x_j^t u/x_i^{\nu_i(u)}\in
    (I:v)$. In conclusion, $(I:v)$ is of Borel type, as required.
\end{proof}

We recall the following definition of Stanley, see \cite{S}.

\begin{dfn}
Let $S=k[x_{1},\ldots,x_{n}]$ and let $M$ be a finitely generated graded $S$-module. The module $M$ is sequentially Cohen-Macaulay if there exists a finite filtration $0 = M_{0} \subset M_{1} \subset \cdots \subset M_{r} = M$ of $M$
by graded submodules of $M$ such that:
\begin{itemize}
    \item $M_{i}/M_{i-1}$ are Cohen-Macaulay for any $i=1,\ldots,r$ and
    \item $dim(M_{1}/M_{0})< dim(M_{2}/M_{1})< \cdots < dim(M_{r}/M_{r-1})$.
\end{itemize}
The above filtration is unique and is called the CM-filtration of $M$

In particular, if $I\subset S$ is a graded ideal then $R=S/I$ is sequentially Cohen-Macaulay if there exists a chain of ideals $I = I_{0} \subset I_{1} \subset \cdots \subset I_{r} = S$ such that $I_{j}/I_{j-1}$ are Cohen-Macaulay and $dim(I_{j}/I_{j-1})<dim(I_{j+1}/I_{j})$ for any $j=1,\ldots,r-1$.
\end{dfn}

Let $I\subset S$ be a monomial ideal. Recursively we define an ascending chain of monomial ideals as follows: We let $I_{0}:=I$. Suppose $I_{\ell}$ is already defined. If $I_{\ell}=S$ then the chain ends. Otherwise, let $n_{\ell} = max\{i:\;x_{i}|u$ for an $u\in G(I_{\ell}) \}$. We set $I_{\ell +1}:= (I_{\ell}:x_{n_{\ell}}^{\infty})$. It is obvious that $n_{\ell}>n_{\ell+1}$, and therefore the chain $I_{0}\subset I_{1} \subset \cdots \subset I_{r}=S$ is finite and has length $r\leq n$. We call this chain of ideals, the \emph{sequential chain} of $I$. Note that if $I$ is an Artinian monomial ideal then $r=1$. The converse is true for ideals of Borel type.
\vspace{5 pt}

\begin{prop}\cite[Corollary 2.5]{hpv}
Let $I$ be a monomial ideal of Borel type. Then $R=S/I$ is sequentially Cohen-Macaulay.
\end{prop}

\begin{proof}
We may assume $I\neq 0$. Let $I=I_0\subset I_1\subset \cdots \subset I_r=S$ be the sequential chain of $I$. Note that, inductively, we get that any ideal $I_{\ell}$ is an ideal of Borel type, since $I_{\ell+1}$ is a quotient of $I_{\ell}$.
In particular, $I_{\ell+1}=(I_{\ell}:(x_1,\ldots,x_{n_{\ell}})^{\infty})$ for all $\ell$. Fix an integer $\ell<r$. Let
$n_j=m(I_j)$ for all j, then the elements of $G(I_j)$ belong to $K[x_1,\ldots,x_{n_{\ell}}]$ for all $j\geq\ell$. Let
$J_{\ell}$ be the ideal generated by $G(I_{\ell})$ in $K[x_1,\ldots,x_{n_{\ell}}]$. Then the saturation
$J_{\ell}^{sat}=(J_{\ell}:(x_1,\ldots,x_{n_{\ell}})^{\infty})$ is generated by the elements of $G(I_{\ell+1})$. It follows that
\[I_{\ell+1}/I_{\ell}
\cong (J_{\ell}^{sat}/J_{\ell})[x_{n_{\ell}+1}, \ldots, x_{n}]\]
is an $(n-n_{\ell})-$ dimensional Cohen-Macaulay $S$-module.
\end{proof}

Let $M$ be a finitely generated graded $S$-module with the minimal graded free resolution $0 \rightarrow F_{s} \rightarrow F_{s-1} \rightarrow \cdots \rightarrow F_{0} \rightarrow M \rightarrow 0$. Let $Syz_{t}(M)=Ker(F_{t} \rightarrow F_{t-1})$. The module $M$ is called $(r,t)$-regular if $Syz_{t}(M)$ is $(r+t)$-regular in the sense that all generators of $F_{j}$ for $t\leq j\leq s$ have degrees $\leq j+r$. The $t$-regularity of $M$ is by definition $(t-reg)(M) = min\{r|\; M\;is\;(r,t)-regular \}$. Obvious $(t-reg)(M)\leq ((t-1)-reg)(M)$. If the equality is strict and $r=(t-reg)(M)$ then $(r,t)$ is called a corner of $M$ and $\beta_{t,r+t}(M)$ is an extremal Betti number of $M$, where $\beta_{ij} = dim_{k}Tor_{i}(k,M)_{j}$ denotes the $ij$-th graded Betti number of $M$. Later, we will use the following result of Popescu:

\begin{teor}\cite[Theorem 3.2]{p}
If $I\subset S$ is a Borel type ideal, then $S/I$ has at most $r+1$-corners among $(n_{\ell},s(J_{\ell}^{sat}/J_{\ell}))$
and the corresponding extremal Betti numbers are
\[ \beta_{n_{\ell},s(J_{\ell}^{sat}/J_{\ell}) + n_{\ell}} (S/I) = dim_{k}(J_{\ell}^{sat}/J_{\ell})_{s(J_{\ell}^{sat}/J_{\ell})}.\]
\end{teor}

\section{Stable properties of Borel type ideals.}

It would be appropriate to recall the definition of the Castelnuovo-Mumford regularity. We refer the reader to \cite{E} for further details on the subject.

\begin{dfn}
Let $K$ be an infinite field, and let $S=K[x_1,...,x_n],n\geq 2$ the polynomial ring over $K$.
Let $M$ be a finitely generated graded $S$-module. The \emph{Castelnuovo-Mumford regularity} $reg(M)$ of $M$ is
\[ \max_{i,j} \{j-i :\; \beta_{ij}(M)\neq 0\}.\]
\end{dfn}

If $M$ is an artinian $S$-module, we denote $s(M)=max\{t:\; M_t\neq 0\}$. Herzog, Popescu and Vl\u adoiu proved the following formula for the regularity of an ideal of Borel type.

\begin{prop}\cite[Corollary 2.7]{hpv}
If $I$ is a Borel type ideal, with the notations of the section $1.1$, we have
\[ reg(I) = max\{s(J_0^{sat}/J_0), \ldots , s(J_{r-1}^{sat}/J_{r-1})\} + 1.\]
\end{prop}

If $I\subset S$ is a monomial ideal, we denote $q(I)=m(I)(deg(I)-1)+1$, where $m(I)$ is the maximal index of a variable
which appear in a monomial from $G(I)$ and $deg(I)$ is the maximal degree of a
monomial from $G(I)$. We cite the following characterization of ideals of Borel type, given
by S.Ahmad and I.Anwar.

\begin{teor}\cite[Theorem 2.2]{saf}
Let $I\subset S$ be a monomial ideal. Then the following statements are equivalent:
\begin{enumerate}
	\item $I$ is an ideal of Borel type.
	\item Each $P\in Ass(S/I)$ has the form $P=(x_1,x_2,...,x_r)$ for some $1\leq r\leq n$.
	\item $I_{\geq q(I)}$ is stable.
\end{enumerate}
\end{teor}

\begin{obs}
{\em Note that the implication $(1)\Rightarrow (2)$ follows immediately from $1.1.3(b)$, since any prime $P\in Ass(S/I)$ can be written as $P=(I:u)$ for some monomial $u\in S$. 
Another proof is given in \cite[Proposition 5.2]{hp-f}.}
\end{obs}

We recall the following result of Eisenbud-Reeves-Totaro.

\begin{prop}\cite[Proposition 12]{ert}
Let $I$ be a monomial ideal with $deg(I)=d$ and let $e\geq d$ such that $I_{\geq e}$ is stable. Then $reg(I)\leq e$.
\end{prop}

\begin{cor}\cite[Corollary 2.4]{saf}
If $I$ is of Borel type then $reg(I)\leq m(I)(deg(I)-1)+1$. 
The same holds for a monomial ideal $I$ with $Ass(S/I)$ totally
ordered by inclusion.
\end{cor}

\begin{proof}By the Theorem $1.2.3$ we have $I_{\geq q(I)}$
stable. As $q(I)\geq deg(I)$ we get $reg(I)\leq q(I)$ by Proposition $1.2.5$.
For the second statement, we renumber the variables $x_i$'s such that $I$
satisfies $(2)$ from Theorem $1.2.3$ and than apply the first statement.
\end{proof}

\begin{obs}
{\em The number $m(I) \cdot (deg(I)-1)+1$ is the best possible linear upper bound for the regularity of a Borel type ideal, $I$. Indeed, if we consider $I=(x_1^2,x_2^2)\subset K[x_1,x_2]$, we have $reg(I)=3 = 2\cdot(deg(I)-1)+1$. See also, \cite[Remark 1.3]{saf}.}
\end{obs}

\begin{lema}
Let $I\subset S$ be a monomial ideal and $I'=IS'$ the extension of $I$ in $S'=S[x_{n+1}]$. If $e\geq deg(I)$, then
$I_{\geq e}$ is stable if and only if $I'_{\geq e}$ is stable.
\end{lema}

\begin{proof}
Let $u\in I'_{\geq e}$ be a monomial. Then $u=x_{n+1}^{k}\cdot v$ for some $v\in I$. If $k>0$ then $m(u)=n+1$ and therefore, for any $i<n+1$, $x_i\cdot u/x_{n+1} = x_{n+1}^{k-1}\cdot x_i\cdot v \in I'_{\geq e}$. If $k=0$ then
$m(u)\leq n$ and since $I_{\geq e}$ is stable, it follows $x_i \cdot u/ x_{m(u)} \in I'_{\geq e}$. Thus $I'_{\geq e}$
is stable. For the converse, simply notice that $G(I_{\geq e})\subset G(I'_{\geq e})$ and since is enough to
check the stable property only for the minimal generators, we are done.
\end{proof}

\begin{lema}
If $I\subset S$ is an Artinian monomial ideal and $e\geq reg(I)$ then $I_{\geq e}$ is stable.
\end{lema}

\begin{proof}
Since $I$ is Artinian, it follows that the length of the sequential chain of $I$ is $r=1$. By $1.2.2$ we get
$reg(I) = s(S/I) + 1$ and therefore, if $e\geq reg(I)$ then $I_{\geq e} = S_{\geq e}$, thus $I_{\geq e}$ is stable.
\end{proof}

\begin{teor}
Let $I\subset S$ be a Borel type ideal and let $e\geq reg(I)$ be an integer. Then $I_{\geq e}$ is stable.
\end{teor}

\begin{proof}
We use induction on $r\geq 1$, where $r$ is the length of the sequential chain of $I$. If $r=1$, i.e. $I$ is an
Artinian ideal, we are done as in the proof of the previous lemma. 

Suppose now $r>1$ and let $I = I_0 \subset I_1 \subset \cdots \subset I_r = S$ be the sequential chain of $I$. Let $n_{\ell}=m(I_{\ell})$ for  $0\leq \ell \leq r$. Let $J_{\ell}\subset S_{\ell} = K[x_{1},\ldots,x_{n_{\ell}}]$ be the ideal generated by $G(I_{\ell})$. Using the induction hypothesis, we may assume $(I_1)_{\geq e}$ stable for $e\geq reg(I_1)$. On the other hand, from $1.2.2$ it follows that $reg(I_1)\leq reg(I)$, thus 
$(I_1)_{\geq e}$ is stable for $e\geq reg(I)$.

Since $J_0^{sat} = I_1 \cap S_{n_0}$, using iteratively Lemma $1.2.8$ it follows that $(J_0^{sat})_{\geq e}$ is stable.
Let $e\geq reg(I)$. Since $reg(I) \geq s(J_0^{sat}/J_0) + 1$ it follows
that $(J_0)_{\geq e} = (J_0^{sat})_{\geq e}$ is stable. Since $I_0 = J_0S$, using again Lemma $1.2.8$, we get
$I_{\geq e}$ stable for $e\geq reg(I)$, as required.
\end{proof}

Theorem $1.2.10$ and Proposition $1.2.5$ yield the following:

\begin{cor}
If $I$ is a Borel type ideal, then $reg(I)=min\{e:\;e\geq deg(I)$ and $I_{\geq e}\;is\; stable \}$. The same conclusion holds, if $I$ is a monomial ideal with $Ass(S/I)$ totally ordered by inclusion.
\end{cor}

\begin{proof}
Denote $f=min\{e:\; I_{\geq e}\;is\; stable \}$. By $1.2.10$, we get $reg(I)\geq f$ and by $1.2.5$, $reg(I)\leq f$.
For the second statement, by renumbering the variables, we can assume that $I$ is of Borel type.
\end{proof}

\begin{exm}{\em
Let $I = (x_1^7,\; x_1^5x_2,\; x_1^2x_2^4,\; x_1x_2^6,\; x_1^5x_3^2,\; x_1x_2^4x_3^2)\subset K[x_1,x_2,x_3,x_4]$.
We construct the sequential chain of $I$. We have $I_0 = I$ and $n_0=m(I_0)=3$, therefore $J_0 = I_0 \cap K[x_1,x_2,x_3]$. Let $I_1 = (I_0 : x_3^{\infty}) = (x_1^5,x_1x_2^4)$. We have $n_1 = m(I_1) = 2$, therefore
$J_1 = I_1 \cap K[x_1,x_2]$. Let $I_2 = (I_1 : x_2^{\infty}) = (x_1)$. We have $n_2 = m(I_2) = 1$, therefore
$J_2 = I_2 \cap K[x_1]$. One can easily compute, $s(J_0^{sat}/J_0) = 7$, $s(J_1^{sat}/J_1)=7$ and $s(J_2^{sat}/J_2)=1$.
Using $1.2.2$, we get:
\[ reg(I) = max\{s(J_0^{sat}/J_0), s(J_1^{sat}/J_1), s(J_2^{sat}/J_2)\} + 1 = 8. \]
We will exemplify the proof of $1.2.10$ for $I$. Let $e\geq 8$ be an integer. Since $(J_2)_{\geq e}$ is obviously stable, it follows from $1.2.8$ that $(I_2)_{\geq e} = (J_2S)_{\geq e}$ is also stable. Note that $I_2 = J_1^{sat}S$ and moreover, $I_2$ and $J_1^{sat}$ have the same minimal set of generators. Therefore, by $1.2.8$, it follows that $(J_1^{sat})_{\geq e}$ is stable. Since $e\geq reg(I)> s(J_1^{sat}/J_1) = 7$, it follows that $(J_1)_{\geq e}$ is stable.
On the other hand, $I_1=J_1S$, therefore $(I_1)_{\geq e}$ is stable. Since $J_0^{sat}S = I_1$ we get, from $1.2.8$,
$(J_0^{sat})_{\geq}$ stable. Since $e\geq reg(I)> s(J_0^{sat}/J_0) = 7$, it follows that $(J_0)_{\geq e}$ is stable.
Finally, since $I = I_0 = J_0S$, we obtain $I_{\geq e}$ stable, as required.
}\end{exm}

\begin{cor}
If $I$ and $J$ are ideals of Borel type, then

(a) $reg(I+J)\leq max\{reg(I),reg(J)\}$;

(b) $reg(I\cap J)\leq max\{reg(I),reg(J)\}$.
\end{cor}

\begin{proof}
Denote $e=max\{reg(I),reg(J)\}$. From Theorem $1.2.10$, it follows that $I_{\geq e}$ and $J_{\geq e}$ stable. Therefore,
$(I+J)_{\geq e} = I_{\geq e} + J_{\geq e}$ is stable, as a sum of stable ideals. By $1.2.5$ it follows that $reg(I+J)\leq e$, thus (a) holds. The proof of (b) is similar.
\end{proof}

\begin{lema}
Let $I,J\subset S$ be two monomial ideals and let $e\geq deg(I)$ and $f\geq deg(J)$ be two integers such that 
$I_{\geq e}$ and $J_{\geq f}$ are stable. Then $(I\cdot J)_{\geq e+f}$ is stable.
\end{lema}

\begin{proof}
Let $u\in (I\cdot J)_{\geq e+f}$. It follows that $u=v\cdot w$ for some monomials $v\in I$ and $w\in J$.
We claim that we can choose $v$ and $w$ such that $v\in I_{\geq e}$ and $w\in J_{\geq e}$. 

Indeed, if we write $v=v'a$ for $v'\in G(I)$ and $a\in S$ a monomial and $w=w'b$ for $w'\in G(J)$ and $b\in S$ a monomial, we can find some new monomials $\bar{a}, \bar{b}\in S$ such that $\bar{a}\bar{b}=ab$, $deg(v'\bar{a})\geq e$ and $deg(w'\bar{b})\geq f$. We are able to do this, since $e\geq deg(I)$ and $f\geq deg(J)$, so $deg(v')\leq e$ and $deg(w')\leq f$. Changing $v$ with $v'\bar{a}$ and $w$ with $w'\bar{b}$, the claim is proved.

Now, let $j<m(u)$ be an integer and suppose $m(u)=m(v)$. Then $x_j u/x_{m(u)} = (x_j v/x_{m(v)})\cdot w \in I\cdot J$,
because $x_j v/x_{m(v)}\in I$ since $I_{\geq e}$ is stable. Analogously, if $m(u)=m(w)$, then
$x_j u/x_{m(u)} = v\cdot (x_j w/x_{m(w)}) \in I\cdot J$. Therefore, $(I\cdot J)_{\geq e+f}$ is stable.
\end{proof}

\begin{teor}
Let $I,J\subset S$ be two monomial ideals of Borel type. Then
\[ reg(I\cdot J) \leq reg(I) + reg(J).\]
\end{teor}

\begin{proof}
Since $I$ and $J$ are ideals of Borel type, if we denote $e:=reg(I)$ and $f=reg(J)$, by Theorem $1.2.10$ it follows that
$I_{\geq e}$ and $J_{\geq f}$ are stable. Using the previous lemma, it follows that $(I\cdot J)_{\geq e+f}$ is stable,
therefore, using Proposition $1.2.5$ we get $reg(I\cdot J)\leq e+f$ as required.
\end{proof}

\begin{cor}
If $I\subset S$ is an ideal of Borel type, then $reg(I^k)\leq k\cdot reg(I)$.
\end{cor}

Note that there are other large classes of graded ideals which have this property, see for instance \cite{CH}, but on the other hand, Sturmfels provided an example of a graded ideal $I\subset S$ with $reg(I^2)>2\cdot reg(I)$ in \cite{S}.

\section{Monomial ideals of strong Borel type.}

Proposition $1.1.2(d)$ suggested us the following definition.

\begin{dfn}
We say that a monomial ideal $I$ is of {\em strong Borel type} (SBT) if for any monomial $u\in I$ and for any $1\leq j<i \leq n$, there exists an integer $t\leq \nu_i(u)$ such that $x_j^{t}u/x_i^{\nu_i(u)}\in I$.
\end{dfn}

\begin{obs}{\em
Obviously, an ideal of strong Borel type is also an ideal of Borel type, but the converse is not true. Take for instance
$I=(x_1^3,x_2^2)\subset K[x_1,x_2]$.

The sum of two ideals of (SBT) is still an ideal of (SBT). Also, the same is true for
an intersection or a product of two ideals of (SBT). The proof is similar with the proof of $1.1.3$, so we skip it.
}\end{obs}

\begin{dfn}
Let $\mathcal A\subset S$ be a set of monomials. We say that $I$ is the {\em (SBT)-ideal} generated by $\mathcal A$, if $I$ is the smallest, with respect to inclusion, ideal of (SBT) containing $\mathcal A$. We write $I=SBT(\mathcal A)$.
In particular, if $\mathcal A=\{u\}$, where $u\in S$ is a monomial, we say that $I$ is the {\em principal (SBT)-ideal}
generated by $u$, and we write $I=SBT(u)$.
\end{dfn}

\begin{lema}
Let $1\leq i_1<i_2<\cdots <i_r\leq n$ be some integers, $\alpha_1,\ldots,\alpha_r$ some positive integers and
$u=x_{i_1}^{\alpha_1}x_{i_2}^{\alpha_2}\cdots x_{i_r}^{\alpha_r} \in S$. Then, the principal (SBT)-ideal generated by $u$, is:
\[ I = SBT(u) = \prod_{q=1}^{r}(\me_q^{[\alpha_q]}),\;where\;\me_q=\{x_1,\ldots,x_{i_q}\}\;and\; 
\me_q^{[\alpha_q]}=\{x_1^{\alpha_q},\ldots,x_{i_q}^{\alpha_q}\}. \]
\end{lema}

\begin{proof}
Denote $I'=\prod_{q=1}^{r}(\me_q^{[\alpha_q]})$. If $v\in G(I')$, then $v=x_{j_1}^{\alpha_1}x_{j_2}^{\alpha_2}\cdots
x_{j_r}^{\alpha_r}$, for some $1\leq j_q\leq i_q$, where $1\leq q\leq r$. Since
\[ v = \frac{x_{j_r}^{\alpha_r}}{x_{i_r}^{\alpha_r}}\cdots \frac{x_{j_2}^{\alpha_2}}{x_{i_2}^{\alpha_2}}\cdot
\frac{x_{j_1}^{\alpha_1}}{x_{i_1}^{\alpha_1}} u, \]
and $I$ is of (SBT) it follows that $v\in I$ and thus $I' \subseteq I$. For the converse, simply notice that $I'$ is itself a (SBT)-ideal. Therefore $I=SBT(u)$ as required.
\end{proof}

\begin{obs}
Let $1\leq i_1<i_2<\cdots <i_r\leq n$ be some integers and $\alpha_1,\ldots,\alpha_r$ some positive integers and
$u=x_{i_1}^{\alpha_1}x_{i_2}^{\alpha_2}\cdots x_{i_r}^{\alpha_r} \in S$. Let $I=SBT(u)$. We describe the sequential chain of $I$. Denote $I_r:=I$. Since $I_r = I = \prod_{q=1}^{r}(\me_q^{[\alpha_q]})$, it follows that $I_{r-1} := (I_r:x_{i_r}^{\infty}) = \prod_{q=1}^{r-1}(\me_q^{[\alpha_q]})$. Analogously, we get $I_{q}=(I_{q+1}:x_{i_{q+1}}^{\infty}) = \prod_{e=1}^{q}(\me_e^{[\alpha_e]})$, for all $0\leq q < r$. Therefore, the sequential chain of $I$ is, \[ I = I_r \subset I_{r-1} \subset \cdots \subset I_1 \subset I_0 = S. \]
Let $J_q$ be the ideal in $S_q = K[x_1,\ldots,x_{i_q}]$ generated by $G(I_q)$, for $1\leq q\leq r$. \linebreak 
If $s_q = s(J_q^{sat}/J_q)$, $1.2.2$ implies $reg(I)=max\{s_q:\; 1\leq q\leq r\}$
\end{obs}

Our next goal is to compute the regularity of a principal (SBT)-ideal, in a special case. In order to do so, we will compute the
sequential chain of $I$ and than apply Proposition $1.2.2$. 

\begin{teor}
Let $1\leq i_1<i_2<\cdots <i_r\leq n$ be some integers, $\alpha_1\geq \alpha_1 \geq \cdots \geq \alpha_r$ some positive integers and
$u=x_{i_1}^{\alpha_1}x_{i_2}^{\alpha_2}\cdots x_{i_r}^{\alpha_r} \in S$. Let $I=SBT(u)$. For $1\leq q\leq r$, we define the numbers:
\[\chi_q = \alpha_1 + \cdots + \alpha_{q-1} + (\alpha_q-1)i_q. \]
With the above notations, we have $reg(SBT(u))= \max\limits_{q=1}^{r}\chi_q + 1$.
\end{teor}

\begin{proof}
We use the notations from $1.3.5$. In order to compute the regularity of $I$ we must determine the numbers $s_q:=s(J_q^{sat}/J_q)$.
We will prove that $s_q = \chi_q$. First of all, note that $J_q = I_q\cap S_q$ and $J_q^{sat} = I_{q-1}\cap S_q$.
We fix $1\leq q\leq r$ and we denote $A:=\{1,2,\ldots,i_q\}\setminus \{i_1,\ldots,i_{q-1}\}$. Let \[w=\prod_{e=1}^{q-1} x_{i_e}^{\alpha_e + \alpha_q -1} \cdot \prod_{j\in A}x_j^{\alpha_q-1} \in S_q. \]
Obviously, $deg(w)=\chi_q$. We claim that $w\in J_q^{sat}$, but $w\notin J_q$ and, therefore, $s_q\geq \chi_q$. Indeed, since $\prod_{e=1}^{q-1} x_{i_e}^{\alpha_e + \alpha_q -1}\in J_q^{sat} = \prod_{e=1}^{q-1}(\me_e^{[\alpha_e]})\subset S_q$ it follows that $w\in J_q^{sat}$. 

We assume, by contradiction,
that $w\in J_q$. Thus, $w=x_{j_1}^{\alpha_1}\cdots x_{j_q}^{\alpha_q}y$, where $1\leq j_e\leq i_e$ for all $e\in \{1,\ldots,q\}$ and $y\in S$ is a monomial. We claim that $\{j_1,j_2,\cdots,j_{q-1}\}=\{i_1,\ldots,i_{q-1}\}$. Let $k\in \{1,\ldots,q-1\}$. If $j_k\in A$, since $x_{j_k}^{\alpha_k}|w$ it follows that $\alpha_k\leq \alpha_q -1$, a contradiction, since $\alpha_k\geq \alpha_q$. Therefore, $j_k\in\{i_1,\ldots,i_{q-1}\}$. 

Using induction on $1\leq e\leq q-1$, we prove that $\{j_1,\ldots,j_e\}=\{i_1,\ldots,i_e\}$ for all $e$. Let $e=1$. Since $x_{j_1}^{\alpha_1}|w$, $j_1\leq i_1$ and for all $j<i_1$ we have $j\in A$, it follows that $j_1=i_1$. Suppose $e\leq q-1$ and $\{j_1,\ldots,j_{e-1}\}=\{i_1,\ldots,i_{e-1}\}$. We have
$j_e\leq i_e$. Suppose $j_e<i_e$. Since $x_{j_e}^{\alpha_e}|w$ it follows that $j_e=i_k$ for some $k<e$. So
$w = (x_{i_1}^{\alpha_1}\cdots x_{i_k}^{\alpha_k}\cdots x_{i_{e-1}}^{\alpha_{e-1}})(x_{j_e}^{\alpha_e}\cdots x_{j_q}^{\alpha_q}\cdot y)$
and it follows $x_{i_k}^{\alpha_e+\alpha_k}|w$. But this is false, since $\alpha_e+\alpha_k>\alpha_k+\alpha_q-1=deg_{x_{i_k}}(w)=$
the exponent in $x_{i_k}$ of $w$. Thus, $j_e=i_e$ and the induction holds.

We get $w=\prod_{e=1}^{q-1} x_{i_e}^{\alpha_e} \cdot x_{j_q}^{\alpha_q}\cdot y$. If $j_q\in A$, obviously, we get a contradiction. Thus $j_q\in \{i_1,\ldots,i_{q-1}\}$. But this, again, cannot be true, since if $j_q=i_e$, it follows that $x_{i_e}^{\alpha_q+\alpha_e}|w$. In conclusion, our assumption is false and thus $w\notin J_q$.

In order to prove that $s_q\leq \chi_q$, we choose a monomial $w\in J_q^{sat}$ such that $w\notin J_q$ and we show that $deg(w)\leq \chi_q$. Since  $w\in J_q^{sat}$, it follows that
\[ w =  \prod_{e=1}^{q-1} x_{j_e}^{\alpha_e}\prod_{j=1}^{i_q} x_j^{\beta_j}, \]
where $1\leq j_e \leq i_e$ and $\beta_j$ are some nonnegative integers. Since $w\notin J_q$ it follows, obviously, that $\beta_j\leq \alpha_q-1$ and, therefore, $deg(w) = \alpha_1+\cdots+\alpha_{q-1} + \sum_{j=1}^{i_q}\beta_j \leq \chi_q$, as required. 
\end{proof}

\begin{exm}{\em Let $u=x_2^{7}x_3^{6}\in S=K[x_1,x_2,x_3]$. We have $i_1=2$, $i_2=3$, $\alpha_1=7$ and $\alpha_2=6$.
     From Lemma $1.3.4$ it follows that $I = SBT(u) =
	   (x_1^7,x_2^7)(x_1^6,x_2^6,x_3^6)$. 	      
	   We have $J_1=(x_1^7,x_2^7)\subset K[x_1,x_2]$ and
	   $J_2 = I$. Also, $J_1^{sat}=K[x_1,x_2]$ and $J_2^{sat} = (x_1^7,x_2^7) \subset S$.
	   $\chi_1=(\alpha_1-1)\cdot 2 = 12$. $\chi_2=\alpha_1 + (\alpha_2-1)\cdot 3 = 7+15=22$. By $1.3.6$, $reg(I)=max\{12,22\}+1=23$.
	}\end{exm}      


\section{$\de$-fixed ideals.}

In the following $\de: 1=d_{0}|d_{1}|\cdots|d_{s}$ is a strictly increasing sequence of positive integers, where $s$ is a positive integer. We say that $\de$ is a $\de$-sequence. We can take also $s=+\infty$, but for convenience we will not do it.

\begin{lema}
Let $\de$ be a $\de$-sequence. Then, for any $a\in \mathbb N$, there exists an unique sequence of positive integers $a_{0},a_{1},\ldots,a_{s}$ such that:
\begin{enumerate}
    \item $a= \sum_{t=0}^{s}a_{t}d_{t}$ and
    \item $0 \leq a_{t} < \frac{d_{t+1}}{d_{t}}$, for any $0\leq t < s$.
\end{enumerate}
Conversely, if $\de: 1=d_{0}<d_{1}<\cdots<d_{s}$ is a sequence of positive integers such that for any
$a \in \mathbb N$ there exists an unique sequence of positive integers $a_{0},a_{1},\ldots,a_{s}$ as before, then
$\de$ is a $\de$-sequence.
\end{lema}

\begin{proof}
Let $a_{s}$ be the quotient of $a$ divided by $d_{s}$. For $0\leq t<s$ let $a_{t}$ be the quotient of $(a-q_{t+1})$ divided by $d_{t}$, where $q_{t+1}=\sum_{j=t+1}^{s}a_{j}d_{j}$.
We will prove that $a_{0},a_{1},\ldots,a_{s}$ fulfill the required conditions. Indeed, it is obvious that $a= \sum_{t=0}^{s}a_{t}d_{t}$. On the other hand, $a-q_{t+1}<d_{t+1}$, therefore, since $a_{t}$ is the quotient of $(a-q_{t+1})$ divided by $d_{t}$, it follows that $a_{t}<\frac{d_{t+1}}{d_{t}}$.

Suppose there exists another decomposition $a = \sum_{j=0}^{s}b_{j}d_{j}$ which also fulfill the conditions $1$ and $2$. Then, we may assume that there exists an integer $0\leq t\leq s$ such that $b_{s}=a_{s}, \cdots, b_{t+1}=a_{t+1}$ and $b_{t}>a_{t}$.
Notice that $d_{t}>\sum_{j=0}^{t-1}a_{j}d_{j}$. Indeed, \[ \sum_{j=0}^{t-1}a_{j}d_{j}\leq \sum_{j=0}^{t-1}(\frac{d_{j+1}}{d_{j}} - 1)d_{j} = (d_{1}-d_{0}) + (d_{2}-d_{1}) + \cdots + (d_{t}-d_{t-1}) = d_{t} - 1<d_{t}. \]
We have $0 = \sum_{j=0}^{s}(b_{j}-a_{j})d_{j} = \sum_{j=0}^{t}(b_{j}-a_{j})d_{j}$, but on the other hand:
\[ (b_{t}-a_{t})d_{t} \geq d_{t} > \sum_{j=0}^{t-1}a_{j}d_{j} \geq \sum_{j=0}^{t-1}(a_{j}-b_{j})d_{j} \]
and therefore $(b_{t}-a_{t})d_{t} - \sum_{j=0}^{t-1}(a_{j}-b_{j})d_{j} = \sum_{j=0}^{t}(b_{j}-a_{j})d_{j} > 0$, which is a contradiction.

For the converse, we use induction on $0 \leq t< s$, the assertion being obvious for $t=0$. Suppose $t>0$ and $d_{0}|d_{1}|\cdots |d_{t}$ and consider the decomposition of $d_{t+1}-1$. Since $d_{t+1}-1<d_{t+1}$, it follows that $d_{t+1}-1 = \sum_{j=0}^{t}a_{t}d_{t}$. On the other hand, since $d_{t+1}-1$ is the largest integer less than  $d_{t+1}$, each $a_{j}$ is maximal between the integers $<d_{j+1}/d_{j}$, for $j<t$. Therefore $a_{j}=d_{j+1}/d_{j}-1$ for $0 \leq j< t$. Thus:
\[ d_{t+1}  = 1 +d_{t+1} - 1  = 1 + a_{0}d_{0} + a_{1}d_{1} + \cdots + a_{t}d_{t} = d_{1} + a_{1}d_{1} + a_{2}d_{2} + \cdots + a_{t}d_{t} = \]\[ = d_{2} + a_{2}d_{2} + \cdots + a_{t}d_{t} = \cdots = (a_{t}+1)d_{t},\;\;so\;\; d_{t}|d_{t+1}.\]
\end{proof}

\begin{dfn}
Let $a,b$ be two positive integers and consider the $\de$-decompositions \linebreak $a=\sum_{j=0}^{s}a_{j}d_{j}$ and $b=\sum_{j=0}^{s}b_{j}d_{j}$. We say that $a\leq_{\de}b$ if $a_{j}\leq b_{j}$ for any $0\leq j\leq s$.
\end{dfn}

\begin{dfn}
We say that a monomial ideal $I\subset S = k[x_{1},\ldots,x_{n}]$ is $\de$-fixed, if for any monomial $u\in I$ and for any indices $1\leq j<i \leq n$, if $t\leq_{\de} \nu_{i}(u)$ (where $\nu_{i}(u)$ denotes the exponent of the variable $x_{i}$ in $u$) then $u \cdot x_{j}^{t}/x_{i}^{t} \in I$.
\end{dfn}

\begin{exm}
(a) Let $\de: 1|p|p^{2}|p^{3}|\cdots$, where $p>0$ is a prime number. Then a $\de$-fixed ideal $I$ is a $p$-Borel ideal.
Therefore, definition $1.4.3$ generalize the definition of a $p$-Borel ideal.

(b) Suppose $\de: 1$, i.e. $s=0$. Let $I\subset S$ be a monomial ideal. Then $I$ is $\de$-fixed, 
if and only if $I$ is strongly stable. Indeed, in the definition of a $\de$-fixed ideal, we can always choose $t=1$,
because the $\de$-decomposition of any positive integer $m$ is in this case $m = m\cdot 1$.

(c) Let $\de: 1=d_0|d_1|\cdots|d_s$ be a $\de$-sequence. If $I,J\subset S$ are $\de$-fixed ideals, then
    $I+J$ and $I\cap J$ are also $\de$-fixed ideals. This is obvious from the definition.
\end{exm}

\pagebreak
\begin{lema}
Let $a,b$ be two positive integers with $a\leq_{\de} b$. Suppose $b=b'+b''$, where $b'$ and $b''$ are positive integers. Then, there exists some positive integers $a'\leqd b'$ and $a''\leqd b''$ such that $a=a'+a''$.
\end{lema}

\begin{proof}
Let $a=\sum_{t=0}^{s}a_{t}d_{t},\; b = \sum_{t=0}^{s}b_{t}d_{t}, \; b'=\sum_{t=0}^{s}b'_{t}d_{t}, \; b''=\sum_{t=0}^{s}b''_{t}d_{t}$. The hypothesis implies $a_{t}\leq b_{t} < d_{t+1}/d_{t}$ and $b'_{t}, b''_{t} < d_{t+1}/d_{t}$ for any $0\leq t<s$. We construct the sequences $a'_{t}, a''_{t}$ using decreasing induction on $t$.
Suppose we have already defined $a'_{j}, a''_{j}$ for $j>t$ such that $\sum_{i=j}^{s}(a'_{i}+a''_{i})d_{i} = \sum_{i=j}^{s}a_{i}d_{i}$ and $b_{t+1} = b'_{t+1}+b''_{t+1}$. This is obvious for $t=s$.

We consider two cases. If $b_{t} = b'_{t} + b''_{t}$, then we choose $a'_{t}\leq b'_{t}$ and $a''_{t}\leq b''_{t}$ such that $a'_{t}+a''_{t}=a_{t}$. We can do this, because $a_{t}\leq b_{t}$. Also, it is obvious from the induction hypothesis that $\sum_{i=t}^{s}(a'_{i}+a''_{i})d_{i} = \sum_{i=t}^{s}a_{i}d_{i}$, so we can pass from $t$ to $t-1$.

If $b_{t} \neq b'_{t} + b''_{t}$ we claim that $b'_{t} + b''_{t} = b_{t}-1$. Indeed, $\sum_{j=0}^{t-1}(b'_{j}+b''_{j})d_{j}<2d_{t}$ and therefore it is impossible to have $b'_{t}+b''_{t} \leq b_{t}-2$, otherwise $\sum_{j=0}^{t}(b'_{j}+b''_{j})d_{j}<b_{t}d_{t}$ and we contradict the equality $b=b'+b''$. Also, since $b'_{t+1} + b''_{t+1} = b_{t+1}$, we cannot have $b'_{t}+b''_{t}>b_{t}$. Similarly we get $b'_{t-1} + b''_{t-1} > b_{t-1}$. By recurrence, we conclude that there exists an integer $u<t$ such that:
$ b'_{u-1}+b''_{u-1} = b_{u-1},\;,b'_{u}+b''_{u} = b_{u} + d_{u+1}/d_{u}$,
$b'_{u+1} + b''_{u+1} = b_{u+1} + d_{u+2}/d_{u+1} - 1,\; \ldots, \; b'_{t-1} + b''_{t-1} = b_{t-1} + d_{t}/d_{t-1} - 1$.

If $a_{j}=b_{j}$ for any $j\in \{u,\ldots,t\}$, we simply choose $a'_{j}=b'_{j}$ and $a''_{j}=b''_{j}$ for any
$j\in \{u,\ldots,t\}$ and the required conditions are fulfilled, so we can pass from $t$ to $u-1$. If this is not the case, then there exists an integer $u \leq q \leq t$ such that $a_{t}=b_{t},\ldots,a_{q+1} = b_{q+1}$ and $a_{q}<b_{q}$. If $q=t$ then for any $j\in \{u,\ldots,t\}$ we can choose $a'_{j}\leq b'_{j}$ and $a''_{j}\leq b''_{j}$ such that $a'_{j}+a''_{j} = a_{j}$. For $j<t$ the previous assertion is obvious because $b'_{j}+b''_{j}\geq b_{j}$, and for $j=t$, since $a_{t}<b_{t}$ we have in fact $a_{t}\leq b'_{t}+b''_{t} = b_{t}-1$ and therefore we can choose again $a'_{t}$ and $a''_{t}$. The conditions are satisfied so we can pass from $t$ to $u-1$.

Suppose $q<t$. For $j\in \{u,\ldots,q-1\}$ we choose $a'_{j}\leq b'_{j}$ and $a''_{j}\leq b''_{j}$ such that $a'_{j}+a''_{j}=a_{j}$. We can do this because $b'_{j}+b''_{j}\geq b_{j} \geq a_{j}$. We choose
$a'_{q}$ and $a''_{q}$ such that $a'_{q}+a''_{q} =a_{q} + d_{q+1}/d_{q}$. We can make this choice, because $a_{q}\leq b_{q}-1$ and $b'_{q}+b''_{q} \geq b_{q} + d_{q+1}/d_{q} - 1$. For $j>q$, we simply put $a'_{j}=b'_{j}$ and $a''_{j}=b''_{j}$. To pass from $t$ to $u-1$ is enough to see that $\sum_{j=u}^{t}a_{j}d_{j} = \sum_{j=u}^{t}(a'_{j}+a''_{j})d_{j}$. Indeed,
\[ \sum_{j=u}^{t}(a'_{j}+a''_{j})d_{j} = \sum_{j=u}^{q-1}(a'_{j}+a''_{j})d_{j} + (a'_{q}+a''_{q})d_{q} + \sum_{j=q+1}^{t}(a'_{j}+a''_{j})d_{j} = \]
\[ = \sum_{j=u}^{q-1}a_{j}d_{j} + (a_{q}+d_{q+1}/d_{q})d_{q} + \sum_{j=q+1}^{t-1}(a_{j}+ d_{j+1}/d_{j}-1)d_{j} + (a_{t}-1)d_{t} =\]\[ = \sum_{j=u}^{t}a_{j}d_{j} + d_{q+1} + \sum_{j=q+1}^{t-1}(d_{j+1}-d_{j}) - d_{t} = \sum_{j=u}^{t}a_{j}d_{j},\]
The induction ends when $t=-1$. Finally, we obtain $a'$ and $a''$ such that $a'+a''=a$, $a'_{t}\leq b'_{t}$ and $a''_{t}\leq b''_{t}$, as required.
\end{proof}

\pagebreak
\begin{cor}
If $I,J\subset S$ are $\de$-fixed ideals then $I\cdot J$ is a $\de$-fixed ideal.
\end{cor}

\begin{proof}
Let $u\in I\cdot J$ be a monomial and let $1\leq j<i\leq n$ be two integers.
We can write $u=v\cdot w$, where $v\in I$ and $w\in J$ are monomials.
Let $t\leq_{\de} \nu_i(u)$ be a positive integer. Since $\nu_i(u)=\nu_i(v)+\nu_i(w)$, by previous lemma, we
can choose two positive integers $t'\leq_{\de} \nu_i(v)$ and $t''\leq_{\de} \nu_i(w)$ such that
$t=t'+t''$. Since $I$ and $J$ are $\de$-fixed, it follows that $x_j^{t'}v/x_i^{\nu_i(v)}\in I$ and
$x_j^{t''}w/x_i^{\nu_i(w)}\in J$. Therefore, $x_j^{t}u/x_i^{\nu_i(u)}\in I\cdot J$ and thus $I\cdot J$ is
$\de$-fixed, as required.
\end{proof}

\begin{dfn}
A $\de$-fixed ideal $I$ is called {\em principal} if it is generated, as a $\de$-fixed ideal by one monomial $u\in S$, i.e. $I$ is the smallest $\de$-fixed ideal which contains $u$. We write $I = <u>_{\de}$.
More generally, if $u_{1},\ldots,u_{r} \in S$ are monomials, the $\de$-fixed ideal generated by $u_{1},\ldots,u_{r}$ is the smallest $\de$-fixed ideal $I$ which contains $u_{1},\ldots,u_{r}$. We write $I = <u_{1},\ldots,u_{r}>_{\de}$.
\end{dfn}

Our next goal is to describe the principal $\de$-fixed ideals. The easiest case is when we have a $\de$-fixed ideal generated by the power of a variable. We denote $\me=\{x_{1},\ldots,x_{n}\}$ and $\me^{[d]}=\{x_{1}^{d},\ldots,x_{n}^{d}\}$, where $d$ is a positive integer. We have the following proposition.

\begin{prop}
If $u=x_{n}^{\alpha}$, then $I= <u>_{\de} = \prod_{t=0}^{s} (\me^{[d_{t}]})^{\alpha_{t}}$,
where $\alpha=\sum_{t=0}^{s}\alpha_{t}d_{t}$.
\end{prop}

\begin{proof}
Let $I' = \prod_{t=0}^{s} (\me^{[d_{t}]})^{\alpha_{t}}$. The minimal generators of $I'$ are monomials of the type
$w = \prod_{t=0}^{s}\prod_{j=1}^{n} x_{j}^{\lambda_{tj}\cdot d_{t}}$, where $0\leq \lambda_{tj}$ and $\sum_{j=1}^{n}\lambda_{tj} = \alpha_{t}$. First, let us show that $I'\subset I$. In order to do this, we choose $w$
a minimal generator of $I'$ (the one bellow). We write $x_{n}^{\alpha}$ like this: $x_{n}^{\alpha} = x_{n}^{\alpha_{0}d_{0} + \alpha_{1}d_{1} + \cdots + \alpha_{s}d_{s}} = x_{n}^{\alpha_{0}d_{0}} \cdot x_{n}^{\alpha_{1}d_{1}} \cdots x_{n}^{\alpha_{s}d_{s}}$.
Since $\lambda_{01}d_{0}\leq_{\de}\alpha_{0}d_{0} + \alpha_{1}d_{1} + \cdots + \alpha_{s}d_{s}$ and $I$ is $\de$-fixed it follows that $x_{1}^{\lambda_{01}d_{0}}x_{n}^{\alpha-\lambda_{01}d_{0}} \in I$. Also, $\lambda_{02}d_{0} < \alpha-\lambda_{01}d_{0} = (\alpha_{0}-\lambda_{01})d_{0} + \alpha_{1}d_{1} + \cdots + \alpha_{s}d_{s}$, and since $I$ is $\de$-fixed it follows that $x_{1}^{\lambda_{01}d_{0}} x_{2}^{\lambda_{02}d_{0}}x_{n}^{\alpha-\lambda_{01}d_{0}-\lambda_{02}d_{0}} \in I$.
Using iteratively this argument, one can easily see that $x_{1}^{\lambda_{01}d_{0}} \cdots x_{n}^{\lambda_{0n}d_{0}}x_{n}^{\alpha-\alpha_{0}d_{0}} \in I$. Also $\alpha - \alpha_{0}d_{0} = \alpha_{1}d_{1} + \cdots +\alpha_{n}d_{n}$. Again, using an inductive argument, we get:
\[ (x_{1}^{\lambda_{01}d_{0}} \cdots x_{n}^{\lambda_{0n}d_{0}})\cdot (x_{1}^{\lambda_{11}d_{1}} \cdots x_{n}^{\lambda_{1n}d_{1}})\cdots (x_{1}^{\lambda_{s1}d_{s}} \cdots x_{n}^{\lambda_{sn}d_{s}}) = w \in I.\]
For the converse, i.e. $I\subset I'$, is enough to verify that $I'$ is $\de$-fixed. In order to do this, is enough to
prove that the minimal generators of $I'$ fulfill the definition of a $\de$-fixed ideal. Let $w$ be a minimal generator of  $I'$. Let $2\leq i \leq n$. Then $\nu_{i}(w) = \sum_{t=0}^{s}\lambda_{ti}d_{t}$. If $\beta\leqd \nu_{i}(w)$ then  $\beta=\sum_{t=0}^{s}\beta_{t}d_{t}$ with $\beta_{t}\leq \lambda_{ti}$. Let $1\leq k<i$. We have
 \[ w\cdot x_{k}^{\beta}/x_{i}^{\beta} = \prod_{t=0}^{s}(\prod_{j\neq i,k}x_{j}^{\lambda_{tj}d_{t}})\cdot x_{i}^{(\lambda_{ti} - \beta_{t})d_{t}}\cdot x_{k}^{(\lambda_{tk} + \beta_{t})d_{t}}. \]
Thus $w\cdot x_{k}^{\beta}/x_{i}^{\beta} \in I'$ and therefore $I'$ is $\de$-fixed. Since $I$ is the smallest $\de$-fixed ideal which contains $x_{n}^{\alpha}$ it follows that $I\subset I'$.
\end{proof}

\pagebreak
\begin{prop}
If $\alpha \leq \beta$ then $<x_{n}^{\beta}>_{\de} \subseteq <x_{n}^{\alpha}>_{\de}$.
\end{prop}

\begin{proof}
The case $\alpha=\beta$ is obvious, so we may assume $\alpha<\beta$. We denote $I=<x_{n}^{\alpha}>_{\de}$ and $I' = <x_{n}^{\beta}>_{\de}$. We write $\alpha=\sum_{t=0}^{s}\alpha_{t}d_{t}$ and $\beta=\sum_{t=0}^{s}\beta_{t}d_{t}$. If $w$ is a minimal generator of $I'$ then $w=\prod_{t=0}^{s}\prod_{i=1}^{n}x_{i}^{\lambda_{ti}d_{t}}$, where $0\leq \lambda_{ti}$ and $\sum_{i=1}^{n}\lambda_{ti}=\beta_{t}$. We claim that $w\in I$ and therefore $I'\subset I$ as required.

Since $\alpha < \beta$ there exists $t\in \{0,\ldots,s\}$ such that $\alpha_{s} = \beta_{s},\ldots,\alpha_{t+1}=\beta_{t+1}$ and $\alpha_{t}<\beta_{t}$. We may assume some $\lambda_{tk}>0$. We have
\[ w = \prod_{j=0}^{s}\prod_{i=1}^{n}x_{i}^{\lambda_{ji}d_{j}} = \prod_{j=0}^{t-1}x_{k}^{\alpha_{j}d_{j}} x_{k}^{(\lambda_{tk}-1)d_{t}} x_{k}^{d_{t}-\sum_{j=0}^{t-1}\alpha_{j}d_{j}} \prod_{i\neq k}^{n}x_{i}^{\lambda_{ti}d_{t}}\prod_{j>t}\prod_{i=1}^{n}x_{i}^{\lambda_{ji}d_{j}} \]
and now it is obvious that $w\in I$.
\end{proof}

\begin{prop}
If $\alpha$ and $\beta$ are two positive integers, then  $<x_n^{\alpha+\beta}>_{\de}\subseteq \linebreak
<x_n^{\alpha}>_{\de} \cdot <x_n^{\beta}>_{\de}$. The equality holds if and only if $\alpha_t+\beta_t < d_{t+1}/d_t$
for all $0\leq t\leq s$, where $\alpha=\sum_{t=0}^s \alpha_td_t$ and $\beta=\sum_{t=0}^s \beta_td_t$
\end{prop}

\begin{proof}
We denote $I=<x_n^{\alpha+\beta}>$ and $I'=<x_n^{\alpha}>\cdot <x_n^{\beta}>$. Also, we
denote $\gamma:=\alpha+\beta$. Let $\gamma=\sum_{t=0}^s \gamma_td_t$. We may assume $\gamma_s\neq 0$, otherwise, we replace $s$ with  $s'=max\{t|\;\gamma_t\neq 0\}$. We use induction on $t=min\{j |\;\alpha_j^2+\beta_j^2\neq 0\}$. It $t=s$, it follows that $\alpha=\alpha_sd_s$, $\beta=\beta_sd_s$ and thus $I=I'$. Suppose $t<s$. We should consider two cases: (I) $\alpha_t+\beta_t = \gamma_t$ or (II) $\alpha_t+\beta_t = \gamma_t + d_{t+1}/d_t$.

(I) Let $\bar{\alpha}:=\alpha-\alpha_td_t$ and $\bar{\beta}=\beta-\beta_td_t$. Let $\bar{\gamma}:=\bar{\alpha} +
\bar{\beta}$. We denote $\bar{I}:=<x_n^{\bar{\gamma}}>_{\de}$ and $\bar{I}':=<x_n^{\bar{\alpha}}>_{\de}\cdot
<x_n^{\bar{\beta}}>_{\de}$. By induction hypothesis, we have $\bar{I} \subseteq \bar{I}'$. Therefore,
since $I=(\me^{[d_t]})^{\gamma_t}\bar{I}$ and $I'=(\me^{[d_t]})^{\gamma_t}\bar{I}'$ it follows that $I\subseteq I'$
as required.

(II)As above, we define $\bar{\alpha}$, $\bar{\beta}$, $\bar{\gamma}$, $\bar{I}$ and $\bar{I'}$.
We notice that $I=(\me^{[d_t]})^{\alpha_t+\beta_t - d_{t+1}/d_t} \cdot (\me^{[d_{t+1}]}) \bar{I}$ and
$I'=(\me^{[d_t]})^{\alpha_t+\beta_t} \bar{I}'$. Using induction hypothesis, it follows that $I\subset I'$. 
Note that in this case, the inclusion is strict, therefore we get the second statement of the proposition.
\end{proof}

We have the general description of a principal $\de$-fixed ideal given by the following proposition. In the proof, we will apply Lemma $1.4.5$.

\begin{prop}
Let $1\leq i_{1}<i_{2}<\cdots<i_{r}=n$ and let $\alpha_{1},\ldots,\alpha_{r}$ be some positive integers.
If $u=x_{i_{1}}^{\alpha_{1}}x_{i_{2}}^{\alpha_{2}} \cdots x_{i_{r}}^{\alpha_{r}}$ then:
\[ I=<u>_{\de} = <x_{i_{1}}^{\alpha_{1}}>_{\de}\cdot <x_{i_{2}}^{\alpha_{2}}>_{\de} \cdots <x_{i_{r}}^{\alpha_{r}}>_{\de} = \prod_{q=1}^{r}\prod_{t=0}^{s}(\me_{q}^{[d_{t}]})^{\alpha_{qt}},\]
where $\me_{q} = \{x_{1},\ldots,x_{i_{q}}\}$, $\me_{q}^{[d_t]} = \{x_{1}^{d_t},\ldots,x_{i_{q}}^{d_t}\}$ and $\alpha_{q}=\sum_{t=0}^{s}\alpha_{qt}d_{t}$.
\end{prop}

\begin{proof}
Let $I' = \prod_{q=1}^{r}\prod_{t=0}^{s}(\me_{q}^{[d_{t}]})^{\alpha_{qt}}$. The minimal generators of $I'$ are monomials of the type $w = \prod_{q=1}^{r} \prod_{t=0}^{s}\prod_{j=1}^{i_{q}} x_{j}^{\lambda_{qtj}\cdot d_{t}}$, where $0\leq \lambda_{qtj}$ and $\sum_{j=1}^{n}\lambda_{qtj} = \alpha_{qt}$. First, we show that $I'\subset I$. In order to do this,
it is enough to prove that by iterative transformations we can modify $u$ such that we obtain $w$.

The idea of this transformations is the same as in the proof of $1.4.8$. Without given all the details, one can see that if we rewrite $u$ as \[ (x_{i_{1}}^{\alpha_{10}d_{0}}x_{i_{1}}^{\alpha_{11}d_{1}}\cdots x_{i_{1}}^{\alpha_{1s}d_{s}}) \cdots (x_{i_{r}}^{\alpha_{r0}d_{0}}x_{i_{r}}^{\alpha_{r1}d_{1}}\cdots x_{i_{r}}^{\alpha_{rs}d_{s}} ),\]
where $\alpha_{q}=\sum_{t=0}^{s}\alpha_{qt}d_{t}$, we can pass to $w$, using the transformations
\[x_{i_{1}}^{\alpha_{10}d_{0}} \mapsto \prod_{j=1}^{i_{1}}x_{j}^{\lambda_{10j}d_{0}},\ldots, x_{i_{1}}^{\alpha_{1s}d_{s}} \mapsto \prod_{j=1}^{i_{1}}x_{j}^{\lambda_{1sj}d_{s}}, \ldots , x_{i_{r}}^{\alpha_{r0}d_{0}} \mapsto \prod_{j=1}^{i_{r}}x_{j}^{\lambda_{r0j}d_{0}}, \ldots, x_{i_{r}}^{\alpha_{rs}d_{s}} \mapsto \prod_{j=1}^{i_{r}}x_{j}^{\lambda_{rsj}d_{s}}.\]
Therefore $w \in I$, and thus $I'\subset I$. For the converse, it is enough to see that $I'$ is a $\de$-fixed ideal. Let $w$ be a minimal generator of $I'$. We choose an index $2 \leq i \leq n$. Then $\nu_{i}(w) = \sum_{q=1}^{r}\sum_{t=0}^{s}\lambda_{qti}d_{t}$. Let $\beta \leq \nu_{i}(w)$. Using Lemma $1.4.5$, we can choose some positive integers $\beta_{1},\ldots,\beta_{r}$ such that:
\[ (a) \beta = \sum_{q=1 , i_{q}\geq i}^{r} \beta_{q} \;and\; (b) \beta_{q}\leqd \sum_{t=0}^{s}\lambda_{qti}d_{t},\; \]
i.e. $\beta_{qt}\leq \lambda_{qti}$, where $\beta_{q} = \sum_{t=0}^{s}\beta_{qt}d_{t}$. Let $k<i$. Then,
\[ w\cdot x_{k}^{\beta}/x_{i}^{\beta} = \prod_{q=1}^{r}  \prod_{t=0}^{s}\left( \prod_{j=1, j\neq k,i}^{i_{q}} x_{j}^{\lambda_{qtj}\cdot d_{t}} \right)
x_{i}^{(\lambda_{qti} - \beta_{qt})d_{t}} x_{k}^{(\lambda_{qtk} + \beta_{qt})d_{t}}.\]
Now, it is easy to see that $w\cdot x_{k}^{\beta}/x_{i}^{\beta} \in I'$, and therefore $I'$
is $\de$-fixed.
\end{proof}

\begin{exm}{\em
Let $\de: 1|2|4|12$.
\begin{enumerate}
    \item Let $u=x_{3}^{21}$. We have $21=1\cdot 1 + 0\cdot 2 + 2 \cdot 4 + 1 \cdot 12$. From $1.4.8$, we get:
    \[ <u>_{\de} = (x_{1},x_{2},x_{3})(x_{1}^{4},x_{2}^{4},x_{3}^{4})^{2}(x_{1}^{12},x_{2}^{12},x_{3}^{12}).\]
    \item Let $u=x_{1}^{2}x_{2}^{9}x_{3}^{16}$. We have $9=1\cdot 1 + 2\cdot 4$ and $16 = 1\cdot 4 + 1\cdot 12$. From
          $1.4.11$, we get
    \[ <u>_{\de} = x_{1}^{2}<x_{2}^{9}>_{\de}<x_{3}^{16}>_{\de} = x_{1}^{2}(x_{1},x_{2})(x_{1}^{4},x_{2}^{4})^{2}
                    (x_{1}^{4},x_{2}^{4},x_{3}^{4})(x_{1}^{12},x_{2}^{12},x_{3}^{12}).\]
\end{enumerate}
}\end{exm}

\begin{obs}{\em
Any $\de$-fixed ideal $I$ is a Borel type ideal.
Indeed, Proposition $1.1.4(d)$ says that an ideal $I$ is of
Borel type if and only if for any $1\leq j<i\leq n$, there exists an
positive integer $t$ such that $x_{j}^{t}(u/x_{i}^{\nu_{i}(u)})\in
I$. Choosing $t=\nu_{i}(u)$, is easy to see that the definition of a
$\de$-fixed ideal implies the condition above.

Let $u = x_{i_{1}}^{\alpha_{1}}x_{i_{2}}^{\alpha_{2}} \cdots x_{i_{r}}^{\alpha_{r}}$ and
$ I=<u>_{\de} =  \prod_{q=1}^{r}\prod_{t=0}^{s}(\me_{q}^{[d_{t}]})^{\alpha_{qt}}$,
where $\alpha_{q}=\sum_{t=0}^{s}\alpha_{qt}d_{t}$.
Let $I_{r-e} = \prod_{q=1}^{e}\prod_{t=0}^{s}(\me_{q}^{[d_{t}]})^{\alpha_{qt}}$. Then $I = I_{0} \subset I_{1} \subset \cdots \subset I_{r}=S$ is the sequential chain of $I$. Let $n_{\ell}=i_{q_{r-\ell}}$.
Indeed, since $x_{n_{\ell}}^{\alpha_{r-e}}I_{\ell+1} \subset I_{\ell} \Rightarrow I_{\ell +1}\subset (I_{\ell}:x_{n_{\ell}}^{\infty})$. For the converse, let $w\in (I_{\ell}:x_{n_{\ell}}^{\infty})$ be any minimal generator. Then there exists an integer $b$ such that $w\cdot x_{n_{\ell}}^{b}\in I_{\ell}$. We may assume that $w$ is a minimal generator of $I_{\ell}$. Then $w\cdot x_{n_{\ell}}^{b} = w'\cdot y$ for a $w'\in I_{\ell+1}$ and $y\in \prod_{j=0}^{t}(\me_{r-\ell}^{[d_{j}]})^{\alpha_{r-\ell,j}}$ with $x_{n_{\ell}}^{b}|y$. Thus $w'|w$, and therefore $w\in I_{\ell+1}$.
}\end{obs}

\section{Socle of factors by principal $\de$-fixed ideals.}

In the following, we suppose $n\geq 2$.

\begin{lema}
Let $\de:1=d_{0}|d_{1}| \cdots |d_{s}$, $\alpha\in \mathbb N$ and $I=<x_{n}^{\alpha}>_{\de} =\prod_{t=0}^{s} (\me^{[d_{t}]})^{\alpha_{t}}$. Let $q_{t}=\sum_{j=t}^{s}\alpha_{j}d_{j}$. Let
\[ J = \sum_{t=0 , \alpha_{t}>0}^{s} (x_{1}\cdots x_{n})^{d_{t}-1} (\me^{[d_{t}]})^{\alpha_{t}-1} \prod_{j>t} (\me^{[d_{j}]})^{\alpha_{j}}.\]
Then:
\begin{enumerate}
    \item $Soc(S/I) = \frac{J+I}{I}$
    \item Let $e$ be a positive integer. Then $(\frac{J+I}{I})_{e}\neq 0 \Leftrightarrow
          e = q_{t} + (n-1)(d_{t}-1) -1$, for some $0\leq t\leq s$ with $\alpha_{t}>0$.
    \item $max\{e| (\frac{J+I}{I})_{e}\neq 0 \} = \alpha_{s}d_{s} + (n-1)(d_{s}-1) - 1$.
\end{enumerate}
\end{lema}

\begin{proof}
1. Firstly, we prove that $\frac{J+I}{I} \subset Soc(S/I)$. Since $Soc(S/I) = (O:_{S/I}\mathbf{m})$, it is enough to show that $\textbf{m}J \subset I$.

We have $J=\sum_{t=0,\;\alpha_{t}>0}^{s}J_{t}$, where $J_{t} = (x_{1}\cdots x_{n})^{d_{t}-1} (\me^{[d_{t}]})^{\alpha_{t}-1} \prod_{j>t} (\me^{[d_{j}]})^{\alpha_{j}}$. It is enough
to prove that $x_{i}J_{t}\subset I$ for any $i$ and any $t$. Suppose $i=1$:
\[x_{1}J_{t} = x_{1}^{d_{t}}(x_{2}\cdots x_{n})^{d_{t}-1} (\me^{[d_{t}]})^{\alpha_{t}-1} \prod_{j>t} (\me^{[d_{j}]})^{\alpha_{j}} \subset (x_{2}\cdots x_{n})^{d_{t}-1} \prod_{j\geq t} (\me^{[d_{j}]})^{\alpha_{j}}. \]
On the other hand, $(x_{2}\cdots x_{n})^{d_{t}-1} \in \prod_{j<t}(\me^{[d_{j}]})^{\alpha_{j}}$, because $d_{t}-1\geq \sum_{j<t}\alpha_{j}d_{j}$. Thus $x_{1}J_{t}\subset I$.

For the converse, we apply induction on $\alpha$. If $\alpha=1$ then $s=0$ and $I=(x_{1},\ldots,x_{n}) = \mathbf{m}$. $J=(x_{1},\ldots,x_{n})^{d_{0}-1}=S$, and obvious $Soc(S/I) = Soc(S/\mathbf{m}) = S/\mathbf{m}$. Let us suppose that
$\alpha>1$. We prove that if $w\in S\setminus I$ is a monomial such that $\mathbf{m}w\subset I$, then $w\in J$. Let $t_{e}=max\{ t:\; x_{e}^{d_{t}-1}|w \}$. Renumbering $x_{1},\ldots,x_{n}$ which does not affect either $I$ or $J$, we may suppose that $t_{1}\geq t_{2}\geq \cdots \geq t_{n}$. We have two cases: (i)$t_{1}>t_{n}$ and (ii)$t_{1}=t_{n}$.
But first, let's make the following remark: $(*)$ If $u=x_{1}^{\beta_{1}}\cdots x_{n}^{\beta_{n}} \in \prod_{j\geq t}\me^{d_{j}}$ and $\beta_{i}<d_{t}$ for certain $i$ then $u/x_{i}^{\beta_{i}} \in \prod_{j\geq t}\me^{[d_{j}]}$ (the proof is similarly to \cite[Lemma 3.5]{hpv}).

In the case (i), there exists an index $e$ such that $t_{e}>t_{e+1} = \cdots = t_{n}$. Then we have $w=(x_{n}\cdots x_{e+1})^{d_{t_{n}}-1}\cdot x_{e}^{d_{t_{e}}-1}\cdot y$, for a monomial $y\in S$. We consider two cases (a) $x_{e}$ does not divide $y$ and (b) $x_{e}$ divide $y$. (a) From $x_{n}w = x_{n}^{d_{t_{n}}}\cdot (x_{n-1}\cdots x_{e+1})^{d_{t_{n}}-1}x_{e}^{d_{t_{e}}-1}\cdot y \in I$ we see that $y \in \prod_{j\geq t_{e}} (\mathbf{m}^{[d_{j}]})^{\alpha_{j}}$, by $(*)$. Therefore $w\in I$, because
$x_{e}^{d_{t_{e}}-1} \in \prod_{j<t_{e}}(\me^{[d_{j}]})^{\alpha_{j}}$, which is an contradiction.

(b) In this case, $w = (x_{n}\cdots x_{e+1})^{d_{t_{n}}-1} x_{e}^{d_{t_{e}}}y'$, where $y'=y/x_{e}$. We claim that there exist $\lambda\leq t_{e}$ such that $\alpha_{\lambda}\neq 0$. Indeed, if all $\alpha_{\lambda}=0$ for $\lambda\leq t_{e}$, then $I = \prod_{j=t_{e}+1}^{s}(\mathbf{m}^{[d_{j}]})^{\alpha_{j}}$ and $x_{n}w\in I$ implies $y'\in I$ because of the maximality of $t_{n}$ and $(*)$. It follows $w\in I$, which is false. Choose $\lambda\leq t_{e}$ maximal possible with $\alpha_{\lambda}\neq 0$. Set $w' = w/ x_{e}^{d_{\lambda}}$. Note that $\mathbf{m}w\subset I$ implies
\[\mathbf{m}w\subset I' = (\mathbf{m}^{[d_{t_{\lambda}}]})^{\alpha_{\lambda}-1} \prod_{j\neq \lambda} (\mathbf{m}^{[d_{t_{j}}]})^{\alpha_{j}}.\]
It is obvious that $x_{q}w'\in I'$ for $q\neq e$. Also, since $x_{e}^{d_{t_{e}+1}}$ does not divide $x_{e}w$ implies $x_{e}w'\in I'$. Choosing $\alpha'=\alpha-d_{\lambda}$, we get $\alpha'_{j}=\alpha_{j}$ for $j\neq \lambda$ and $\alpha'_{\lambda} = \alpha_{\lambda} - 1$ and therefore we can apply our induction hypothesis for $I'$ (because $\alpha'<\alpha$) and for the ideal $J'$ associated to $I'$, which has the form:
\[ J' = \sum_{q=0,\alpha'_{q}\neq 0} (x_{1}\cdots x_{n})^{d_{q}-1} (\mathbf{m}^{[d_{q}]})^{\alpha'_{q}-1} \prod_{j>q} (\mathbf{m}^{[d_{j}]})^{\alpha'_{j}}, \]
and so $w = x_{e}^{d_{\lambda}}w' \in x_{e}^{d_{\lambda}}J' \subset J$.

It remains to consider the case (ii) in which we have in fact $t_{1}=t_{2}=\cdots = t_{n}$. If $y=w/ (x_{1}\cdots x_{n})^{d_{t_{n}}-1} \in \mathbf{m}$, then there exists $e$ such that $x_{e}|y$, and we apply our induction hypothesis as in the case $(b)$ above. Thus we may suppose $y=1$, i.e. $w=(x_{1}\cdots x_{n})^{d_{t_{n}}-1}$. Since $\mathbf{m}w\subset I$, we see that $\alpha_{j}=0$ for $j>t_{n}$ and $\alpha_{t_{n}}=1$ (otherwise $w\in I$, which is absurd). Thus $w\in J$.

2. Let $v=x_{1}^{q_{t}-1}(x_{2}\cdots x_{n})^{d_{t}-1}$. Then $deg(v)=q_{t} + (n-1)(d_{t}-1) -1$. But $v\in J$ and $v\notin I$, therefore $v\neq 0$ in $Soc(S/I) = \frac{J+I}{I}$.


3. Let $e_{t}= q_{t} + (n-1)(d_{t}-1) -1$ for $0\leq t\leq s$. Let $t<s$. Then
\[ e_{t+1}-e_{t} = q_{t+1}-q_{t} + (n-1)(d_{t+1}-d_{t}) = -\alpha_{t}d_{t} + (n-1)(d_{t+1}-d_{t}) \geq d_{t+1} - (\alpha_{t}+1)d_{t} \geq 0,\;so \]
\[ max\{e| ((J+I)/I)_{e}\neq 0 \} = e_{s} = \alpha_{s}d_{s} + (n-1)(d_{s}-1) - 1. \]
\end{proof}

\begin{obs}{\em
From the proof of the above lemma, we may easily conclude that for
$n\geq 3$, $e_{t}=e_{t'}$ if and only if $t=t'$, and if $n=2$, then $e_{t}=e_{t'}$ ($t<t'$) if and only if \linebreak $\alpha_{t'-1}=d_{t'}/d_{t'-1},\ldots,\alpha_{t}=d_{t+1}/d_{t}$.
}\end{obs}

\begin{cor}
With the notations of previous lemma and remark, let $0\leq t\leq s$ be an integer such that $\alpha_{t}\neq 0$. Let $h_{t}=dim_{K}((I+J_{t})/I)$. Then:
\begin{enumerate}
    \item $G(J_{t})\cap (I+J_{t'}) = 0$ for $0\leq t'\leq s$, $t'\neq t$.
    \item $h_{t}=\binom{n+\alpha_{t}-2}{n-1}\prod_{j>t}\binom{n+\alpha_{j}-1}{n-1}$.
    \item $dim_{K}(Soc(S/I)_{e}) =
\begin{cases}
    h_{q},& if\; n\geq 3\;and\;e=e_{q}\;for\;a\;q\leq s\;with\;\alpha_{q}\neq 0. \\
    \sum_{q}h_{q},& if\; n=2\;and\;q\in\{\epsilon|e=e_{\epsilon}\;for\;\epsilon\leq s\;with\;\alpha_{\epsilon}\neq 0\}.\\
    0,& otherwise.
\end{cases}.$
\end{enumerate}
\end{cor}

\begin{proof}
1. First suppose $t'<t$. A minimal generator $x^{\beta}=x_{1}^{\beta_{1}}\cdots x_{n}^{\beta_{n}}$ of $J_{t}$ has the form \[(x_{1}\cdots x_{n})^{d_{t}-1}\prod_{j\geq t}(x_{1}^{\lambda_{1j}d_{j}}\cdots x_{n}^{\lambda_{nj}d_{j}}),\;where \;\sum_{\nu=1}^{n}\lambda_{\nu j} = \begin{cases} \alpha_{j},& if\; j>t,\\ \alpha_{t}-1,&if\; j=t.  \end{cases}.\] Thus, $\beta_{i} = d_{t}-1 + \sum_{j=t}^{s}\lambda_{ij}d_{j}$. On the other hand, $d_{t}-1 = \sum_{j=0}^{t-1} (d_{j+1}/d_{j}-1)d_{j}$, so $\beta_{i}$ has the writing $\sum_{j=0}^{s}\beta_{ij}d_{j}$, where $\beta_{ij}=d_{j+1}/d_{j}-1$ for $j<t$ and $\beta_{ij}=\lambda_{ij}$ for $j\geq t$.

Assume that $x^{\beta}\in I+J_{t'}$ for a certain $t'<t$. Then there exists $\gamma\in\Nn$ such that $x^{\gamma}\in G(I)$ (or $x^{\gamma}\in G(J_{t'})$) and $x^{\gamma}|x^{\beta}$, that is $\gamma_{i}\leq \beta_{i}$ for all $1\leq i\leq n$. Let $\gamma_{i} = \sum_{j=0}^{s}\gamma_{ij}d_{j}$, the $\de-$ decomposition of $\gamma_{i}$. We notice that $(\beta_{is},\ldots,\beta_{i0})\geq (\gamma_{is},\ldots,\gamma_{i0})$ in the lexicographic order.

Note that all minimal generators $x^{\gamma}$ of $I$ have the same degree $\alpha< e_{t}$ and $\sum_{i=1}^{n}\gamma_{iq} = \alpha_{q}$ for each $0\leq q \leq s$. Also
all minimal generators $x^{\gamma}$ of $J_{t'}$ have the same degree $e_{t'}<e_{t}$ and $\sum_{i=1}^{n}\gamma_{iq}=\alpha_{q}$ for each $t\leq q\leq s$. It follows $deg(x^{\beta})>deg(x^{\gamma})$ and so $\beta_{i}>\gamma_{i}$ for some $i$. Choose a maximal $q<s$ such that $\beta_{iq}>\gamma_{iq}$ for some $i$. Thus $\beta_{ij}=\gamma_{ij}$ for $j>q$. It follows $\beta_{iq}\geq \gamma_{iq}$ since $(\beta_{is},\ldots,\beta_{i0})\geq_{lex} (\gamma_{is},\ldots,\gamma_{i0})$. If $q\leq t$ then we have
\[ \alpha_{q} = \sum_{i=1}^{n}\gamma_{iq} < \sum_{i=1}^{n}\beta_{iq} = \sum_{i=1}^{n}\lambda_{iq}\leq \alpha_{q}, \]
which is not possible. It follows $q<t$ and so $\beta_{it} = \gamma_{it}$ for each $i$. But this is not possible because we get $\alpha_{t} = \sum_{i=1}^{n} \gamma_{it} = \sum_{i=1}^{n} \lambda_{it} = \alpha_{t}-1$. Hence $x^{\beta}\notin I+J_{t'}$.

Suppose now $t'>t$. If $e_{t'}>e_{t}$, then $G(J_{t})\cap G(J_{t'}) = \emptyset$ by degree reason. Assume $e_{t}=e_{t'}$. If follows $n=2$ by the previous remark. If $x_{1}^{\beta_{1}}x_{2}^{\beta_{2}}\in G(J_{t})\cap J_{t'}$ we necessarily get
$x_{1}^{\beta_{1}}x_{2}^{\beta_{2}}\in G(J_{t})\cap G(J_{t'})$ again by degree reason. But this is not possible since it implies that $\alpha_{t'}-1 = \beta_{1t'} + \beta_{2t'} = \alpha_{t'}$.

2. and 3. follows from 1.
\end{proof}

\begin{teor}
Let $u = \prod_{q=1}^{r}x_{i_{1}}^{\alpha_{q}}$, where $2\leq i_{1}<i_{2}<\cdots <i_{r}\leq n$. Let
\[ I=<u>_{\de} = \prod_{q=1}^{r}\prod_{j=0}^{s}(\me_{q}^{[d_{j}]})^{\alpha_{qj}},\] where $\alpha_{q} = \sum_{j=0}^{s}\alpha_{qj}d_{j}$. Suppose $i_{r}=n$. Let $1\leq a \leq r$ be an integer and
\[ P_{a}(I):=\{ (\lambda,t)\in \mathbb N^{a}\times \mathbb N^{a}|
\;1\leq \lambda_{1}<\cdots < \lambda_{a}=r,
t_{a}>\cdots >t_{1}, \alpha_{\lambda_{\nu} t_{\nu}}\neq 0,\; for\; 1\leq \nu \leq a \}.  \]
Let $J = \sum_{a=1}^{r}\sum_{(\lambda,t)\in P_{a}(I)}J_{(\lambda,t)}$, where $J_{(\lambda,t)}$ is the ideal
\[ \prod_{e=1}^{a}(x_{i_{\lambda_{e}}} \cdots x_{i_{\lambda_{e-1}}+1})^{d_{t_{e}}-1}
\prod_{\nu=1}^{a}\me_{\lambda_{\nu}}^{[d_{t_{\nu}+1}]}\prod_{j>t_{\nu}} (\me_{\lambda_{\nu}}^{[d_{j}]})^{\alpha_{\lambda_{\nu}j}} (\me_{\lambda_{\nu}}^{[d_{t_{\nu}}]})^{\alpha_{\lambda_{\nu},t_{\nu}}-1} \prod_{q=\lambda_{\nu-1}+1}^{\lambda_{\nu}-1} \prod_{j\geq t_{\nu}} (m_{q}^{[d_{j}]})^{\alpha_{qj}}, \]
where we denote $\me^{[d_{t_{a+1}}]} = S$. Then $Soc(S/I) = (J+I)/I$.
\end{teor}

\begin{proof}
The proof will be given by induction on $r$, the case $r=1$ being done in Lemma $1.5.1$. Suppose that $r>1$. For $1\leq q \leq r$, let: $I_{q}=\prod_{e=1}^{q}\prod_{j=0}^{s}(\me_{e}^{[d_{j}]})^{\alpha_{ej}}$ and $S_{q}=k[x_{1},x_{2},\ldots,x_{i_{q}}]$
For $t$ with $\alpha_{rt}\neq 0$, denote:
\[ I^{(t)} = \me_{r-1}^{[d_{t}]}\prod_{j<t}(\me_{r-1}^{[d_{j}]})^{\alpha_{r-1,j}}I_{r-2}.\]
Let $J^{(t)}$ be an ideal in $S_{r-1}$ such that $Soc(S_{r-1}/I^{(t)}) = (J^{(t)}+I^{(t)})/I^{(t)}$. The induction step is given in the following lemma:
\begin{lema}
Suppose $i_{r}=n$ and let
\[ J = \sum_{t=0,\alpha_{rt}\neq 0} (x_{n}\cdots x_{i_{r-1}+1})^{d_{t}-1} \prod_{j>t}(\me_{r}^{[d_{j}]})^{\alpha_{r}j}
 \prod_{j\geq t}(\me_{r-1}^{[d_{j}]})^{\alpha_{r-1,j}}(\me_{r}^{[d_{t}]})^{\alpha_{rt}-1}J^{(t)}. \]
 Then $Soc(S/I) = (J+I)/I$.
\end{lema}
\begin{proof}
Let $w \in S\setminus I$ be a monomial such that $\me_{r}w \subset I$. As in the proof of lemma $1.5.1$, we choose for each $1\leq \rho \leq n$, $e_{\rho}= max\{e:\; x_{\rho}^{d_{e}-1}|w\}$. Renumbering variables $\{x_{n}, \ldots, x_{i_{r-1}+1} \}$ (it does not affect $I$, $J$ and $I^{(t)}$), we may suppose $e_{n}\leq e_{n-1}\leq \cdots \leq e_{i_{r-1}+1}$. Set $t=e_{n}$. We claim that $\alpha_{rt}\neq 0$. Indeed, if $\alpha_{rt}=0$ then from $x_{n}w\in I$ we get $x_{n}w/x_{n}^{d_{t}-1} \in \widetilde{I} = \prod_{j>t}(\me_{r}^{[d_{j}]})^{\alpha_{rj}}I_{r-1}$ because $x_{n}^{d_{t}-1} \in \prod_{j<t}(\me_{r}^{[d_{j}]})^{\alpha_{rj}}$. Since $t=e_{n}$ is maximal chosen, we get $w/x_{n}^{d_{t}-1} \in \widetilde{I}$ and so $w\in I$ a contradiction.

Reduction to the case that $x_{n}^{d_{t}}$ does not divide $w$. Suppose that
$w =x_{n}^{d_{t}}\widetilde{w}$ and set
\[ \widetilde{I} = (\me_{r}^{[d_{t}]})^{\alpha_{rt}-1} \prod_{\epsilon \leq 0,\epsilon\neq t}(\me_{r}^{[d_{\epsilon}]})^{\alpha_{r\epsilon}}I_{r-1}.\]
We see that $\me w\in I \Leftrightarrow \me \widetilde{w}\in \widetilde{I}$. Replacing $w$ and $I$ with $\widetilde{w}$ and $\widetilde{I}$, we reduce our problem to a new $\widetilde{t}<t$. The above argument implies that $\widetilde{\alpha}_{r\widetilde{t}} \neq 0$, where $\widetilde{\alpha}$ is the 'new' $\alpha$ of $\widetilde{I}$.

Reduction to the case when $\alpha_{rj} = \alpha_{r-1,j} = 0$ for $j>t$, $\alpha_{rt}=1$ and $\alpha_{r-1,t}=0$. From $x_{n}w\in I$, we see that there exists $\rho<n$ such that $x_{\rho}^{d_{j}}|w$ for $j>t$ if $\alpha_{rj}\neq 0$, or $j=t$ if $\alpha_{rt}>1$. Choose such maximal possible $\rho$. Set $w' = w/x_{\rho}^{d_{j}}$,
\[ I' = (\me_{r}^{[d_{j}]})^{\alpha_{rj}-1} \prod_{\epsilon\geq 0, \epsilon\neq j} (\me_{r}^{[d_{\epsilon}]})^{\alpha_{r\epsilon}}I_{r-1}. \]
We see that $\me w \subset I \Leftrightarrow \me w' \subset I'$, because from $x_{n}w\in I$, we get $x_{n}w'\in I'$ from the maximality of  $\rho$.

Let $\alpha'_{rj}=\alpha_{rj}-1$ and $\alpha'_{q\epsilon} = \alpha_{q\epsilon}$ for $(q,\epsilon) \neq (r,j)$. $\alpha'$ is the 'new' $\alpha$ for $I'$. If we show that
\[ w'\in J' = \sum_{e \geq 0, \alpha'_{re}\neq 0} (x_{n}\cdots x_{i_{r-1}+1})^{d_{e}-1} \prod_{\epsilon >e} (\me_{r}^{[d_{\epsilon}]})^{\alpha_{r}\epsilon}
 \prod_{j\geq e}(\me_{r-1}^{[d_{j}]})^{\alpha_{r-1,j}}(\me_{r}^{[d_{\epsilon}]})^{\alpha_{re}-1}J^{(t)}, \]
then $w = x_{\rho}^{d_{j}}w' \in \me_{r}^{[d_{j}]}J'\subset J$.
Using this procedure, by recurrence we arrive to the case $\alpha_{rj}=0$ for $j>t$ and $\alpha_{rt}=1$. Again from $x_{n}w \in I$, we note that there exists $\rho < i_{r-1}$ such that $x_{\rho}^{d_{j}}|w$ for $j\geq t$ with $\alpha_{r-1,j}\neq 0$. Choose such maximal possible $\rho$ and note that $\me w \subset I$ if and only if $\me w'' \in I''$ for $w'' = w/ x_{\rho}^{d_{j}}$, where \[ I'' = (\me_{r-1}^{d_{j}})^{\alpha_{r-1,j-1}} \prod_{\epsilon\geq 0,\;\epsilon\neq j} (\me_{r-1}^{[d_{\epsilon}]})^{\alpha_{r-1, \epsilon}} \prod_{\epsilon \geq 0} (\me_{r}^{[d_{\epsilon}]})^{\alpha_{r\epsilon}}I_{r-2}.\]
As above, we reduce our problem to $I''$ and the $\alpha''$, which is the new $\alpha$ of $I''$, is given by
$\alpha''_{r-1,j} = \alpha_{r-1,j-1}$, $\alpha''_{q\epsilon} = \alpha_{q\epsilon}$ for $(q,\epsilon) \neq (r-1,j)$.
Using this procedure, by recurrence we end our reduction.

Case $\alpha_{rj} = \alpha_{r-1j} = 0$ for $j>t$, $\alpha_{rt}=1$ , $\alpha_{r-1t}=0$ and $x_{n}^{d_{t}}$ does not divide $w$. Let express $w=(x_{n}\cdots x_{i_{r-1}+1})^{d_{t}-1}y$. We will show that $y$ does not depend on $\{x_{n},\ldots,x_{i_{r-1}+1}\}$. Indeed, if $n=i_{r-1}+1$ then there is nothing to show since $x_{n}^{d_{t}}$ does not divide $w$. Suppose that $n>i_{r-1}+1$, then from $x_{n}w\in I$ we get $y\in I_{r-1}$ because $x_{n-1}^{d_{t}-1}\in \prod_{j<t} (\me_{r}^{d_{j}})^{\alpha_{rj}}$ and the variables $x_{n}, \ldots, x_{i_{r-1}+1}$ are regular on $S/I_{r-1}S$. If $y=x_{\eta}y'$ for $\eta>i_{r-1}$, then as above $y' \in I_{r-1}$. Thus $w\in x_{\eta}^{d_{t}}x_{\rho}^{d_{t}-1}y' \subset I$ for any $\rho \neq \eta, i_{r-1} < \rho \leq n$, a contradiction.

Note that $\me_{r}w\in I \Rightarrow \me_{r-1}y \in I^{(t)}$ and so $w\in (x_{n}\cdots x_{i_{r-1}+1})^{d_{t}-1}J^{(t)}$. Since $\alpha_{rj} = \alpha_{r-1j} = 0$ for $j>t$ and $\alpha_{rt}=1$ and $\alpha_{r-1,t} = 0$, we get $w\in J$. Conversely, if $y\in J^{(t)}$, then it is clear that $w\in J$.
\end{proof}

We see by the above lemma that:
\[(*)\;J = \sum_{e\geq 0,\alpha_{re}\neq 0} (x_{n} \cdots x_{i_{r-1}+1})^{d_{e}-1} \prod_{j>e} (\me_{r-1}^{[d_{j}]})^{\alpha_{rj}} \prod_{j\geq e} (\me_{r-1}^{[d_{j}]})^{\alpha_{r-1,j}} (\me_{r}^{[d_{e}]})^{\alpha_{re}-1}J^{(e)}.\]
Since $\lambda_{a}=r-1$, by the induction hypothesis applied to $I^{(e)}$ we get:
\[ J^{(e)} = \sum_{a=1}^{r-1} [\sum_{(\lambda,t)\in P_{a}(I^{(e)}),t_{a}=e} \prod_{s=1}^{a} (x_{i_{\lambda_{s}}} \cdots x_{i_{\lambda_{s-1}+1}} )^{d_{t_{s}}-1} \cdot J'_{(\lambda,t)} + \]\[ +
\sum_{(\lambda,t)\in P_{a}(I^{(e)}),t_{a}<e} \prod_{s=1}^{a} (x_{i_{\lambda_{s}}} \cdots x_{i_{\lambda_{s-1}+1}} )^{d_{t_{s}}-1} \cdot J''_{(\lambda,t)} ],\;where\; \]
\[ J'_{(\lambda,t)} = \prod_{q=\lambda_{a-1}+1}^{\lambda_{a}-1}\prod_{j\geq e} (\me_{q}^{[d_{j}]})^{\alpha_{qj}}\widetilde{J}_{(\lambda,t)} \;and\;\]
\[ J''_{(\lambda,t)} = \me_{r-1}^{[d_{e}]} \prod_{j>t_{a}}^{e-1}(\me_{\lambda_{a}}^{[d_{j}]})^{\alpha_{\lambda_{a},j}} (\me_{\lambda_{a}}^{[d_{t_{a}}]})^{\alpha_{\lambda_{a},t_{a}-1}}\cdot \prod_{q=\lambda_{a-1}+1}^{\lambda_{a}-1} \prod_{j\geq t_{a}} (\me_{q}^{[d_{j}]})^{\alpha_{qj}}\widetilde{J}_{(\lambda,t)} ,\;and\]
\[ \widetilde{J}_{(\lambda,t)} = \prod_{\nu=1}^{a-1}\me_{\lambda_{\nu}}^{[d_{t_{\nu+1}}]}
\prod_{j>t_{\nu}}(\me_{\lambda_{\nu}}^{[d_{j}]})^{\alpha_{\lambda_{\nu},j}} (\me_{\lambda_{\nu}}^{[d_{t_{\nu}}]})^{\alpha_{\lambda_{\nu},t_{\nu}-1}}\cdot
 \prod_{q=\lambda_{\nu-1}+1}^{\lambda_{\nu}-1}\prod_{j\geq t_{\nu}} (\me_{q}^{[d_{j}]})^{\alpha_{qj}}.\]

If $t_{a}=e$, set $\lambda'_{\nu} = \lambda_{\nu}$ for $\nu<a$, $\lambda'_{a}=r$ and see that $(\lambda',t) \in P_{a}(I)$. If $t_{a}<e$, then put $\lambda''_{\nu} = \lambda_{\nu}$ for $\nu\leq a$, $\lambda''_{a+1}=r$, $t''_{\nu} = t_{\nu}$ for $\nu\leq a$ and $t''_{a+1}=e$ and then $(\lambda'',t) \in P_{a+1}(I)$. Substituting $J^{(e)}$ in $(*)$, we get the following expression for $J$:
\[\sum_{a=1}^{r-1} \sum_{(\lambda',t)\in P_{a}(I)} \prod_{\nu=1}^{a} (x_{i_{\lambda'_{\nu}}}\cdots x_{i_{\lambda'_{\nu-1}+1}})^{d_{t_{\nu}}-1} \cdot [ \prod_{j>e} (\me_{\lambda'_{a}}^{[d_{j}]})^{\alpha_{\lambda'_{a}j}} (\me_{\lambda'_{a}}^{[d_{e}]})^{\alpha_{\lambda'_{a}e}-1}
 \cdot \prod_{q=\lambda'_{a-1}+1}^{\lambda'_{a}-1}\prod_{j\geq e} (\me_{q}^{[d_{j}]})^{\alpha_{qj}}]\cdot \] \[ \cdot \widetilde{J}_{(\lambda,t)} +
\sum_{a=1}^{r-1} \sum_{(\lambda'',t'')\in P_{a+1}(I)} \prod_{\nu=1}^{a+1} (x_{i_{\lambda''_{\nu}}}\cdots x_{i_{\lambda''_{\nu-1}+1}})^{d_{t_{\nu}}-1} \cdot  [\prod_{j>e}(\me_{\lambda''_{a+1}}^{[d_{j}]})^{\alpha_{\lambda''_{a+1},j}}
(\me_{\lambda''_{a+1}}^{[d_{t''_{a}+1}]})^{\alpha_{\lambda''_{a+1}t''_{a+1}}-1}] \]\[
[ \me_{\lambda''_{a}}^{[d_{t''_{a+1}}]} \prod_{j\geq t''_{a}} (\me_{\lambda''_{a}}^{[d_{j}]})^{\alpha_{\lambda''_{a+1},j}}
(\me_{\lambda''_{a}}^{[d_{t''_{a}}]})^{\alpha_{\lambda''_{a}t''_{a}}-1}
\prod_{q=\lambda''_{a-1}+1}^{\lambda''_{a}-1}\prod_{j\geq t''_{a}} (\me_{q}^{[d_{j}]})^{\alpha_{qj}}] \cdot \widetilde{J}_{(\lambda,t)}. \]
Since all the pairs of $P_{b}(I)$ have the form $(\lambda',t)$ or $(\lambda'',t'')$ for a pair $(\lambda,t)\in P_{b}(I)$ or $(\lambda,t)\in P_{b-1}(I)$ respectively, it is not hard to see that the expression above is the formula of $J$ as stated.
\end{proof}

\noindent
Let $s_{q}=max\{j|\alpha_{qj}\neq 0\}$, $d_{qt} = \sum_{e=1}^{q}\sum_{j\geq t}^{s_{q}}\alpha_{ej}d_{j}$,
$D_{q}=d_{q,s_{q}} + (i_{q}-1)(d_{s_{q}}-1)$ for $1\leq q\leq r$.

\begin{cor}
With the notation and hypothesis of above theorem, for $(\lambda,t)\in P_{a}(I)$ let:
\[d_{(\lambda,t)} = \sum_{\nu=1}^{a} \sum_{q=\lambda_{\nu-1}+1}^{\lambda_{\nu}} \sum_{j\geq t_{\nu}} \alpha_{qj}d_{j}.\; Then:\]
\begin{enumerate}
    \item $Soc(I_{r-1}S/I) = Soc(S/I)$.
    \item $((J+I)/I)_{e}\neq 0$, if and only if
          $ e = d_{(\lambda,t)} + \sum_{\nu=1}^{a}(i_{\lambda_{\nu}} - i_{\lambda_{\nu-1}})(d_{t_{\nu}}-1) - d_{t_{1}}$,
          for some $1\leq a\leq r$ and $(\lambda,t)\in P_{a}(I)$.
    \item $c=max\{e|((J+I)/I)_{e}\neq 0  \} = d_{r,s_{r}} + (n-1)(d_{s_{r}} -1) - 1$.
\end{enumerate}
\end{cor}

\begin{proof}
1.Note that $J_{(\lambda,t)}$ is contained in
\[ \prod_{q=1,q\notin \{\lambda_{1},\ldots,\lambda_{q}\}}^{r} (\me_{q}^{[d_{j}]})^{\alpha_{qj}}
\prod_{\nu=1}^{a}[\prod_{j\neq t_{\nu}} (\me_{\lambda_{\nu}}^{d_{j}})^{\alpha_{\lambda_{\nu},j}}
(\me_{\lambda_{\nu}}^{t_{\nu}})^{\alpha_{\lambda_{\nu},t_{\nu}-1}}]\prod_{\epsilon=1}^{a-1}\me_{\lambda_{\epsilon+1}}^{d_{t_{\epsilon}}+1}.\]
Since $\me_{\lambda_{\epsilon}}^{d_{t_{\epsilon}+1}} \subset \me_{\lambda_{\epsilon}}^{d_{t_{\epsilon}}}$ for $t_{\epsilon+1}>t_{\epsilon}$ and $\lambda_{a}=r$ if follows that
\[ J \subset \prod_{j\neq t_{a}} (\me_{r}^{[d_{j}]})^{\alpha_{rj}} (\me_{r}^{[d_{t_{a}}]})^{\alpha_{rt_{a}}-1}I_{r-1},\]
as desired.

2.If $((J+I)/I)_{e}\neq 0$ then there exists a monomial $u\in J\setminus I$ of degree $e$. But $u\in J$, implies that there exists $a\in \{1,\ldots,r\}$ and $(\lambda,t)\in P_{a}(I)$ such that $u\in J_{(\lambda,t)}$. Thus the degree of $u$ is $e = d_{(\lambda,t)} + \sum_{\nu=1}^{a}(i_{\lambda_{\nu}} - i_{\lambda_{\nu-1}})(d_{t_{\nu}}-1) - d_{t_{1}}$, as required.

Conversely, let $e = d_{(\lambda,t)}+ \sum_{\nu=1}^{a}(i_{\lambda_{\nu}} - i_{\lambda_{\nu-1}})(d_{t_{\nu}}-1) - d_{t_{1}}$ for some $a\in \{1,\ldots,r\}$ and $(\lambda,t)\in P_{a}(I)$. We show that the monomial
\[w = \prod_{\nu=1}^{a} (x_{i_{\lambda_{\nu}}}\cdots x_{i_{\lambda_{\nu-1}+1}})^{d_{t_{\nu}}-1}\cdot x_{1}^{d_{(\lambda,t)}-d_{t_{1}}} \in J\setminus I. \]
Obvious $w\in J$. Let us assume that $w\notin I$. Then $w/x_{i_{\lambda_{a}}}^{d_{t_{a}}-1} \in \prod_{j\geq t_{a}} (\me_{\lambda_{a}}^{[d_{j}]})^{\alpha_{\lambda_{a}j}}I_{\lambda_{a}-1}$ because $x_{i_{\lambda_{a}}}^{d_{t_{a}}-1} \in \prod_{j<t_{a}}(\me_{\lambda_{a}}^{[d_{j}]})^{\alpha_{\lambda_{a}j}}$ and $x_{i_{\lambda_{a}}} \notin \me_{j}$ for $j<\lambda_{a}$. Inductively we get that:
\[ w/(x_{i_{\lambda_{a}}} \cdots x_{i_{\lambda_{a-1}+1}})^{d_{t_{a}}-1} \in \prod_{q=\lambda_{a-1}+1}^{\lambda_{a}}  \prod_{j\geq t_{a}} (\me_{\lambda_{a}}^{[d_{j}]})^{\alpha_{qj}}I_{\lambda_{a}-1}. \]
Following the same reduction and using that $t_{a}>\cdots>t_{1}$ we obtain that:
\[ x_{1}^{d_{(\lambda,t)}-d_{t_{1}}} \in \prod_{\nu=1}^{a} \prod_{q=\lambda_{\nu-1}+1}^{\lambda_{\nu}}  \prod_{j\geq t_{\nu}} (\me_{\lambda_{a}}^{[d_{j}]})^{\alpha_{qj}}. \]
So $d_{(\lambda,t)}-d_{t_{1}} \geq d_{(\lambda,t)}$, a contradiction.

3.Note that $c=d_{(\lambda',t')}$ for $(\lambda',t')\in P_{1}(I)$ with $\lambda' = \lambda_{1} = r$ and $t'=t_{1}=s_{r}$. We have to show that:
\[ c = d_{r,s_{r}} + (n-1)(d_{s_{r}} -1) - 1 \leq d_{(\lambda,t)} + \sum_{\nu=1}^{a}(i_{\lambda_{\nu}} - i_{\lambda_{\nu-1}})(d_{t_{\nu}}-1) - d_{t_{1}}, \]
for any $1\leq a\leq r$ and $(\lambda,t)\in P_{a}(I)$. Since $d_{s_{r}}-1 \leq (d_{t_{\nu}}-1) + \sum_{j\leq t_{\nu}}^{s_{r}-1}\alpha_{qj}d_{j}$ for all $q$ with $i_{\nu-1}<q \leq i_{\nu}$, we see that:
\[ d_{r,s_{r}} + (n-1)(d_{s_{r}} -1)-1 \geq d_{(\lambda,t)} + \sum_{\nu=2}^{a}(i_{\lambda_{\nu}} - i_{\lambda_{\nu-1}})(d_{t_{\nu}}-1) + (i_{\lambda_{1}}-1)(d_{t_{1}}-1) -1. \]
On the other hand, $(i_{\lambda_{1}}- i_{\lambda_{0}})(d_{t_{1}}-1) = (i_{\lambda_{1}}- 1)(d_{t_{1}}-1) + d_{t_{1}}-1$,
and replacing that in the above relation we obtained what is required.
\end{proof}

\begin{exm}{\em
Let $\de: 1|2|4|12$.
\begin{enumerate}
    \item Let $u=x_{3}^{21}$. We have $\alpha_{0}=1$, $\alpha_{1}=0$,  $\alpha_{2}=2$ and $\alpha_{3}=1$ so:
    \[ I = <u>_{\de} = (x_{1},x_{2},x_{3})(x_{1}^{4},x_{2}^{4},x_{3}^{4})^{2}(x_{1}^{12},x_{2}^{12},x_{3}^{12}).\]
    Let $J= \sum_{t=0,\alpha_{t}>0} J_{t}$, where $J_{t} = (x_{1}x_{2}x_{3})^{d_{t}-1}(x_{1}^{d_{t}},
                x_{2}^{d_{t}}, x_{3}^{d_{t}})^{\alpha_{t}-1}\prod_{j>t}(x_{1}^{d_{j}},
                x_{2}^{d_{j}}, x_{3}^{d_{j}})^{\alpha_{j}}$.

    $J_{0} = (x_{1}x_{2}x_{3})^{1-1}\cdot (x_{1},x_{2},x_{3})^{1-1} \cdot \prod_{j>t}(x_{1}^{d_{j}},
                x_{2}^{d_{j}}, x_{3}^{d_{j}})^{\alpha_{j}} = (x_{1}^{4},x_{2}^{4},x_{3}^{4})^{2}(x_{1}^{12},x_{2}^{12},x_{3}^{12}).$

    $J_{2} =(x_{1}x_{2}x_{3})^{4-1}(x_{1}^{4},x_{2}^{4},x_{3}^{4})^{2-1}(x_{1}^{12},x_{2}^{12},x_{3}^{12}) = (x_{1}x_{2}x_{3})^{3}(x_{1}^{4},x_{2}^{4},x_{3}^{4})(x_{1}^{4},x_{2}^{4},x_{3}^{12})$ and

    $ J_{3}=(x_{1}x_{2}x_{3})^{12-1} = (x_{1}x_{2}x_{3})^{11}.$ From $1.5.1$ , $Soc(S/I) = (J+I)/I$.

    \item Let $u=x_{2}^{9}x_{3}^{16}$. We have $r=2$, $i_{1}=2$ and $i_{2}=3$. Also $\alpha_{10}=1$, $\alpha_{12}=2$,
          $\alpha_{22}=1$, $\alpha_{23}=1$ and the other components of $\alpha$ are zero. Then
    \[ I = <u>_{\de} = <x_{2}^{9}>_{\de}<x_{3}^{16}>_{\de} = (x_{1},x_{2})(x_{1}^{4},x_{2}^{4})^{2}
                    (x_{1}^{4},x_{2}^{4},x_{3}^{4})(x_{1}^{12},x_{2}^{12},x_{3}^{12}).\]
    We have two possible partitions: (a) $(2)$ and (b) $(1<2)$.

    (a)$\lambda=\lambda_{1}=2$, $t=t_{1}$ such that $\alpha_{2t}\neq 0$. We have two possible $t$: $t=2$ or $t=3$.

    (i)For $t=2$ we obtain (according to the Theorem $1.5.4$) the following part of the socle:
    \[ J_{(2,2)} = (x_{1}x_{2}x_{3})^{3} (x_{1}^{12}, x_{2}^{12}, x_{3}^{12}) (x_{1}^{4},x_{2}^{4})^{4} \]

  (ii)For $t=3$ we obtain:
    \[ J_{(2,3)} = (x_{1}x_{2}x_{3})^{11}  \]

  (b)$1=\lambda_{1}<\lambda_{2}=2$, $t=(t_{1},t_{2})$ such that $\alpha_{\lambda_{e},t_{e}}\neq 0$ for $1\leq e\leq 2$
  and $t_{1}<t_{2}$. According to our expressions for $\alpha_{i}$ we have three possible cases: $t_{1}=0,t_{2}=2$ or
  $t_{1}=0,t_{2}=3$ or $t_{1}=2,t_{2}=3$.

  (i)For $t_{1}=0$ and $t_{2}=2$ we obtain:
  \[ J_{(1,2),(0,2)} = x_{3}^{3}(x_{1}^{4},x_{2}^{4})(x_{1}^{4},x_{2}^{4})^{2}(x_{1}^{12},x_{2}^{12},x_{3}^{12}). \]

  (ii)For $t_{1}=0$ and $t_{2}=3$ we obtain:
  \[ J_{(1,2),(0,3)} = x_{3}^{11}(x_{1}^{12},x_{2}^{12})(x_{1}^{4},x_{2}^{4})^{2} \]

  (iii)For $t_{1}=2$ and $t_{2}=3$ we obtain:
  \[ J_{(1,2),(2,3)} =  x_{1}^{3}x_{2}^{3}x_{3}^{11}(x_{1}^{12},x_{2}^{12})(x_{1}^{4},x_{2}^{4})\]

  From $2.4$ it follows that if $J = J_{(2,2)} + J_{(2,3)} + J_{(1,2),(0,2)} + J_{(1,2),(0,3)} + J_{(1,2),(2,3)}$ then
  $Soc(S/I)=(I+J)/J$.
\end{enumerate}
}\end{exm}

\section{A generalization of Pardue's formula.}

In this section, we give a generalization of a theorem proved by Aramova-Herzog \cite{ah} and Herzog-Popescu \cite{hp} which is known as "Pardue's formula".

Let $1\leq i_{1}<i_{2}<\cdots<i_{r}=n$ and let $\alpha_{1},\ldots,\alpha_{r}$ be some positive integers.
Let $u=\prod_{i=1}^{r}x_{i_{q}}^{\alpha_{q}} \in S=K[x_{1},\ldots,x_{n}]$. Let $I = <u>_{\de}$ the principal
$\de$-fixed ideal generated by $u$. From Proposition $1.4.11$ it follows that 
$I = \prod_{r=1}^{q}\prod_{j=0}^{s}(\me_{q}^{[d_{j}]})^{\alpha_{qj}}$, 
where $\alpha_{q}=\sum_{j=0}^{s}\alpha_{qj}d_{j}$. If $i_{1}=1$, it follows that $I=x_{1}^{\alpha_{1}}I'$, where
$I' = \prod_{r=2}^{q}\prod_{j=0}^{s}(\me_{q}^{[d_{j}]})^{\alpha_{qj}}$, and therefore $reg(I) = \alpha_{1} + reg(I')$. Thus, we may assume $i_{1}\geq 2$.

If $N$ is a graded $S$-module of finite length, we denote $s(N)=max\{i|N_{i}\neq 0\}$. Let $s_{q}=max\{j|\alpha_{qj}\neq 0\}$ and $d_{qt}= \sum_{e=1}^{q}\sum_{j\geq t}^{s_{e}}\alpha_{ej}d_{j}$.
Let $D_{q} = d_{qs_{q}} + (i_{q}-1)(d_{s_{q}}-1)$, for $1\leq q \leq r$. With this notations we have:

\begin{teor}
$reg(I) = max_{1\leq q \leq r} D_{q}$. In particular, if $I = <x_{n}^{\alpha}>_{\de}$ and $\alpha=\sum_{t=0}^{s} \alpha_{t}d_{t}$ with $\alpha_{s}\neq 0$ then $reg(I) = \alpha_{s}d_{s} + (n-1)(d_{s}-1)$.
\end{teor}

\begin{proof}
Let $I_{\ell}=\prod_{q=1}^{r-\ell}\prod_{j=0}^{s}(\me_{q}^{[d_{j}]})^{\alpha_{qj}}$, for $0 \leq \ell \leq r$.
Then $I=I_{0}\subset I_{1}\subset \cdots \subset I_{r} = S$ is the sequential chain of ideals of $I$, i.e.
$I_{\ell+1} = (I_{\ell}:x_{n_{\ell}}^{\infty})$, where $n_{\ell} = i_{r-\ell}$. 
Let $S_{\ell}=k[x_{1},\ldots,x_{n_{\ell}}]$ and $m_{\ell} = (x_{1},\ldots,x_{n_{\ell}})$.
Let $J_{\ell}\subset S_{\ell}$ be the ideal generated by $G(I_{\ell})$.

The Corollary $1.5.6$ implies that $c_{e}=D_{e}-1$ is the maximal degree for a nonzero element of $Soc(S_{\ell}/J_{\ell})$. Proposition $1.2.2$ implies $reg(I) = max\{s(J_{\ell}^{sat}/J_{\ell})\;|\; \ell=0,\ldots,r-1)\} + 1$. Also, from the Corollary $1.5.6$, we get
\[ s(J_{\ell}^{sat}/J_{\ell}) = s(Soc(J_{\ell}^{sat}/J_{\ell})) = s(Soc(S_{\ell}/J_{\ell})) = D_e -1,\]
which complete the proof. Indeed, for the first equality, if $u\in J_{\ell}^{sat}\setminus J_{\ell}$ and 
$deg(u)=s(J_{\ell}^{sat}/J_{\ell})$ it follows that $u\in Soc(J_{\ell}^{sat}/J_{\ell})$, since $\me u\subset J_{\ell}$
by degree reasons.
\end{proof}

\begin{cor}
$reg(I) \leq n\cdot (deg(u)-1) + 1$.
\end{cor}

\begin{cor}
$S/I$ has at most $r$-corners among $(i_{q},D_{q}-1)$ for $1\leq q\leq r$. If $i_{1}=1$ we replace $(i_{1},D_{1}-1)$ with $(1,\alpha_{1})$. The corresponding extremal Betti numbers are $\beta_{i_{q},D_{q}+i_{q}-1}$.
\end{cor}

\begin{proof}
By Theorem $1.1.6$ combined with the proof of Theorem $1.6.1$, $S/I$ has at most $r$-corners among
$(n_{\ell},s(I_{\ell+1}S_{\ell}/I_{\ell}S_{\ell}))$ and is enough to apply Corollary $1.5.6$.
\end{proof}

\begin{exm}{\em
Let $\de: 1|2|4|12$.
\begin{enumerate}
    \item Let $u=x_{3}^{21} \in k[x_{1},x_{2},x_{3}]$. We have $21=1\cdot 1 + 0\cdot 2 + 2 \cdot 4 + 1 \cdot 12$. From
    $3.1$, we get: \[ reg(<u>_{\de}) = 1\cdot 12 + (3-1)\cdot (12 - 1) = 34.\]
    \item Let $u=x_{1}^{2}x_{2}^{16}x_{3}^{9}$. Then $reg(<u>_{\de}) = 2 + reg(<u'>_{\de})$, where $u' = u/x_{1}^{2}$.
          We compute $reg(<u'>_{\de})$. With the notations above, we have $i_{1}=2$, $i_{2}=3$, $r=2$, $\alpha_{1}=16$ and
          $\alpha_{2}=9$. We have $\alpha_{1} = 1\cdot 4 + 1\cdot 12$ and $\alpha_{2} =1\cdot 1 + 2\cdot 4$, thus
          $s_{1}=3$ and $s_{2}=2$. $D_{1} = d_{13} + (2-1)(d_{3}-1) = 12 + 11 = 23$ and
          $D_{2} = d_{22} + (3-1)(d_{2}-1) = 24 + 6 = 30$. In conclusion, $reg(<u>_{\de}) = 2 + max\{23,30\} = 32$.
    \item Let $u=x_a^{\alpha}x_b^{\beta}$, with $1<a<b\leq n$ and
$\beta<\alpha$, $\beta|\alpha$ be two integers. Set $d_1=\beta$,
$d_2=\alpha$ and let $I$ be the principal $\de$-fixed ideal generated by
$u$. Obviously $I=SBT(u)=(x_1^{\alpha},\ldots,
x_a^{\alpha})(x_1^{\beta},\ldots,x_b^{\beta})$. We have $i_1=a$,
$i_2=b$, $s_1=2$, $s_2=1$, $d_{1s_1}=\alpha$,
$d_{2s_2}=\alpha+\beta$,
$D_1=\alpha+(a-1)(\alpha-1)=a(\alpha-1)+1=\chi_1+1$,
$D_2=\alpha+\beta+(b-1)(\beta-1)=\alpha+b(\beta-1)+1=\chi_2+1$ in
the notations of Theorem $1.3.6$. Note that $1.3.6$ and $1.6.1$ give
the same regularity of $I$, namely $max\{D_1,D_2\}$.	      
\end{enumerate}
}\end{exm}

\begin{dfn}
We say that a monomial ideal $I\subset S$ is a \Dei, if $I$ is a sum of $\de$-fixed ideals, 
for various $\de$-sequences.
\end{dfn}

Since any $\de$-fixed ideal is a Borel type ideal and a sum of Borel type ideals is stil a Borel type ideal, it follows
that any \Dei $I$ is a Borel type ideal. Therefore, from Corollary $1.2.11$ we get the next:

\begin{cor}
If $I$ is a \Di then $reg(I)=min\{e:\;e\geq deg(I), I_{\geq e}\;is\; stable \}$.
\end{cor}

We mention that this result was first obtained as a consequence of the proof of Pardue's formula by Herzog-Popescu \cite{hp} in the special case of a principal $p$-Borel ideal.

\section{\dis generated by powers of variables.}

Firstly, let fix some notations. Let $u_{1},\ldots,u_{m}\in S$ be some monomials. We say that $I$ is the \di generated by $u_{1},\ldots,u_{m}$, if $I$ is the smallest \di, w.r.t inclusion, which contained $u_{1},\ldots,u_{m}$, and we write $I=<u_{1},\ldots,u_{m}>_{\de}$. In particular, if $m=1$, we say that $I$ is the principal \di generated by $u=u_{1}$ and we write $I=<u>_{\de}$. 

\begin{lema}
If $1\leq j\leq j'\leq n$ and $\alpha\geq \beta$ are positive integer, then $<x_{j}^{\alpha}>\subset <x_{j'}^{\beta}>$.
\end{lema}

\begin{proof}
Indeed, using Proposition $1.4.9$ it is enough to notice that $<x_{j}^{\alpha}> \subset <x_{j'}^{\alpha}>$ which is true because $x_{j}^{\alpha} \in <x_{j'}^{\alpha}>$.
\end{proof}

Our next goal is to give the minimal set of generators for a \di generated by some powers of variables. 
Using the previous lemma, we had reduced to the next case:

\begin{prop}
Let $n\geq 2$ and let $1\leq i_{1} < i_{2} < \cdots < i_{r}=n$ be some integers. Let $\alpha_{1}< \alpha_{2}< \cdots < \alpha_{r}$ some positive integers. Then
\[ I=<x_{i_{1}}^{\alpha_{1}},x_{i_{2}}^{\alpha_{2}},\ldots,x_{i_{r}}^{\alpha_{r}}>_{\de} = \sum_{q=1}^{r} I^{(q)},
with \; I^{(q)} =
   \sum_{ \footnotesize \begin{array}{c} 0\leq \gamma_{1},\ldots,\gamma_{q}\leq_{\de}\alpha_{q},\\ \;\gamma_{1}+\cdots+\gamma_{i}<\alpha_{i},\;for\;i<q\\ 
   \gamma_{1}+\cdots+\gamma_{i}<_{d}\alpha_{q},\;for\;i<q\\
   \gamma_{1}+\cdots+\gamma_{q}=\alpha_{q} \end{array} \normalsize }   
   \prod_{e=1}^{q}\prod_{t=0}^{s} (\mathbf{n}_{e}^{[d_{t}]})^{\gamma_{et}},\]
   where $\mathbf{n}_{e} = \{x_{i_{e-1}+1},\ldots,x_{i_{e}}\}$, $\mathbf{n}_{e}^{[d_{t}]} = \{x_{i_{e-1}+1}^{d_{t}}, \ldots,x_{i_{e}}^{d_{t}}\}$, $i_{0}=0$ and $\gamma_{e}=\sum_{t=0}^{s}\gamma_{et}$.
\end{prop}

\begin{proof}
Let $\mathbf{m}_{q} = \{x_{1},\ldots,x_{i_{q}}\}$ for $1\leq q\leq r$. Obviously, $\en_{q}=\me_{q}\setminus \me_{q-1}$ for $q>0$ and $\en_{1}=\me_{1}$. Using the obvious fact that $I$ is the sum of principal \dis generated by the $\de$-generators of $I$ together with Proposition $1.4.8$ we get:
\[ I 
= \sum_{q=1}^{r} \prod_{t=0}^{s} (\me_{q}^{[d_{t}]})^{\alpha_{qt}},\;where\;\alpha_{q}=\sum_{t=0}^{s}\alpha_{qt}d_{t}.\]
Denote $S_{q}=K[x_{1},\ldots,x_{i_{q}}]$ for $1\leq q\leq r$. In order to obtain the required formula, we use induction on $r\geq 1$, the case $r=1$ being obvious. Let $r>1$ and
assume that the assertion is true for $r-1$, i.e 
\[ I'=<x_{i_{1}}^{\alpha_{1}},\ldots,x_{i_{r-1}}^{\alpha_{r-1}}>_{\de} = \sum_{q=1}^{r-1}  
   \sum_{\footnotesize \begin{array}{c} 0\leq \gamma_{1},\ldots,\gamma_{q}\leq_{\de}\alpha_{q},\\ \;\gamma_{1}+\cdots+\gamma_{i}<\alpha_{i},\;for\;i<q\\
   \gamma_{1}+\cdots+\gamma_{i}<_{d}\alpha_{q},\;for\;i<q\\
   \gamma_{1}+\cdots+\gamma_{q}=\alpha_{q} \end{array} \normalsize} \prod_{e=1}^{q}\prod_{t=0}^{s} (\mathbf{n}_{e}^{[d_{t}]})^{\gamma_{et}}  \subset S_{r-1}. \]
Obviously, $I=I'S + <x_{n}^{\alpha_{r}}>_{\de} = I'S + \prod_{t=0}^{s}(\me_{r}^{[d_{t}]})^{\alpha_{rt}}$. Also, $I'S$ 
and $I'$ have the same set of minimal generators and none of the minimal generators of $I'S$ is in $I^{(r)}$,
because of degree.

But, a minimal generator of $<x_{n}^{\alpha_{r}}>_{\de}$ is of the form $w=\prod_{t=0}^{s}\prod_{j=1}^{n}x_{j}^{\lambda_{tj}d_{t}}$ with $0\leq \lambda_{tj}$ and $\sum_{j=1}^{n}\lambda_{tj}=\alpha_{rt}$. Suppose $w\notin I'S$. In order to complete the proof, we will 
show that $w\in I^{(r)}$. Let $v_{q}=\prod_{t=0}^{s}\prod_{j=i_{q-1}+1}^{i_{q}}x_{j}^{\lambda_{tj}d_{t}}$ and
let $w_{q}=\prod_{e=1}^{q}v_{e}$. Obvious, $w = v_{1}\cdots v_{r} = w_{r}$. Since $w\notin I'$ it follows that
$w_{q}\notin I^{(q)}$ for any $1\leq q\leq r-1$. But $w_{q}\notin I^{(q)}$ implies 
$(*)\;\sum_{t=0}^{s}\sum_{j=1}^{i_{q}} \lambda_{tj}d_{t} < \alpha_{q}$, otherwise $w_{q} \in <x_{i_{q}}^{\alpha_{q}}S_{q}>_{\de}S_{r-1} \subset I'$ and thus $w\in I'$, a contradiction.
We choose $\gamma_{e} = \sum_{t=0}^{s}\sum_{j=i_{e-1}+1}^{i_{e}}\lambda_{tj}d_{t}$ for $1\leq e\leq r$.
For $1\leq q<r$, $(*)$ implies $\gamma_{1}+\cdots +\gamma_{q}<\alpha_{q}$. On the other hand, it is obvious that 
$\gamma_{1}+\cdots + \gamma_{e}\leq_{d}\alpha_{r}$ for any $1\leq e\leq r$ and $\gamma_{1}+\cdots+\gamma_{r}=\alpha_{r}$. Thus $w\in I^{(r)}$ as required.
\end{proof}

\begin{exm}{\em
Let $\de:1|2|4|12$ and let $I=<x_{2}^{7},x_{3}^{10},x_{5}^{17}> \subset K[x_{1},\ldots,x_{5}]$. We have $7=1\cdot 1 + 1\cdot 2 + 1\cdot 4$, $10=1\cdot 2 + 2\cdot 4$, $17=1\cdot 1 + 1\cdot 4 + 1\cdot 12$. 
We have 

\[ I^{(1)}= <x_2^{7}>_{\de} = (x_{1},x_{2})(x_{1}^{2},x_{2}^{2})(x_{1}^{4},x_{2}^{4}).\]

In order to compute $I^{(2)}$, we need to find all the pairs $(\gamma_1,\gamma_2)$ such that $\gamma_1<7$,
$\gamma_1<_{\de}10$ and $\gamma_2=10-\gamma_1$. We have $4$ pairs, namely $(0,10)$, $(2,8)$, $(4,6)$ and $(6,4)$, thus

\[ I^{(2)}=(x_{1}^{2},x_{2}^{2})(x_{1}^{4},x_{2}^{4})x_{3}^{4} + (x_{1}^{4},x_{2}^{4})x_{3}^{6} + 
	      (x_{1}^{2},x_{2}^{2})x_{3}^{8} + (x_{3}^{10}).\]

In order to compute	$I^{(3)}$, we need to find all $(\gamma_1,\gamma_2,\gamma_3)$ such that $\gamma_1<7$,
$\gamma_1+\gamma_2<10$, $\gamma_1<_{\de}17$, $\gamma_1+\gamma_2<_{\de}17$ and $\gamma_3 = 17-\gamma_1+\gamma_2$.
If $\gamma_1=0$ then, the pair $(\gamma_2,\gamma_3)$ is one of the following:$(0,17)$,$(1,16)$,$(4,13)$ or $(5,12)$.
If $\gamma_1=1$ then, the pair $(\gamma_2,\gamma_3)$ is one of the following:$(0,16)$ of $(4,12)$.
If $\gamma_1=4$ then, the pair $(\gamma_2,\gamma_3)$ is one of the following:$(0,13)$ of $(1,12)$. If $\gamma_1=5$ then, the pair $(\gamma_2,\gamma_3)$ is $(0,12)$. Thus 
\[I^{(3)} = (x_{1},x_{2})(x_{1}^{4},x_{2}^{4})(x_{4}^{12},x_{5}^{12}) +
	      (x_{1}^{4},x_{2}^{4})x_{3}(x_{4}^{12},x_{5}^{12}) + (x_{1}^{4},x_{2}^{4})(x_{4},x_{5})(x_{4}^{12},x_{5}^{12})+ \]\[ + (x_{1},x_{2})x_{3}^{4}(x_{4}^{12},x_{5}^{12}) + (x_{1},x_{2})(x_{4}^{4},x_{5}^{4})(x_{4}^{12},x_{5}^{12}) +
 x_{3}(x_{4}^{4},x_{5}^{4})(x_{4}^{12},x_{5}^{12}) + \]\[ + x_{3}^{4}(x_4,x_5)(x_{4}^{12},x_{5}^{12}) +
 x_{3}^{5}(x_{4}^{12},x_{5}^{12}) + (x_{4},x_{5})(x_{4}^{4},x_{5}^{4})(x_{4}^{12},x_{5}^{12}).\]
From Proposition $1.7.2$, we get $I=I^{(1)}+I^{(2)}+I^{(3)}$.
}\end{exm}

\begin{obs}{\em
For any $1\leq q\leq r$ and any nonnegative integers $\gamma_{1},\ldots,\gamma_{q} \leq_{\de} \alpha_{q}$ such that
\linebreak
$\gamma_{1}+\cdots +\gamma_{i}<\alpha_{i}$, $\gamma_{1}+\cdots +\gamma_{i}<_{\de}\alpha_{q}$ for $1\leq i<q$ and $\gamma_{1}+\cdots +\gamma_{q}=\alpha_{q}$ we denote
\[I^{(q)}_{\gamma_{1},\ldots,\gamma_{q}} = \prod_{e=1}^{q}\prod_{t=0}^{s} (\mathbf{n}_{e}^{[d_{t}]})^{\gamma_{et}}.
Proposition\; 1.7.2\; implies:\;
I = \sum_{q=1}^{r} \sum_{ \gamma_{1},\ldots,\gamma_{q}}I^{(q)}_{\gamma_{1},\ldots,\gamma_{q}}.\]
Let $\me=(x_{1},\ldots,x_{n})\subset S$ the irrelevant ideal of $S$. We have:
\[ (I:_{S}\me) = \bigcap_{j=1}^{n}(I:x_{j}) = \bigcap_{j=1}^{n}( (\sum_{q=1}^{r} \sum_{ \gamma_{1},\ldots,\gamma_{q}}I^{(q)}_{\gamma_{1},\ldots,\gamma_{q}}) : x_{j}) =
   \bigcap_{j=1}^{n}( \sum_{q=1}^{r} \sum_{ \gamma_{1},\ldots,\gamma_{q}} (I^{(q)}_{\gamma_{1},\ldots,\gamma_{q}} : x_{j})). \]
On the other hand, if $x_{j}\in \mathbf{n}_{p}$ for some $1\leq p\leq q$ then 
\[ J^{(q),j}_{\gamma_{1},\ldots,\gamma_{q}} := (I^{(q)}_{\gamma_{1},\ldots,\gamma_{q}} : x_{j}) = \prod_{e\neq p}^{q}\prod_{t=0}^{s} (\mathbf{n}_{e}^{[d_{t}]})^{\gamma_{et}}\mathbf{n_{p,\hat{j}}}^{[d_t]} (\mathbf{n}_p^{[d_t]})^{\gamma_{pt}-1}(\sum_{\gamma_{pt}>0}\prod_{j\neq t}(\mathbf{n}_{e}^{[d_{t}]})^{\gamma_{jt}}), \]
where $\mathbf{n_{p,\hat{j}}}^{[d_t]}= (x_{i_{p-1}+1}^{d_t},\ldots,x_{j}^{d_t-1},\ldots,x_{i_{p}}^{d_t})$
and $\mathbf{n_{p,\hat{j}}}^{[d_t]} (\mathbf{n}^{[d_t]})^{\gamma_{pt}-1}:=S$ if $\gamma_{pt}=0$. Thus
\[ (I:_{S}\me) = \left( \sum_{q^{1}=1}^{r} \sum_{ \gamma_{1}^{1},\ldots,\gamma_{q}^{1}} \right)
\cdots \left( \sum_{q^{n}=1}^{r} \sum_{ \gamma_{1}^{n},\ldots,\gamma_{q}^{n}} \right) \bigcap_{j=1}^{n} J^{(q^{j}),j}_{\gamma_{1}^{j},\ldots,\gamma_{q}^{j}}, \]
where for a given $q^{j}$, we have that $\gamma_{1}^{j},\ldots,\gamma_{q}^{j}\leq_{\de} \alpha_{q}$ are some 
nonnegative integers such that $\gamma_{1}^{j}+\cdots +\gamma_{i}^{j}<\alpha_{i}$, $\gamma_{1}^{j}+\cdots +\gamma_{i}^{j}<_{\de}\alpha_{q}$ for $1\leq i<q^{j}$ and $\gamma_{1}^{j}+\cdots +\gamma_{q}^{j}=\alpha_{q}$.
}\end{obs}

\begin{prop}
Let $n\geq 2$ and let $1\leq i_{1} < i_{2} < \cdots < i_{r}=n$ some integers. Let $\alpha_{1}< \alpha_{2}< \cdots < \alpha_{r}$ some positive integers. Let $s_q=max\{t|\; \alpha_{qt}>0\}$ for any $1\leq q\leq r$. We consider the ideal
$I=\sum_{q=1}^{r}I_{q}$, where $I_q = <x_{i_{q}}^{\alpha_{q}}>_{\de}$. Then, we have: $reg(I)\leq reg(I_r)$.
\end{prop}

\begin{proof}
Denote $e=max\{reg(I_1),\ldots,reg(I_r)\}$. Since a $\de$-fixed ideal is an ideal of Borel type, by Corollary $1.2.13$ it follows that $reg(I)\leq e$. On the other hand, if we denote $s_q=max\{t|\; \alpha_{qt}>0\}$ for any $1\leq q\leq r$, from Theorem $1.6.1$ we get $reg(I_q) = \alpha_{qs_q}d_{s_q} + (i_q - 1)(d_{s_q}-1)$, thus $max\{reg(I_1),\ldots,reg(I_r)\}=reg(I_r)$. In conclusion, $reg(I)\leq reg(I_r)$.
\end{proof}

\begin{prop}
With the above notations, for any $1\leq q\leq r$ we have:
\[ (I_q:\me_{q}) + (I_1 + \cdots +I_q) \subset ((I_1+\cdots+I_q):\me_{q}) \subset
((I_1+\cdots+I_q):\en_{q}) = (I_q:\en_{q}) + (I_1 + \cdots +I_q).\]
\end{prop}

\begin{proof}
Fix $1\leq q\leq r$. The first two inclusions are obvious. In order to prove the last equality, it is enough to
show that $((I_1+\cdots+I_q):x_j) \subset (I_q:x_j) + (I_1 + \cdots +I_q)$ for any $x_j\in \en_q$. Indeed, suppose 
$u\in ((I_1+\cdots+I_q):x_j)$, therefore $x_j\cdot u \in I_1+\cdots+I_q$. If $x_j\cdot u \notin I_q$ it follows that
$x_j\cdot u \in I_e$ for some $e<q$. Thus $u\in I_e$, since $x_j$ does not appear in a minimal generators
of $I_e$.
\end{proof}

\chapter{Generic initial ideal for complete intersections.}

\section{Main results.}

Let $S=K[x_1,\ldots,x_n]$ and let $I=(f_1,\ldots,f_n)\subset S$ be an ideal generated by a regular sequence of homogeneous polynomials. We say that a homogeneous polynomial $f$ of degree $d$ is \emph{semiregular} for $S/I$ if the maps $(S/I)_{t}\stackrel{\cdot f}{\longrightarrow} (S/I)_{t+d}$ are either injective, either surjective for all $t\geq 0$. We say that $S/I$ has the \emph{weak Lefschetz property} (WLP) if there exists a linear form $\ell\in S$, semiregular on $S/I$, in which case we say that $\ell$ is a weak Lefschetz element for $S/I$.
We say that $S/I$ has the \emph{strong Lefschetz property} $(SLP)$ if there exists a linear form $\ell\in S$ such that $\ell^{b}$ is semiregular on $S/I$ for all integer $b\geq 1$. In this case, we say that $\ell$ is a strong Lefschetz element for $S/I$.

We say that a property $(P)$ holds for a \emph{generic} sequence of homogeneous polynomials 
$f_1,f_2,\ldots, f_n \in S=K[x_1,x_2,\ldots,x_n]$ of given degrees $d_1,d_2,\ldots,d_n$ if there exists a nonempty open Zariski subset $U\subset S_{d_1}\times S_{d_2}\times \cdots \times S_{d_n}$ such that for every $(f_{1},f_{2},\ldots,f_{n})\in U$ the property (P) holds.

\vspace{5 pt}
\noindent
\textbf{Conjecture}.(Moreno) If $f_1,f_2,\ldots,f_n\in S=K[x_1,\ldots,x_n]$ is a generic sequence of homogeneous polynomials of given degrees $d_1,d_2,\ldots,d_n$, $I=(f_1,\ldots,f_n)$ and $J$ is the initial ideal of $I$ with respect to the revlex order, then $J$ is an almost revlex ideal, i.e. if $u\in J$ is a minimal generator of $J$ then every monomial of the same degree which preceeds $u$ must be in $J$ as well. 
\vspace{5 pt}

\begin{teo}
If $f_1,f_2,f_3\in S=K[x_1,x_2,x_3]$ is a regular sequence of homogeneous polynomials of given degrees $d_1,d_2,d_3$ and $I=(f_1,f_2,f_3)$ such that $S/I$ has the (SLP) then $J=Gin(I)$ is uniquely determined and
is an almost reverse lexicographic ideal.
\end{teo}

\begin{proof}
The theorem is a direct consequence of the Propositions $2.2.3$, $2.2.8$, $2.3.3$, $2.3.8$, $2.3.13$ and $2.3.17$.
\end{proof}

\pagebreak
\begin{teo}
The conjecture Moreno is true for $n=3$ (and $char(K)=0$).
\end{teo}

\begin{proof}
Notice that (SLP) is an open condition. Also, the condition that a sequence of homogeneous polynomial is regular 
is an open condition. It follows, using Theorem $2.1.1$, that for a generic sequence $f_1,f_2,f_3$ of homogeneous polynomials of given degrees $d_1,d_2,d_3$, $J=Gin(I)$ is almost revlex, where $I=(f_1,f_2,f_3)$. But the definition of the generic initial ideal implies to choose a generic change of variables, and therefore for a generic sequence $f_1,f_2,f_3$ of homogeneous polynomials of given degrees $in(I)$ is almost revlex, as required.
\end{proof}

From now on, we fix some integers $2\leq d_1\leq d_2\leq d_3$. Let $f_1,f_2,f_3\in S=K[x_1,x_2,x_3]$ be a regular sequence of homogeneous polynomials of given degrees $d_1,d_2,d_3$ and $I=(f_1,f_2,f_3)$ such that $S/I$ has the (SLP) and let $J=Gin(I)$ be the generic initial ideal of $I$, with respect to revers lexicographic order. We will use these conventions in the sections $2$ and $3$ of this chapter.

\begin{rem}
In order to compute $J$ we will use the fact that $J$ is a strongly stable ideal, i.e. for any monomial $u\in J$ and any indices $j<i$, if $x_{i}|u$ then $x_{j}u/x_{i}\in J$. Also, a theorem of Wiebe (see \cite{W}) states that $S/I$ has (SLP) if and only if $x_{3}$ is a strong Lefschetz element for $S/J$.
We need to consider several cases: I. $d_1+d_2\leq d_3+1$ with $2$ subcases $d_{1}=d_{2}<d_{3}$, $d_{1}<d_{2}<d_{3}$ (section $2$) and II. $d_1+d_2>d_3+1$ with $4$ subcases: $d_{1}=d_{2}=d_{3}$, $d_{1}=d_{2}<d_{3}$, $d_{1}<d_{2}=d_{3}$, $d_{1}<d_{2}<d_{3}$ (section $3$). 

The construction of $J$ in all cases, follows the next procedure.
For any nonnegative integer $k$, we denote by $J_{k}$ the set of monomials of degree $k$ in $J$. We can easily compute the cardinality of each $J_{k}$ from the Hilbert series of $S/J$. We denote $Shad(J_{k})=\{x_{i}u:\;1\leq i\leq n,\;u\in J_{k}\}$.  We begin with $J_{0}=\emptyset$ and we pass from $J_{k}$ to $J_{k+1}$ noticing that $J_{k+1} = Shad(J_{k})\cup $ eventually some new monomial(s) (exactly $|J_{k+1}|-|Shad(J_{k})|$ new monomials). The fact that $J$ is strongly stable and that $x_{3}$ is a strong Lefschetz element for $S/J$ tell us what we need to add to $Shad(J_{k})$ in order to obtain $J_{k+1}$. We continue this procedure until $k=d_1 + d_2 +d_3 -2$ since $J_{k} = S_{k}$(=the set of all monomials of degree $k$) for any $k\geq d_1 + d_2 +d_3 -2$ and so we cannot add any new monomials in larger degrees. $J$ is the ideal generated by all monomials added to $Shad(J_{k})$ at some step $k$. \vspace{2 pt}
We will present detailed this construction only in the subcase $d_{1}=d_{2}<d_{3}$ of the case $I.d_1+d_2\leq d_3+1$ (Proposition $2.2.3$), the other cases being presented in sketch, but the reader can easily complete the proofs.
The condition that $S/I$ has (SLP) is needed. Indeed, there are monomial ideals $J$ such that $S/J$
and $S/I$ have the same Hilbert series, $J$ is strongly stable and $x_3$ is a weak Lefschetz element for $S/J$, but
is not strong Lefschetz. For example, 
\[ J=(x_1\{x_1,x_2\}^{2},\;x_2^4,\;x_1^2x_3^2,\;x_2^3,x_3^2,\;x_1x_2x_3^3,\;x_1x_3^4,\;x_2^2x_3^4,\;x_2x_3^5,\;x_3^7)\]
has $H(S/J,t)=(1+t+t^{2})^{3}$, i.e. the Hilbert series of a complete intersection $S/(f_1,f_2,f_3)$, with $f_1,f_2,f_3$ homogeneous of degree $3$. Also, $J$ is strongly stable and $x_3$ is a weak Lefschetz element for $S/J$, but not strong Lefschetz, since the map $(S/J)_2 \stackrel{\cdot x_3^2}{\rightarrow} (S/J)_4$ is not injective and $|(S/J)_2|=|(S/J)_4|=6$.
\end{rem}

\section{Case $d_1+d_2\leq d_3+1$.}

\medskip
\begin{itemize}
\item Subcase $d_1=d_2<d_3$.
\end{itemize} 
\medskip

\begin{prop}
\em Let $2\leq d:=d_1=d_2<d_3$ be positive integers such that $2d\leq d_3+1$. The Hilbert function of the standard graded complete intersection $A=K[x_{1},x_{2},x_{3}]/I$, where $I$ is the ideal generated by $f_{1}$, $f_{2}$, $f_{3}$, with $f_i$ homogeneous polynomials of degree $d_i$, for all $i$, with $1\leq i\leq 3$, has the form:  
\begin{enumerate}
	\item $H(A,k) = \binom{k+2}{k}$, for $k\leq d-1$.
	\item $H(A,k) = \binom{d+1}{2} +  \sum_{i=1}^{j}(d-i)$, 
	      for $k=d-1+j$, where $0\leq j \leq d-1$.
	\item $H(A,k) = d^{2}$, for $2d-2\leq k \leq d_{3}-1 $.
	\item $H(A,k) = H(A,2d+d_{3}-3-k)$ for $k\geq d_{3}$.
\end{enumerate}
\end{prop}

\begin{proof}
It follows from \cite[Lemma 2.9(a)]{PV}.
\end{proof}

\begin{cor}
\em In the conditions of Proposition $2.2.1$, let $J=\Gin(I)$ be the generic initial ideal of $I$ with respect to the reverse lexicographic order. If we denote by $J_{k}$ the set of monomials of $J$ of degree $k$, then:
\begin{enumerate}
	\item $J_{k} = \emptyset$, for $k\leq d-1$.
	\item $|J_{k}| = j(j+1)$, for $k=d-1+j$, where $0\leq j \leq d-1$.
	\item $|J_{k}| = d(d-1) + 2dj + \frac{j(j-1)}{2}$, for $k=2d-2+j$, where $0\leq j \leq d_3+1-2d$.
  \item $|J_{k}| = \frac{d_{3}(d_{3}+1)}{2} -d^2 +jd_{3} + j(j+1)$, for $k=d_{3}-1+j$, where $0\leq j \leq d-1$.
  \item $|J_{k}| = \frac{d_3(d_3-1)}{2}+d(d_3-1)+j(d_{3}+2d)$, for $k=d+d_{3}-2+j$, where $0\leq j \leq d$.
  \item $J_{k} = S_{k}$, for $k\geq 2d+d_{3}-2$, where $S_{k}$ is the set of monomials of degree $k$.
\end{enumerate}
\end{cor}

\begin{proof}
Using that $|J_k| = |S_k| - H(S/J,k)$, together with the general fact that $H(S/J,k)=H(S/I,k)$, the proof follows immediately from Proposition $2.2.1$. 
\end{proof}

I would like to express my gratitude to Dr.Marius Vladoiu for his help on the proof of the 
following proposition and, in general, for his help on the sections $2.2$ and $2.3$.

\begin{prop}
\em Let $2\leq d:=d_1=d_2<d_3$ be positive integers such that $2d\leq d_3+1$. Let $f_{1},f_{2},f_{3}\in K[x_{1},x_{2},x_{3}]$ be a regular sequence of homogeneous polynomials of degrees $d_{1},d_{2},d_{3}$. If  $I=(f_{1},f_{2},f_{3})$, and $J=\Gin(I)$, the generic initial ideal with respect to the reverse lexicographic order, and
$S/I$ has (SLP), then:
\[ J = ( x_{1}^{d}, x_{1}^{d-j-1}x_{2}^{2j+1}\;for\;0\leq j\leq d-1, 
x_{2}^{2d-2j-2}x_{3}^{d_{3}-2d+2j+2}\{x_{1},x_{2}\}^{j}\;for\;0\leq j\leq d-2,\]
\[ x_{3}^{d_{3}+2j-2}\{x_{1},x_{2}\}^{d-j}\; for\; 1\leq j\leq d). \]
\end{prop}

\begin{proof}
We have $|J_{d}|=2$, hence 
$J_{d}=\{ x_{1}^{d-1}\{x_{1},x_{2}\} \}$, since $J$ is a strongly stable ideal. Therefore:
\[ 
Shad(J_{d}) = \{x_{1}^{d-1}\{x_{1},x_{2}\}^{2} ,x_{1}^{d-1}x_{3}\{x_{1},x_{2}\}\},
\]
 Now we have two possibilities to analyze: $d=2$ and $d\geq 3$. First, suppose $d\geq 3$.

Using the formulas from Corollary $2.2.2$ we have $|J_{d+1}|-|Shad(J_{d})|= 1$, so there is only one generator to add to the set $Shad(J_{d})$ in order to obtain $J_{d+1}$.  Since $J$ is strongly stable, we have only two possibilities: $x_{1}^{d-2}x_{2}^{3}$ or $x_{1}^{d-1}x_{3}^{2}$. We cannot have $x_{1}^{d-1}x_{3}^{2}$, since otherwise the application
\[ 
(S/J)_{d-1} \stackrel{\cdot x_{3}^{2}}{\longrightarrow} (S/J)_{d+1},
\]
with $|(S/J)_{d-1}|<|(S/J)_{d+1}|$ (see Proposition $2.2.1$) would not be injective ($0\neq x_{1}^{d-1}\in (S/J)_{d-1}$ and is mapped to $0$), which is a contradiction to the fact that $x_{3}$ is a strong Lefschetz element for $S/J$. Hence:
\[
J_{d+1} = \{x_{1}^{d-2}\{x_{1},x_{2}\}^{3} ,x_{1}^{d-1}x_{3}\{x_{1},x_{2}\}\}.
\]
We prove by induction on $j$, with $1\leq j \leq d-2$, that:
\[ 
J_{d+j} = Shad(J_{d+j-1})\cup \{x_{1}^{d-j-1}x_{2}^{2j+1}\} =
\]
\[
= \{x_{1}^{d-j-1}\{x_{1},x_{2}\}^{2j+1}, x_{1}^{d-j}x_{3}\{x_{1},x_{2}\}^{2j-1},\ldots,x_{1}^{d-1}x_{3}^{j}\{x_{1},x_{2}\}\}.
\]
The assertion was checked above for $j=1$. Assume now that the statement is true for some $j<d-2$. Then $Shad(J_{d+j})$ is the following set:
\[
\{x_{1}^{d-j-1}\{x_{1},x_{2}\}^{2j+2},x_{1}^{d-j-1}x_{3}\{x_{1},x_{2}\}^{2j+1}, x_{1}^{d-j}x_{3}^2\{x_{1},x_{2}\}^{2j-1},\ldots,x_{1}^{d-1}x_{3}^{j+1}\{x_{1},x_{2}\}\}. 
\]
We have $|J_{d+j+1}| - |Shad(J_{d+j})| = 1$, so we must add only one generator to $Shad(J_{d+j})$ to get $J_{d+j+1}$. The ideal $J$, being strongly stable, allows only two possibilities, namely $x_{1}^{d-j-2}x_{2}^{2j+3}$ or $x_{1}^{d-2}x_{2}^2x_3^{j+1}$. The second one is not allowed because the application
\[
(S/J)_{d} \stackrel{\cdot x_{3}^{j+1}}{\longrightarrow} (S/J)_{d+j+1}
\]
would not be injective ($0\neq x_{1}^{d-2}x_2^2\in (S/J)_{d}$ and is mapped to $0$), $x_{3}$ being a strong-Lefschetz element for $S/J$. Therefore, we must add $x_{1}^{d-j-2}x_{2}^{2j+3}$ and our claim is proved. In particular, we obtain
\[
J_{2d-2} =\{x_{1}\{x_{1},x_{2}\}^{2d-3}, x_{1}^{2}x_{3}\{x_{1},x_{2}\}^{2d-5},\ldots,x_{1}^{d-1}x_{3}^{d-2}\{x_{1},x_{2}\}\}.
\]

In order to compute $J_{k}$, with $2d-2\leq k \leq d_3$, we must consider two posibilities. 

\begin{itemize}
	\item 1. $d_3=2d-1$. 
\end{itemize}

Since $|J_{2d-1}|-|Shad(J_{2d-2})| = 2$, there are two generators to add to $Shad(J_{2d-2})$. We prove that these generators are
$x_2^{2d-1},x_2^{2d-2}x_3$. Assuming by contradiction that we have other generators, since $J$ is strongly stable, it follows that
there is at least one generator from the set $\{x_{1}x_{2}^{2d-4}x_{3}^{2},x_{1}^{2}x_{2}^{2d-6}x_{3}^{3},\ldots, x_{1}^{d-2}x_{2}^{2}x_{3}^{d-1}\}$. Then, the application $(S/J)_{2d-3} \stackrel{\cdot x_{3}^{2}}{\longrightarrow} (S/J)_{2d-1}$
would not be injective, a contradiction since $x_3$ is a strong Lefschetz element for $S/J$. Therefore
\[
J_{2d-1} = J_{d_3} = \{ \{x_{1},x_{2}\}^{2d-1},  x_{3}\{x_{1},x_{2}\}^{2d - 2},\ldots,x_{1}^{d-1}x_{3}^{d-1}\{x_{1},x_{2}\}\}.
\]

\begin{itemize}
	\item 2. $d_3>2d-1$. 
\end{itemize}	

Since $|J_{2d-1}|-|Shad(J_{2d-2})| = 1$, there is only one generator to add to $Shad(J_{2d-2})$, which can be selected from the set $\{x_2^{2d-1},x_1x_2^{2d-4}x_3^2,x_1^2x_2^{2d-6}x_3^3,\ldots,x_1^{d-2}x_2^2x_3^{d-1}\}$ because $J$ is strongly stable. In a similar manner to what we have done above can be shown that, $x_3$ being a strong Lefschetz  element, leaves us as unique possibility $x_{2}^{2d-1}$, therefore:
\[
   J_{2d-1}=\{\{x_{1},x_{2}\}^{2d-1}, x_{1}x_{3}\{x_{1},x_{2}\}^{2d-3},\ldots,x_{1}^{d-1}x_{3}^{d-1}\{x_{1},x_{2}\}\}. 
\] 
One can easily show, using induction on $1\leq j \leq d_{3}-2d$, if case, that $|J_{2d-1+j}| = |Shad(J_{2d-2+j})|$ and $J_{2d-1+j}$ is the set \[ \{ \{x_{1},x_{2}\}^{2d-1+j}, \ldots ,x_{3}^{j}\{x_{1},x_{2}\}^{2d-1},
x_{3}^{j+1}x_{1}\{x_{1},x_{2}\}^{2d-2},\ldots,x_{3}^{d+j-1}x_{1}^{d-1}\{x_{1},x_{2}\}\} .\]
In particular, we obtain that $J_{d_{3}-1}$ is the set
\[ \{\{x_{1},x_{2}\}^{d_{3}-1}, \ldots,x_{3}^{d_{3}-2d}\{x_{1},x_{2}\}^{2d-1}, x_{3}^{d_{3}-2d+1}x_{1}\{x_{1},x_{2}\}^{2d-3},\ldots,x_{3}^{d_{3}-d-1}x_{1}^{d-1}\{x_{1},x_{2}\}\}.\]

Since $|J_{d_3}|-|Shad(J_{d_{3}-1})|=1$, the generator which has to be add to $Shad(J_{d_{3}-1})$ can be selected from the set
$\{x_{2}^{2d-2}x_3^{d_3 - 2d+2}, x_{1}x_2^{2d-4}x_3^{d_3 - 2d+3}, \ldots ,x_{1}^{d-2}x_{2}^{2}x_{3}^{d_{3}-d}\}$ such that $J$
is strongly stable. The generator is $x_{2}^{2d-2}x_3^{d_3 - 2d+2}$, otherwise the application $(S/J)_{2d-3} \stackrel{\cdot x_{3}^{d_{3}-2d+3}}{\longrightarrow} (S/J)_{d_{3}}$ is not injective, a contradiction, since $x_3$ is a strong Lefschetz element
for $S/J$. Hence, we get that $J_{d_{3}}$ is
\[  \{\{x_{1},x_{2}\}^{d_{3}}, \ldots, x_{3}^{d_{3}-2d+2}\{x_{1},x_{2}\}^{2d-2},x_{1}^{2}x_{3}^{d_{3}-2d+3}\{x_{1},x_{2}\}^{2d-5},\ldots,x_{3}^{d_{3}-d}x_{1}^{d-1}\{x_{1},x_{2}\}\},\] and one can check that is the same formula as in 1.($d_3 = 2d-1$).


Now, we show by induction on $1\leq j\leq d-2$ that
\[ 
J_{d_{3}+j} = Shad(J_{d_{3}-1+j}) \cup \{x_{2}^{2d-2j-2}x_{3}^{d_{3}-2d+2j+2}\{x_{1},x_{2}\}^{j} \}. 
\]
Indeed, for $j=1$, $|J_{d_{3}+1}|-|Shad(J_{d_{3}})|=2$ and the generators which must be added are $x_{1}x_{2}^{2d-4}x_{3}^{d_{3}-2d+4}, x_{2}^{2d-3}x_{3}^{d_{3}-2d+4}$. If not, since $J$ is strongly stable, then at least one of the generators belongs to the set
$\{ x_{1}^{2}x_{2}^{2d-6}x_{3}^{d_{3}-2d+5},\ldots,x_{1}^{d-2}x_{2}^{2}x_{3}^{d_{3}-d+1} \}$ (for $d=3$ this is the emptyset).
but then the map $(S/J)_{2d-4} \stackrel{x_{3}^{d_{3}-2d+5}}{\longrightarrow} (S/J)_{d_{3}+1}$ is not injective, contradiction.

Assume now that we proved the assertion for some $j<d-2$. Then $|J_{d_{3}+j+1}|-|Shad(J_{d_{3}+j})| = j+2$ and the new generators are $x_{2}^{2d-2j-4}x_{3}^{d_{3}-2d+2j+4}\{x_{1},x_{2}\}^{j+1}$. Indeed, if not, since $J$ is strongly stable, then at least one of the generators belongs to the set
$\{x_{1}^{j+2}x_{2}^{2d-2j-6}x_{3}^{d_{3}-2d+2j+5},\ldots,x_{1}^{d-2}x_{2}^{2}x_{3}^{ d_{3}-d+j+1}\}$ (for $d=3$ this is the emptyset...) but then the map $(S/J)_{2d-j-4} \stackrel{x_{3}^{d_{3}-2d+2j+5}}{\longrightarrow} (S/J)_{d_{3}+j+1}$ is not injective, contradiction, and we are done. Hence, 
\[J_{d_{3}+d-2} = \{ \{x_{1},x_{2}\}^{d_{3}+d-2},x_{3}\{x_{1},x_{2}\}^{d+d_{3}-3},\ldots ,x_{3}^{d_{3}-2}\{x_{1},x_{2}\}^{d}\}.\]

We prove by induction on $1\leq j \leq d$ that $J_{d+d_{3}-2+j} = Shad(J_{d+d_{3}-3+j})\cup x_{3}^{d_{3}+2j-2}\{x_{1},x_{2}\}^{d-j}$.
If $j=1$ then $|J_{d+d_{3}-1}|- |Shad(J_{d+d_{3}-2})| = d$ so we must add $d$ generators, which are precisely the elements of the set $x_{3}^{d_{3}}\{x_{1},x_{2}\}^{d-1}$. Indeed, if we have a generator which does not belong to the set it is divisible by $x_{3}^{d_{3}+1}$ and therefore the map $(S/J)_{d-2} \stackrel{\cdot x_{3}^{d_{3}+1}}{\longrightarrow} (S/J)_{d_{3}+d-1}$ is not injective, which is a contradiction with $x_{3}$ is a strong Lefschetz element for $S/J$ (the map has to be bijective). The induction
step is similar and finally we obtain that $J_{d_{3}+2d-2}=S_{d_{3}+2d-2}$ and thus we cannot add new minimal generators of $J$ in degree $>d_{3}+2d-2$.

In order to complete the proof we must consider now $d=2$. The hypothesis implies $d_{3}\geq 3$. We already seen that $J_{2}=\{x_{1}^{2},x_{1}x_{2}\}$ and $Shad(J_{2}) = \{x_{1}\{x_{1},x_{2}\}^{2} ,x_{1}x_{3}\{x_{1},x_{2}\}\}$.

Using the formulas from Corollary $2.2.2$ we have $|J_{3}|-|Shad(J_{2})|= 1$, so there is only one generator to add to the set $Shad(J_{d})$ in order to obtain $J_{d+1}$. 

Since $J$ is strongly stable, we have only two possibilities: $x_{2}^{3}$ or $x_{1}x_{3}^{2}$. We can not have $x_{1}x_{3}^{2}$, since otherwise the application
\[ 
(S/J)_{1} \stackrel{\cdot x_{3}^{2}}{\longrightarrow} (S/J)_{3},
\]
with $|(S/J)_{1}|<|(S/J)_{3}|$ (see Proposition $2.2.1$) would not be injective ($0\neq x_{1}\in (S/J)_{1}$ and is mapped to $0$), which is a contradiction to the fact that $x_{3}$ is a strong Lefschetz element for $S/J$. Hence:
\[
J_{3} = \{\{x_{1},x_{2}\}^{3} ,x_{1}x_{3}\{x_{1},x_{2}\}\}.
\]
Assume now $d_{3}\geq 4$. One can easily show, using induction on $1\leq j \leq d_{3}-4$, if case, that $|J_{3+j}| = |Shad(J_{2+j})|$ and $J_{3+j}$ is the set \[ \{ \{x_{1},x_{2}\}^{3+j}, \ldots ,x_{3}^{j}\{x_{1},x_{2}\}^{3},
x_{3}^{j+1}x_{1}\{x_{1},x_{2}\}\}.\]
In particular, we obtain that \[J_{d_{3}-1} = \{\{x_{1},x_{2}\}^{d_{3}-1}, \ldots,x_{3}^{d_{3}-4}\{x_{1},x_{2}\}^{3}, x_{3}^{d_{3}-3}x_{1}\{x_{1},x_{2}\}\}.\]

Since $|J_{d_3}|-|Shad(J_{d_{3}-1})|=1$, the generator which has to be add to $Shad(J_{d_{3}-1})$ is exactly $x_{2}^{2}x_3^{d_3 - 2}$ such that $J$ is strongly stable. Hence, we get
\[J_{d_{3}} =  \{\{x_{1},x_{2}\}^{d_{3}}, \ldots, x_{3}^{d_{3}-3}\{x_{1},x_{2}\}^{3} , x_{1}x_{3}^{d_{3}-2}\{x_{1},x_{2}\}\},\] and one can check that is the same formula as in the case $d_3 = 3$.

Since $|J_{d_3+1}|-|Shad(J_{d_{3}})|=1$, there is only one generator to add to the set $Shad(J_{d_{e}})$ in order to obtain $J_{d_3+1}$. Since $J$ is strongly stable, we have only two possibilities: $x_{2}^{2}x_{3}^{d_{3}-1}$ or $x_{1}x_{3}^{d_{3}}$. We can not have $x_{1}x_{3}^{d_{3}}$, since otherwise the application
\[ 
(S/J)_{1} \stackrel{\cdot x_{3}^{d_{3}}}{\longrightarrow} (S/J)_{d_{3}+1},
\]
with $|(S/J)_{1}|=|(S/J)_{3}|$ (see Proposition $2.2.1$) would not be injective ($0\neq x_{1}\in (S/J)_{1}$ and is mapped to $0$), which is a contradiction to the fact that $x_{3}$ is a strong Lefschetz element for $S/J$. Hence:
\[
J_{d_{3}+1} = \{\{x_{1},x_{2}\}^{d_{3}+1}, \ldots, x_{3}^{d_{3}-1}\{x_{1},x_{2}\}^{2} \}.
\] 
Since $|J_{d_3+2}|-|Shad(J_{d_{3}+1})|=2$ and $J$ is strongly stable, we must add $x_{1}x_{3}^{d_{3}+1}$ and
$x_{2}x_{3}^{d_{3}+1}$ at $Shad(J_{d_{3}+1})$ in order to obtain $J_{d_3+2}$. Hence 
$J_{d_{3}+2}=S_{d_{3}+2}\setminus \{x_{3}^{d_{3}+2}\}$. Finally, since $J_{d_{3}+3}=S_{d_{3}+3}$ we add $x_{3}^{d_{3}+3}$ at $Shad(J_{d_{3}+2})$ and thus we cannot add new minimal generators of $J$ in degree $>d_{3}+2$.
\end{proof}

\begin{cor}
In the conditions of the above proposition, the number of minimal generators of $J$ is $d^{2}+d+1$.
\end{cor}

\begin{ex}
Let $d_{1}=d_{2}=3$ and $d_{3}=9$. Proposition $2.2.3$ implies:
\[ J=(x_{1}^{3},\; x_{1}^{2}x_{2},\; x_{1}x_{2}^{3},\; x_{2}^{5},\; x_{2}^{4}x_{3}^{5},\; x_{1}x_{2}^{2}x_{3}^{7},\;
      x_{2}^{3}x_{3}^{7},\;x_{3}^{9}\{x_{1},x_{2}\}^{2},\;x_{3}^{11}\{x_{1},x_{2}\},\;x_{3}^{13}    ). \]
\end{ex}

\medskip
\begin{itemize}
\item Subcase $d_1<d_2<d_3$.
\end{itemize} 
\medskip

\begin{prop}
\em Let $2\leq d_1<d_2<d_3$ be positive integers such that \linebreak $d_1 + d_2\leq d_3+1$. The Hilbert function of the standard graded complete intersection $A=K[x_{1},x_{2},x_{3}]/I$, where $I$ is the ideal generated by $f_{1}$, $f_{2}$, $f_{3}$, with $f_i$ homogeneous polynomials of degree $d_i$, for all $i$, with $1\leq i\leq 3$, has the form:  
\begin{enumerate}
	\item $H(A,k) = \binom{k+2}{k}$, for $k\leq d_{1}-1$.
	\item $H(A,k) = \binom{d_{1}+1}{2} + jd_{1}$, for $k=j + d_1 - 1$, where $0\leq j\leq d_2 - d_1$.
	\item $H(A,k) = \binom{d_{1}+1}{2} + d_{1}(d_{2}-d_{1}) + \sum_{i=1}^{j}(d_{1}-i)$, 	      
	      for $k=j+d_{2}-1$, where $0 \leq j \leq d_{1}-1$.
	\item $H(A,k) = d_{1}d_{2}$, for $d_{1}+d_{2}-2\leq k \leq d_{3}-1 $.
	\item $H(A,k) = H(A,d_{1}+d_{2}+d_{3}-3-k)$ for $k\geq d_{3}$.
\end{enumerate}
\end{prop}

\begin{proof}
It follows from \cite[Lemma 2.9(a)]{PV}.
\end{proof}

\begin{cor}
\em In the conditions of Proposition $2.2.6$, let $J=\Gin(I)$ be the generic initial ideal of $I$ with respect to the reverse lexicographic order. If we denote by $J_{k}$ the set of monomials of $J$ of degree $k$, then:
\begin{enumerate}
	\item $|J_{k}| = 0$, for $k\leq d-1$.
	\item $|J_{k}| = j(j+1)/2$, for $k=j+d_1-1$, where $0\leq j\leq d_2 - d_1$.	
	\item $|J_{k}| = \frac{(d_{2}-d_{1})((d_{2}-d_{1}-1))}{2} + j(d_{2}-d_{1}) + j(j+1)$, 
	      for $k=j+d_{2}-1$, where $0 \leq j \leq d_{1}-1$.	      
  \item $|J_{k}| = \frac{d_{1}^{2}+d_{2}^{2}-d_{1}-d_{2}}{2} + j(d_{1}+d_{2}) + \frac{j(j-1)}{2}$, 
        for $k=j+d_{1}+d_{2}-2$, where $0 \leq j \leq d_{3}-d_{1}-d_{2}+1$.        
  \item $|J_{k}| = \frac{d_{3}^{2}+d_{3}-2d_{1}d_{2}}{2} +jd_{3} + j(j+1)$, 
        for $k=j+d_{3}-1$ where $0 \leq j \leq d_{1}-1$.
  
  \item $|J_{k}| = \frac{ (d_{1}+d_{3})(d_{1}+d_{3}-1) +d_{1}^{2} - d_{1} -2d_{1}d_{2}}{2} +j(d_{3}+2d_{1}) +
         \frac{j(j-1)}{2}$, for $k = j + d_1 + d_3 - 2$, where $0 \leq j \leq d_{2}-d_{1}$ .
  \item $|J_{k}| = \frac{(d_{2}+d_{3})(d_{2}+d_{3}-1) +d_{1}(d_{1}-1)}{2} + j(d_{1}+d_{2}+d_{3})$, 
        for $k=d_{2}+d_{3}-2$, where $0 \leq j \leq d_{1}-1$.	
  \item $J_{k} = S_{k}$, for $k\geq 3d-2$.
\end{enumerate}
\end{cor}

\begin{prop}
\em Let $2\leq d_1<d_2<d_3$ be positive integers such that $d_1 + d_2\leq d_3+1$. Let $f_{1},f_{2},f_{3}\in K[x_{1},x_{2},x_{3}]$ be a regular sequence of homogeneous polynomials of degrees $d_{1},d_{2},d_{3}$. If $I=(f_{1},f_{2},f_{3})$ ,$J=\Gin(I)$, the generic initial ideal with respect to the reverse lexicographic order, and $S/I$ has (SLP),then:
\[ J = ( x_{1}^{d_{1}}, x_{1}^{d_{1}-j} x_{2}^{d_{2}-d_{1}+2j-1}\;for\;1\leq j\leq d_1 - 1,\; x_{2}^{d_1+d_2 -1},\;
x_{2}^{d_1+d_2 -2}x_{3}^{d_{3}-d_1-d_2+2}, \]
\[ x_{3}^{d_{3}-d_{1}-d_{2}+2j+2}x_{2}^{d_{1}+d_{2}-2j-2}\{x_{1},x_{2}\}^{j}\;for\;1\leq j\leq d_1-2,\; \]\[
x_{3}^{d_{3}+d_{1}-d_{2}-2+2j}x_{2}^{d_{2}-d_{1}+1-j}\{x_{1},x_{2}\}^{d_1 - 1}\;for\;1\leq j\leq d_2-d_1,\;\]\[
x_{3}^{d_{3}+d_2-d_1+2j-2}\{x_{1},x_{2}\}^{d_1-j}\;for\;1\leq j\leq d_1 ).\]
\end{prop}

\begin{proof}
We have $|J_{d_1}|=1$, hence $J_{d_1}=\{x_{1}^{d_1}\}$, since $J$ is a strongly stable ideal. Therefore:
\[ 
Shad(J_{d_1}) = \{x_{1}^{d_1}\{x_{1},x_{2}\} ,x_{1}^{d_1}x_{3}\}.
\]
Assume $d_2>d_1+1$. Since $|J_{d_1+1}|=|Shad(J_{d_1})|$ from the formulas of $2.2.7$, if follows $J_{d_1+1}=Shad(J_{d_1})$. We prove by induction on $1\leq j\leq d_2-d_1-1$ that 
\[ J_{d_1+j}=Shad(J_{d_1+j-1})= \{x_{1}^{d_{1}}\{x_{1},x_{2}\}^{j}, x_{3}x_{1}^{d_{1}}\{x_{1},x_{2}\}^{j-1},\ldots,x_{3}^{j}x_{1}^{d_{1}} \}. \]
Indeed, the case $j=1$ is already proved. Suppose the assertion is true for some $j<d_2-d_1-1$. Since
$|J_{d_1+j+1}|-|Shad(J_{d_1+j})|=0$ it follows that 
\[ J_{d_1+j+1}=Shad(J_{d_1+j})= \{x_{1}^{d_{1}}\{x_{1},x_{2}\}^{j+1}, x_{3}x_{1}^{d_{1}}\{x_{1},x_{2}\}^{j},\ldots,x_{3}^{j+1}x_{1}^{d_{1}} \}\]
thus we are done. In particular, we get
\[ J_{d_{2}-1} = \{x_{1}^{d_{1}}\{x_{1},x_{2}\}^{d_{2}-d_{1}-1}, x_{3}x_{1}^{d_{1}}\{x_{1},x_{2}\}^{d_{2}-d_{1}-2},\ldots,x_{3}^{d_{2}-d_{1}-1}x_{1}^{d_{1}}\} \]
which is the same formula as in the case $d_2=d_1+1$.

We have $|J_{d_2}| - |Shad(J_{d_{2}-1})| = 1$ so we must add a new generator to $Shad(J_{d_{2}-1})$ to obtain $J_{d_{2}}$. Since $J$ is strongly stable and $x_3$ is a strong Lefschetz element for $S/J$ this new generator is $x_{1}^{d_{1}-1}x_{2}^{d_{2}-d_{1}+1}$, therefore
\[ J_{d_{2}} = \{x_{1}^{d_{1}-1}\{x_{1},x_{2}\}^{d_{2}-d_{1}+1}, x_{3}x_{1}^{d_{1}}\{x_{1},x_{2}\}^{d_{2}-d_{1}-1},\ldots,x_{3}^{d_{2}-d_{1}}x_{1}^{d_{1}}\}.\]
Assume $d_1>2$. We prove by induction on $1\leq j\leq d_{1}-1$ that
\[ J_{d_{2}-1+j} = Shad(J_{d_{2}-2+j}) \cup \{ x_{1}^{d_{1}-j} x_{2}^{d_{2}-d_{1}+2j-1}\} =  \{x_{1}^{d_{1}-j}\{x_{1},x_{2}\}^{d_{2}-d_{1}+2j-1},\]\[ x_{3}x_{1}^{d_{1}-j+1}\{x_{1},x_{2}\}^{d_{2}-d_{1}+2j-3}, \ldots, x_{3}^{j}x_{1}^{d_{1}}\{x_{1},x_{2}\}^{d_{2}-d_{1}-1},\ldots,x_{3}^{d_{2}-d_{1}+j-1}x_{1}^{d_{1}}\}.\]
The assertion was proved for $j=1$. Suppose $1\leq j<d_1-1$ and the assertion is true for $j$. We have
$|J_{d_{2}+j}|-|Shad(J_{d_{2}-1+j})|=1$, thus we must add a new generator to $Shad(J_{d_{2}-1+j})$ in order to
obtain $J_{d_{2}+j}$ and since $J$ is strongly stable and $x_3$ is a strong Lefschetz element for $S/J$, this is $x_{1}^{d_{1}-j-1} x_{2}^{d_{2}-d_{1}+2j+1}$ and we are done. In particular, we obtain:
\[ J_{d_{1}+d_{2}-2} = \{x_{1}\{x_{1},x_{2}\}^{d_{1}+d_{2}-3}, x_{3}x_{1}^{2}\{x_{1},x_{2}\}^{d_{1}+d_{2}-5}, \ldots,\] 
\[ x_{3}^{d_{1}-1}x_{1}^{d_{1}}\{x_{1},x_{2}\}^{d_{2}-d_{1}-1},
x_{3}^{d_{1}}x_{1}^{d_{1}}\{x_{1},x_{2}\}^{d_{2}-d_{1}-2},\ldots,x_{3}^{d_{2}-2}x_{1}^{d_{1}} \}\]
and one can check that is the same expression as in the case $d_1=2$.

In order to compute $J_{k}$, with $2d-2\leq k \leq d_3$, we must consider two possibilities. 

\begin{itemize}
	\item 1.$d_3=d_1+d_2-1$. 
\end{itemize}

Since $|J_{d_1+d_2-1}|-|Shad(J_{d_1+d_2-2})| = 2$, there are two generators to add to $Shad(J_{d_1+d_2-2})$ to get $J_{d_1+d_2-1}$, but on the other hand $J$ is strongly stable and $x_3$ is a strong Lefschetz element for $S/J$ so these generators must be $x_2^{d_1+d_2-1},x_2^{d_1+d_2-2}x_3$. Therefore:
\[J_{d_1+d_2-1} = J_{d_3} = \{ \{x_{1},x_{2}\}^{d_1+d_2-1},  x_{3}\{x_{1},x_{2}\}^{d_1 +d_2 - 2}, \ldots,x_{1}^{d_1}x_{3}^{d_2-1}\}.\] 

\begin{itemize}
	\item 2.$d_3>d_1+d_2-1$. 
\end{itemize}	

Since $|J_{d_1+d_2-1}|-|Shad(J_{d_1+d_2-2})| = 1$, there is only one generator to add to $Shad(J_{2d-2})$, which is precisely $x_{2}^{d_1+d_2-1}$ since $J$ is strongly stable and $x_3$ is a strong Lefschetz element for $S/J$. Therefore
\[ J_{d_1+d_2-1} = \{\{x_{1},x_{2}\}^{d_1+d_2-1}, x_{1}x_{3}\{x_{1},x_{2}\}^{d_1+d_2-3},\ldots, x_{1}^{d_1}x_{3}^{d_2-1} \}.\] 
One can easily show, using induction on $1\leq j \leq d_{3}-d_1-d_2$, if case, that \linebreak
$|J_{d_1+d_2-1+j}| = |Shad(J_{d_1+d_2-2+j})|$ and $J_{d_1+d_2-1+j}$ is the set 
\[\{ \{x_{1},x_{2}\}^{j+d_{1}+d_{2}-1}, x_{3}\{x_{1},x_{2}\}^{j+d_{1}+d_{2}-2}, \ldots, x_{3}^{j}\{x_{1},x_{2}\}^{d_{1}+d_{2}-1},\]
\[ x_{3}^{j+1}x_{1}\{x_{1},x_{2}\}^{d_{1}+d_{2}-3},\ldots,x_{3}^{d_{1}+j}x_{1}^{d_{1}}\{x_{1},x_{2}\}^{d_{2}-d_{1}-1},
\ldots, x_{3}^{j+d_{2}-1}x_{1}^{d_{1}} \} \]

\[So\;J_{d_{3}-1}=\{ \{x_{1},x_{2}\}^{d_{3}-1}, x_{3}\{x_{1},x_{2}\}^{d_{3}-2}, \ldots, x_{3}^{d_{3}-d_{1}-d_{2}}\{x_{1},x_{2}\}^{d_{1}+d_{2}-1},\]
\[ x_{3}^{d_{3}-d_{1}-d_{2}+1}x_{1}\{x_{1},x_{2}\}^{d_{1}+d_{2}-3},\ldots,
x_{3}^{d_{3}-d_{2}}x_{1}^{d_{1}}\{x_{1},x_{2}\}^{d_{2}-d_{1}-1},
\ldots, x_{3}^{d_{3}-d_{1}-1}x_{1}^{d_{1}} \} \]

Since $|J_{d_3}|-|Shad(J_{d_{3}-1})|=1$, $J$ is strongly stable and $x_3$ is a strong Lefschetz element for $S/J$, the generator which has to be added to $Shad(J_{d_{3}-1})$ is $x_{2}^{d_1+d_2-2}x_3^{d_3 - d_1-d_2 +2}$. Hence, we get 
\[ J_{d_3} = \{ \{x_{1},x_{2}\}^{d_{3}}, x_{3}\{x_{1},x_{2}\}^{d_{3}-1}, \ldots, x_{3}^{d_{3}-d_{1}-d_{2}+2}\{x_{1},x_{2}\}^{d_{1}+d_{2}-2}, \]\[
x_{3}^{d_{3}-d_{1}-d_{2}+3}x_{1}^{2}\{x_{1},x_{2}\}^{d_{1}+d_{2}-5},\ldots,
x_{3}^{d_{3}-d_{2}+1}x_{1}^{d_{1}}\{x_{1},x_{2}\}^{d_{2}-d_{1}-1},
\ldots, x_{3}^{d_{3}-d_{1}}x_{1}^{d_{1}}, \} \]
and one can check that is the same formula as in 1.($d_3 = d_1+d_2-1$).

Assume now $d_1>2$. 
We show by induction on $1\leq j \leq d_1-2$ that
\[ J_{d_3+j} = Shad(J_{d_3-1+j})\cup \{x_{1}^{j}x_{3}^{d_{3}-d_{1}-d_{2}+2j+2}x_{2}^{d_{1}+d_{2}-2j-2}, \ldots,   x_{3}^{d_{3}-d_{1}-d_{2}+2j+2}x_{2}^{d_{1}+d_{2}-2-j}\}. \]
Indeed, for $j=1$, $|J_{d_3+1}|-|Shad(J_{d_3})|=2$ so we must add two generators to $Shad(J_{d_3})$ in order to
obtain $J_{d_3+1}$. Since $J$ is strongly stable and $x_3$ is a strong Lefschetz element for $S/J$, these new generators are $x_{1}x_{3}^{d_{3}-d_{1}-d_{2}+4}x_{2}^{d_{1}+d_{2}-2j-4}$ and $x_{3}^{d_{3}-d_{1}-d_{2}+4}x_{2}^{d_{1}+d_{2}-2j-3}$.
Assume now that we proved the assertion for some $j<d_1-2$. Then $|J_{d_3+j+1}|-|Shad(J_{d_3+j})| = j+2$ and since
$J$ is strongly stable and $x_3$ is a strong Lefschetz element for $S/J$, the new generators are $x_{3}^{d_{3}-d_{1}-d_{2}+2j+4}x_{2}^{d_{1}+d_{2}-2j-4}\{x_{1},x_{2}\}^{j+1}$ as required. Hence, we get
\[ J_{d_{1}+d_{3}-2} = \{ \{x_{1},x_{2}\}^{d_{1}+d_{3}-2},\ldots,x_{3}^{d_{3}+d_{1}-d_{2}-2}\{x_{1},x_{2}\}^{d_{2}},\]
\[ x_{1}^{d_{1}}x_{3}^{d_{3}+d_{1}-d_{2}-1}\{x_{1},x_{2}\}^{d_{2}},\ldots,x_{1}^{d_{1}}x_{3}^{d_{3}-2}\},\]
and one can check that is the same formula as in the case $d_1 = 2$.

We have $|J_{d_{1}+d_{3}-1}|-|Shad( J_{d_{1}+d_{3}-2})|=d_1$ so we must add $d_1$ new generators to
$Shad( J_{d_{1}+d_{3}-2})$ and since $J$ is strongly stable and $x_3$ is a strong Lefschetz element for $S/J$, they are $x_{3}^{d_{3}+d_{1}-d_{2}}x_{2}^{d_{2}-d_{1}}\{x_{1},x_{2}\}^{d_1 - 1}$. Therefore $J_{d_{1}+d_{3}-1}$ is the set
\[\{ \{x_{1},x_{2}\}^{d_{1}+d_{3}-1}, .., x_{3}^{d_{3}+d_{1}-d_{2}}\{x_{1},x_{2}\}^{d_{2}-1},
x_{1}^{d_{1}}x_{3}^{d_{3}+d_{1}-d_{2}+1}\{x_{1},x_{2}\}^{d_{2}-d_{1}-2},.., x_{1}^{d_{1}}x_{3}^{d_{3}-1}\}. \]

Suppose $d_1>2$. We prove by induction on $1\leq j\leq d_{2}-d_{1}$ that:
\[J_{d_{1}+d_{3}-2+j} = Shad(J_{d_1+d_3-3+j})\cup \{x_{3}^{d_{3}+d_{1}-d_{2}-2+2j}x_{2}^{d_{2}-d_{1}+1-j}
\{x_{1},x_{2}\}^{d_1 - 1}\} = \]
\[ = \{\{x_{1},x_{2}\}^{d_{1}+d_{3}-2+j},\ldots, x_{3}^{d_{3}+d_{1}-d_{2}-2+2j}\{x_{1},x_{2}\}^{d_{2}-j},\]\[
x_{1}^{d_{1}}x_{3}^{d_{3}+d_{1}-d_{2}-1+2j}\{x_{1},x_{2}\}^{d_{2}-d_{1}-j-1},\ldots ,x_{1}^{d_{1}}x_{3}^{d_{3}-2+j}\}.\]
We already proved this for $j=1$. Suppose the assertion is true for some $j<d_2-d_1$. Since
$|J_{d_{1}+d_{3}-1+j}|-|Shad(J_{d_{1}+d_{3}-2+j})|=d_1$ we must add $d_1$ new generators to $Shad(J_{d_{1}+d_{3}-2+j})$
and these new generators are $x_{3}^{d_{3}+d_{1}-d_{2}-2+2j}x_{2}^{d_{2}-d_{1}+1-j}\{x_{1},x_{2}\}^{d_1-1}$ because $J$ is strongly stable and $x_3$ is a strong Lefschetz element for $S/J$. In particular, we get:
\[J_{d_{2}+d_{3}-2} = \{ \{x_{1},x_{2}\}^{d_{2}+d_{3}-2},x_{3}\{x_{1},x_{2}\}^{d_{2}+d_{3}-3},\ldots ,x_{3}^{d_{3}+d_{2}-d_{1}-2}\{x_{1},x_{2}\}^{d_{1}}\},\]
which is the same formula as in the case $d_1=2$.

We prove by induction on $1\leq j \leq d_1$ that \[ J_{d_2+d_{3}-2+j} = Shad(J_{d_2+d_{3}-3+j})\cup \{ x_{3}^{d_{3}+d_2-d_1+2j-2}\{x_{1},x_{2}\}^{d_1-j} \}.\]
If $j=1$ then $|J_{d_2+d_{3}-1}|- |Shad(J_{d_2+d_{3}-2})| = d_1$ so we must add $d_1$ generators, which are precisely the elements of the set $x_{3}^{d_{3}}\{x_{1},x_{2}\}^{d-1}$ since $J$ is strongly stable and $x_3$ is a strong Lefschetz element for $S/J$. The induction step is similar and finally we obtain that $J_{d_{3}+2d-2}=S_{d_{3}+2d-2}$ and thus we cannot add new minimal generators of $J$ in degrees $>d_{3}+2d-2$.
\end{proof}

\begin{cor}
In the conditions of the above proposition, the number of minimal generators of $J$ is $1+d_1 + d_1d_2$.
\end{cor}

\begin{ex}
Let $d_{1}=3$, $d_{2}=4$ and $d_{3}=9$. Then 
\[ J = (x_{1}^{3},\;\; x_{1}^{2}x_{2}^{2},\;\; x_{1}x_{2}^{4},\;\; x_{2}^{6},\;\; x_{2}^{5}x_3,\;\; x_{3}^{6}x_{2}^{3}\{x_{1},x_{2}\},\;\; \] 
\[ x_{3}^{8}x_{2}\{x_{1},x_{2}\}^{2},\;\; x_{3}^{10}\{x_{1},x_{2}\}^{2},\;\; x_{3}^{12}\{x_{1},x_{2}\},\;\; x_{3}^{14}).\]
\end{ex}

\section{Case $d_1+d_2 > d_3+1$.}

\medskip
\begin{itemize}
\item Subcase $d_1=d_2=d_3$.
\end{itemize} 
\medskip

\begin{prop}
\em Let $2\leq d:=d_1=d_2=d_3$ be positive integers. The Hilbert function of the standard graded complete intersection $A=K[x_{1},x_{2},x_{3}]/I$, where $I$ is the ideal generated by $f_{1}$, $f_{2}$, $f_{3}$, with $f_i$ homogeneous polynomials of degree $d_i$, for all $i$, with $1\leq i\leq 3$, has the form:  
\begin{enumerate}
	\item $H(A,k) = \binom{k+2}{k}$, for $k\leq d-1$.
	\item $H(A,k) = \binom{k+2}{2} - \frac{3j(j+1)}{2}$, for $k=j+d-1$, 
	      where $0\leq j \leq \left\lfloor \frac{d-1}{2} \right\rfloor$.	      
  \item $H(A,k) = H(A,3d-k-3)$, for $k\geq \left\lceil  \frac{3d-3}{2} \right\rceil$.
\end{enumerate}
\end{prop}

\begin{proof}
It follows from \cite[Lemma 2.9(b)]{PV}.
\end{proof}

\begin{cor}
\em In the conditions of Proposition $2.3.1$, let $J=\Gin(I)$ be the generic initial ideal of $I$ with respect to the reverse lexicographic order. If we denote by $J_{k}$ the set of monomials of $J$ of degree $k$, then:
\begin{enumerate}
	\item $|J_{k}| = 0$, for $k\leq d-1$.
	\item $|J_{k}| = \frac{3j(j+1)}{2}$, for $k=d-1+j$, where $0 \leq j \leq [\frac{3d-1}{2}]$.
	\item If $d$ is even, then $|J_{k}| = \frac{3d^{2}+3d(4j+2)}{8} + \frac{3j(j+1))}{2}$, for
	      $k = j + \frac{3d-2}{2}$, where $0\leq j\leq \frac{d-2}{2}$.        
	      
	      \noindent
	      If $d$ is odd, then $|J_{k}| = \frac{3(d^{2}-1)+12jd}{8} + \frac{3j^{2}}{2}$, for
        $k = j + \frac{3d-3}{2}$, where $0\leq j\leq \frac{d-1}{2}$.        
  \item $|J_{k}| = \frac{3d(d-1)}{2} + 3jd$, for $k=j+2d-2$, where $0 \leq j \leq d-1$.
  \item $J_{k} = S_{k}$, for $k\geq 3d-2$.
\end{enumerate}
\end{cor}

\begin{prop}
\em Let $2\leq d:=d_1=d_2=d_3$ be positive integers. Let $f_{1},f_{2},f_{3}\in K[x_{1},x_{2},x_{3}]$ be a regular sequence of homogeneous polynomials of degrees $d_{1},d_{2},d_{3}$. If $I=(f_{1},f_{2},f_{3})$ , $J=\Gin(I)$, the generic initial ideal with respect to the reverse lexicographic order, and $S/I$ has (SLP), then:
\[ J = ( x_{1}^{d-2}\{x_{1},x_{2}\}^{2}, x_{1}^{d-2j-1}x_{2}^{3j+1}, x_{1}^{d-2j-2}x_{2}^{3j+2}\; for\; 1\leq j \leq  \frac{d-3}{2} , \; x_{2}^{\frac{3d-1}{2}}, x_{3}x_{2}^{\frac{3d-3}{3}}, \]\[
x_{3}^{2j+1}x_{1}^{2j}x_{2}^{\frac{3d-3}{2}-3j}, \ldots,  x_{3}^{2j+1}x_{2}^{\frac{3d-3}{2}-j} , 
1\leq j \leq \frac{d-3}{2} \;
 ,x_{3}^{d-2+2j}\{x_{1},x_{2}\}^{d-j}, 1\leq j\leq d) \]
if $d$ is odd, or
\[ J = ( x_{1}^{d-2}\{x_{1},x_{2}\}^{2}, x_{1}^{d-2j-1}x_{2}^{3j+1}, x_{1}^{d-2j-2}x_{2}^{3j+2}\; for\; 1\leq j \leq  \frac{d-4}{2}, \; x_{1}x_{2}^{\frac{3d-4}{2}}, x_{2}^{\frac{3d-2}{2}},\; \] \[
x_{3}^{2j}x_{1}^{2j-1}x_{2}^{\frac{3d}{2}-3j}, \ldots, x_{3}^{2j}x_{2}^{\frac{3d-2}{2}-j},  1\leq j\leq\frac{d-2}{2}
 ,x_{3}^{d-2+2j}\{x_{1},x_{2}\}^{d-j}, 1\leq j\leq d) \]
if $d$ is even.
\end{prop}

\begin{proof}
We have $|J_{d}|=3$, hence 
$J_{d}=\{ x_{1}^{d-2}\{x_{1},x_{2}\}^{2} \}$, since $J$ is strongly stable and $x_3$ is strong Lefschetz for $S/J$. Therefore:
\[ 
Shad(J_{d}) = \{x_{1}^{d-2}\{x_{1},x_{2}\}^{3} ,x_{1}^{d-2}x_{3}\{x_{1},x_{2}\}^{2}\}.
\]
Now we have four possibilities to analyze: $d=2$, $d=3$, $d=4$ and $d\geq 5$.

$\mathbf{d=2}$. Using the formulas from Corrolary $2.3.2$ we have $|J_{3}|-|Shad(J_2)|=2$ so there are two
generators to add to $Shad(J_{2})$ to obtain $J_{3}$. Since $J$ is strongly stable and $x_3$ is a strong Lefschetz element for $S/J$ these new generators are $x_{3}^{2}x_{1}$ and $x_{3}^{2}x_{2}$. Therefore
\[ J_3 = \{\{x_{1},x_{2}\}^{3} ,x_{3}\{x_{1},x_{2}\}^{2}, x_{3}^{2}\{x_{1},x_{2}\}\}.\]
Since $|J_4|-|Shad(J_3)|=1$ there is only one generator to add to $Shad(J_3)$ and this is precisely $x_{3}^{4}$.
It follows $J_4=S_4$ and thus we cannot add new minimal generators of $J$ in degree $>5$.

$\mathbf{d=3}$. We have $|J_{4}| - |Shad(J_3)|= 2$ so there are two generators to add to $Shad(J_{3})$ to obtain $J_{4}$. Since $J$ is strongly stable and $x_3$ is a strong Lefschetz element for $S/J$ these new generators are $x_{2}^{4}$ and $x_{3}x_{2}^{3}$. Therefore
\[ J_4 = \{\{x_{1},x_{2}\}^{4} ,x_{3}\{x_{1},x_{2}\}^{3}\}. \]
Since $|J_5|-|Shad(J_4)| = 3$ there are three new monomials to add to $Shad(J_4)$ in order to obtain $J_5$. Since
$J$ is strongly stable and $x_3$ is a strong Lefschetz element for $S/J$ these new generators are $x_{3}^{3}\{x_{1},x_{2}\}^{2}$. Analogously, we must add two new monomial to $Shad(J_5)$ in order to obtain $J_6$ and these are $x_{3}^{5}\{x_{1},x_{2}\}$. Finally, we will add $x_{3}^{7}$ and thus we cannot add new minimal generators of $J$ in degree $>7$.

$\mathbf{d=4}$. We have $|J_{4}| - |Shad(J_3)|= 2$ so there are two generators to add to $Shad(J_{4})$ to obtain $J_{5}$. Since $J$ is strongly stable and $x_3$ is a strong Lefschetz element for $S/J$ these new generators are $x_{1}x_{2}^{4}$ and $x_{2}^{5}$. Therefore
\[ J_5 = \{\{x_{1},x_{2}\}^{5} ,x_{1}^{2}x_{3}\{x_{1},x_{2}\}^{2}\}. \]
Since $|J_6|-|Shad(J_5)| = 2$ there are two new monomials to add to $Shad(J_5)$ in order to obtain $J_6$ and using the usual argument these new monomials are $x_{3}^{2}x_{1}x_{2}^{2},x_{3}^{2}x_{2}^{4}$. It follows
\[ J_6 = \{\{x_{1},x_{2}\}^{6} , x_{3}\{x_{1},x_{2}\}^{5},  x_{3}^{2}\{x_{1},x_{2}\}^{4}\}. \]
Finally, we will add consequently $x_{3}^{4}\{x_{1},x_{2}\}^{3}$, $x_{3}^{6}\{x_{1},x_{2}\}^{2}$, 
$x_{3}^{8}\{x_{1},x_{2}\}$ and $x_{3}^{10}$.

Suppose now $d\geq 5$. We have $|J_{d+1}|-|Shad(J_{d})|=2$ so there are two generators to add to $Shad(J_{d})$ to obtain $J_{d+1}$. Since $J$ is strongly stable and $x_3$ is a strong Lefschetz element for $S/J$ these new generators are
$x_{1}^{d-3}x_{2}^{4}, x_{1}^{d-4}x_{2}^{5}$. It follows
\[ J_{d+1} = \{ x_{1}^{d-4}\{x_{1},x_{2}\}^{5}, x_{1}^{d-2}x_{3}\{x_{1},x_{2}\}^{3}, x_{1}^{d-2}x_{3}^{2}\{x_{1},x_{2}\}^{2}\}. \]
We prove by induction on $j$, with $1\leq j \leq \left\lfloor \frac{d-3}{2} \right\rfloor$ that
\[ J_{d+j} = Shad(J_{d+j-1})\cup \{ x_{1}^{d-2j-1}x_{2}^{3j+1}, x_{1}^{d-2j-2}x_{2}^{3j+2} \} = \]
\[ = \{x_{1}^{d-2j-2}\{x_{1},x_{2}\}^{3j+2}, x_{1}^{d-2j+1}x_{3}\{x_{1},x_{2}\}^{3j-1}, \ldots, x_{1}^{d-2}x_{3}^{j}\{x_{1},x_{2}\}^{2}\}, \]
the assertion being checked for $j=1$. Assume now that the statement is true for 
some $j \leq \left\lfloor \frac{d-3}{2} \right\rfloor$. Then
\[Shad(J_{d+j}) = \{x_{1}^{d-2j-2}\{x_{1},x_{2}\}^{3j+3}, x_{1}^{d-2j-2}x_{3}\{x_{1},x_{2}\}^{3j+2},  \ldots, x_{1}^{d-2}x_{3}^{j+1}\{x_{1},x_{2}\}^{2}\}.\]
Since $|J_{d+j+1}|-|Shad(J_{d+j})|=2$ we must add two generators to $Shad(J_{d+j})$ to obtain $J_{d+j+1}$. Using the
fact that $J$ is strongly stable and $x_3$ is a strong Lefschetz element for $S/J$ it follows that these new generators are
$x_{1}^{d-2j-3}x_{2}^{3j+4}, x_{1}^{d-2j-4}x_{2}^{3j+5}$, so the induction step is fulfilled.

We must consider now two possibilities.

1. $d$ is odd. We obtain
	\[ J_{\frac{3d-3}{2}} = \{x_{1}\{x_{1},x_{2}\}^{\frac{3d-5}{2}}, x_{1}^{3}x_{3}\{x_{1},x_{2}\}^{\frac{3d-11}{2}}, 
	\ldots, x_{1}^{d-2}x_{3}^{\frac{d-3}{2}}\{x_{1},x_{2}\}^{2}\}.\]
	Since $|J_{\frac{3d-1}{2}}|-|Shad(J_{\frac{3d-3}{2}})|=2$ there are two generators to add to
	$Shad(J_{\frac{3d-3}{2}})$ to obtain $J_{\frac{3d-1}{2}}$, and they must be $x_{2}^{\frac{3d-1}{2}}$,
	$x_{3}x_{2}^{\frac{3d-3}{3}}$ using the usual argument. Therefore,
	\[ J_{\frac{3d-1}{2}} = \{ \{x_{1},x_{2}\}^{\frac{3d-1}{2}}, x_{3}\{x_{1},x_{2}\}^{\frac{3d-3}{2}}, 
	\ldots, x_{1}^{d-2}x_{3}^{\frac{d-1}{2}}\{x_{1},x_{2}\}^{2}\}.\]
  Since $|J_{\frac{3d+1}{2}}|-|Shad(J_{\frac{3d-1}{2}})|=1$ we must add a $3$ new generators to 
  $Shad(J_{\frac{3d-1}{2}})$ to obtain $J_{\frac{3d+1}{2}}$ and since $J$ is strongly stable and $x_3$ is a strong Lefschetz element for $S/J$, they are
  $x_{1}^{2}x_{3}^{3}x_{2}^{\frac{3d-9}{2}},x_{1}x_{3}^{3}x_{2}^{\frac{3d-7}{2}},x_{3}^{3}x_{2}^{\frac{3d-5}{2}}$.
  
We prove by induction on $j$, with $1\leq j \leq  \frac{d-3}{2}$ that
\[J_{\frac{3d-1}{2}+j} = \{Shad(J_{\frac{3d-3}{2}+j})\}\cup \{x_{3}^{2j+1}x_{1}^{2j}x_{2}^{\frac{3d-3}{2}-3j}, \ldots,  x_{3}^{2j+1}x_{2}^{\frac{3d-3}{2}-j}\} =\]\[ = \{ \{x_{1},x_{2}\}^{\frac{3d-1}{2}+j}, x_{3}\{x_{1},x_{2}\}^{\frac{3d-3}{2}+j}, \ldots, x_{3}^{2j+1}\{x_{1},x_{2}\}^{\frac{3d-3}{2}-j},\] 
\[ x_{3}^{2j+2}x_{1}^{2j+3}\{x_{1},x_{2}\}^{\frac{3d-11}{2}-3j},\ldots,x_{1}^{d-2}x_{3}^{\frac{d-1}{2}+j}\{x_{1},x_{2}\}^{2}\}. \]
This assertion is proved for $j=1$. Assume the assertion is true for some $j< \frac{d-3}{2} $.
Since $|J_{\frac{3d+1}{2}+j}| - |Shad(J_{\frac{3d-1}{2}+j})| = 2j+3$ we must add $2j+3$ new generators to 
$Shad(J_{\frac{3d-1}{2}+j})$ in order to obtain $J_{\frac{3d+1}{2}+j}$. The usual argument implies that those new generators are $x_{3}^{2j+3}x_{1}^{2j+2}x_{2}^{\frac{3d-3}{2}-3j-3}, \ldots, x_{3}^{2j+3}x_{2}^{\frac{3d-3}{2}-j-1}$, which conclude the induction.

2. $d$ is even. We obtain
\[ J_{\frac{3d-4}{2}} = \{x_{1}^{2}\{x_{1},x_{2}\}^{\frac{3d-8}{2}}, x_{1}^{4}x_{3}\{x_{1},x_{2}\}^{\frac{3d-14}{2}}, 
\ldots, x_{1}^{d-2}x_{3}^{\frac{d-4}{2}}\{x_{1},x_{2}\}^{2}\}.\]
We have $|J_{\frac{3d-2}{2}}|-|Shad(J_{\frac{3d-4}{2}})|=2$ so we must add two new generators to 
$Shad(J_{\frac{3d-4}{2}})$ to obtain $J_{\frac{3d-2}{2}}$. Since $J$ is strongly stable and $x_3$ is a strong Lefschetz element for $S/J$, they are $x_{1}x_{2}^{\frac{3d-4}{2}}$ and $x_{2}^{\frac{3d-2}{2}}$, therefore
\[ J_{\frac{3d-2}{2}} = \{\{x_{1},x_{2}\}^{\frac{3d-2}{2}}, x_{1}^{2}x_{3}\{x_{1},x_{2}\}^{\frac{3d-8}{2}}, \ldots, x_{1}^{d-2}x_{3}^{\frac{d-2}{2}}\{x_{1},x_{2}\}^{2}\}.\]
Since $|J_{\frac{3d}{2}}|-|Shad(J_{\frac{3d-2}{2}})|=2$ we must add two new generators to $Shad(J_{\frac{3d-2}{2}})$
in order to obtain $J_{\frac{3d}{2}}$ and since $J$ is strongly stable and $x_3$ is a strong Lefschetz element for $S/J$, they are $x_{1}x_{3}^{2}x_{2}^{\frac{3d-6}{2}}$ and
$x_{3}^{2}x_{2}^{\frac{3d-4}{2}}$.

We prove by induction on $j$, with $1\leq j \leq  \frac{d-2}{2} $ that
\[ J_{\frac{3d-2}{2} + j} = Shad(J_{\frac{3d-4}{2}+j}) \cup \{ x_{3}^{2j}x_{1}^{2j-1}x_{2}^{\frac{3d}{2}-3j}, \ldots, x_{3}^{2j}x_{2}^{\frac{3d-2}{2}-j}\}.\]
The assertion has been proved for $j=1$ and suppose it is true for some \linebreak $j< \frac{d-2}{2} $.
Since $|J_{\frac{3d}{2} + j}| - |Shad(J_{\frac{3d-2}{2} + j})| = 2j+2$ we must add $2j+2$ generators to
$Shad(J_{\frac{3d-2}{2} + j})$ in order to obtain $J_{\frac{3d}{2} + j}$ and since $J$ is strongly stable and $x_3$ is a strong Lefschetz element for $S/J$ they must be 
$x_{3}^{2j+2}x_{1}^{2j+1}x_{2}^{\frac{3d}{2}-3j-3}, \ldots, x_{3}^{2j+2}x_{2}^{\frac{3d-4}{2}-j}$, which conclude
the induction.

Either if $d$ is even, either if $d$ is odd, we obtain
\[J_{2d-2} = \{\{x_{1},x_{2}\}^{2d-2},x_{3}\{x_{1},x_{2}\}^{2d-2},\ldots,x_{3}^{d-2}\{x_{1},x_{2}\}^{d}\} = \{\{x_{1},x_{2}\}^{d} \{x_{1},x_{2},x_{3}\}^{d-2}\}.\]
Since $|J_{2d-1}|-|Shad(J_{2d-2})|=d$ we must add $d$ new generators to $Shad(J_{2d-2})$ to obtain $J_{2d-1}$. But
$J$ is strongly stable and $x_3$ is a strong Lefschetz element for $S/J$, so we must add $x_{3}^{d}\{x_{1},x_{2}\}^{d-1}$, therefore $J_{2d-1}=\{\{x_{1},x_{2}\}^{d-1} \{x_{1},x_{2},x_{3}\}^{d}\}$. Using induction on $j\leq d$, we prove that
 \[ J_{2d-2+j} = Shad(J_{2d-3+j}) \cup \{x_{3}^{d-2+2j}\{x_{1},x_{2}\}^{d-j}\} = \{x_{1},x_{2}\}^{d-j}\{x_{1},x_{2},x_{3}\}^{d-2+2j}.\]
For $j=1$ we already proved. Suppose the assertion is true for $j<d$. We have $|J_{2d-1+j}|-|Shad(J_{2d-2+j})|=d-j$ so
we must add $d-j$ new monomials to $Shad(J_{2d-2+j})$ to obtain $J_{2d-1+j}$ and from the usual argument, these new monomials are $x_{3}^{d+2j}\{x_{1},x_{2}\}^{d-j-1}$. Finally, since $J_{3d-2}=S_{3d-2}$ we cannot add new minimal generators of $J$ in degree $>3d-2$.
\end{proof}

\begin{cor}
In the conditions above, the number of minimal generators of $J$ is 
$1 + \frac{d(d+1)}{2} + \frac{(d+1)^{2}}{4}$ when $d$ is odd, 
or $1 + \frac{d(d+1)}{2} + \frac{d(d+2)}{4}$ when $d$ is even.
\end{cor}

\begin{ex}
\begin{enumerate}
	\item Let $d_{1}=d_{2}=d_{3}=5$. Then
	\[ J = (x_{1}^{3}\{x_1,x_2\}^{2},\; x_1^{2}x_2^{4},\;x_1x_2^{5},\; x_2^{7},\; x_3 x_2^{6},\;
	        x_{2}^{3}x_{3}^{3}\{x_1,x_2\}^{2},\; \]
	\[ x_{3}^{5}\{x_1,x_2\}^{4},\;x_{3}^{7}\{x_1,x_2\}^{3},\; x_{3}^{9}\{x_1,x_2\}^{2}, \; x_{3}^{11}\{x_1,x_2\}, \;
	x_{3}^{13}).\]
	\item Let $d_{1}=d_{2}=d_{3}=6$. Then
	\[ J = (x_1^{4}\{x_{1},x_{2}\}^{2},\; x_{1}^{3}x_{2}^{4},\; x_{1}^{2}x_{2}^{5},\; x_{1}x_{2}^{7},\; x_{2}^{8},\;
	        x_2^{6}x_3^{2}\{x_1,x_2\},\;x_2^{3}x_3^{4}\{x_1,x_2\}^{3},\;\]
  \[ x_{3}^{6}\{x_{1},x_2\}^{5},\; x_{3}^{8}\{x_{1},x_2\}^{4},\; x_{3}^{10}\{x_{1},x_2\}^{3},\;
	         x_{12}^{6}\{x_{1},x_2\}^{2},\; x_{14}^{6}\{x_{1},x_2\},\; x_{16}^{6}) \]
\end{enumerate}
\end{ex}

\medskip
\begin{itemize}
\item Subcase $d_1=d_2<d_3$.
\end{itemize} 
\medskip

\begin{prop}
\em Let $2\leq d:=d_1=d_2<d_3$ be positive integers such that $d_1+d_2 > d_3 + 1$. The Hilbert function of the standard graded complete intersection $A=K[x_{1},x_{2},x_{3}]/I$, where $I$ is the ideal generated by $f_{1}$, $f_{2}$, $f_{3}$, with $f_i$ homogeneous polynomials of degree $d_i$, for all $i$, with $1\leq i\leq 3$, has the form:  
\begin{enumerate}
	\item $H(A,k) = \binom{k+2}{k}$, for $k\leq d-1$.
	\item $H(A,k) = \binom{d+1}{2} + \sum_{i=1}^{j}(d-i)$, for $k=j+d-1$, where $0\leq j \leq d_{3}-d$.
	\item $H(A,k) = \binom{d+1}{2} + \sum_{i=1}^{d_{3}-d}(d-i) + \sum_{i=1}^{j}(2d-d_{3}-2i)$, 
	      for $k=j + d_3 - 1$, where $0\leq j \leq \left\lfloor \frac{2d-d_{3}-1}{2} \right\rfloor$.
	\item $H(A,k) = H(A,d_{3}+2d-3-k)$,for $k\geq \left\lceil \frac{d_{3}+2d-3}{2} \right\rceil$.
\end{enumerate}
\end{prop}

\begin{proof}
It follows from \cite[Lemma 2.9(b)]{PV}.
\end{proof}

\begin{cor}
\em In the conditions of Proposition $2.3.6$, let $J=\Gin(I)$ be the generic initial ideal of $I$ with respect to the reverse lexicographic order. If we denote by $J_{k}$ the set of monomials of $J$ of degree $k$, then:
\begin{enumerate}
	\item $|J_{k}| = 0$, for $k\leq d-1$.
	\item $|J_{k}| = j(j+1)$, for $k = j + d - 1$ where $0 \leq j \leq d_{3}-d$.
	\item $|J_{k}| = d_{3}^{2}+d_{3}-d^{2}-d-2dd_{3} + j(2d_{3}-d) + \frac{3j(j+1)}{2}$,
	      for $k=j+d_{3}-1$, where $0\leq j \leq \left\lfloor \frac{2d-d_{3}-1}{2} \right\rfloor$.
	
	\item If $d_{3}$ is even then $|J_{k}| = \frac{4d^{2}+3d_{3}^{2}-4dd_{3}+4d}{8} + \frac{j(2d+d_{3})}{2} +
	      \frac{3j(j+1)}{2}$, for $k = j + \frac{2d+d_{3}-2}{2}$, where $0 \leq j \leq \frac{2d-d_{3}-2}{2}$.
	      	
	      If $d_{3}$ is odd then $|J_{k}| = \frac{3d_{3}^{2}+4d^{2}-4dd_{3}-3}{2} + \frac{j(2d+d_{3}-3)}{2} +
	      \frac{3j(j+1)}{2}$, for $k = j + \frac{2d+d_{3}-3}{2}$, where $0 \leq j \leq \frac{2d-d_{3}-1}{2}$.	      		
  \item $|J_{k}| =  3d^{2} - 2d + \frac{d_{3}(d_{3}+1)}{2} - 2dd_{3} + (4d-d_{3})j + j(j-1)$, for $k=j+2d-2$,
        where $0 \leq j \leq d_{3}-d$.        
  \item $|J_{k}| = \frac{(d+d_{3})(d+d_{3}-1)}{2} - \frac{d(d+1)}{2} + j(2d+d_{3})$, for $k=j+d_{3}+d-2$,
        where $0 \leq j \leq d-1$.	
  \item $J_{k} = S_{k}$, for $k\geq 3d-2$.
\end{enumerate}
\end{cor}

\begin{prop}
\em Let $2\leq d:=d_1=d_2<d_3$ be positive integers such that $2d>d_3+1$. Let $f_{1},f_{2},f_{3}\in K[x_{1},x_{2},x_{3}]$ be a regular sequence of homogeneous polynomials of degrees $d_{1},d_{2},d_{3}$. If $I=(f_{1},f_{2},f_{3})$, $J=\Gin(I)$, the generic initial ideal with respect to the reverse lexicographic order and
$S/I$ has (SLP), then if $d_3$ is even, we have:
	\[ J = (x_{1}^{d},x_{1}^{d-1}x_{2}, x_{1}^{d-j-1}x_{2}^{2j+1} \;for\; 1\leq j\leq d_{3}-d-1, \]
	\[ x_{1}^{2d-d_{3}-2j+1}x_{2}^{2d_{3}-2d+3j-2}, x_{1}^{2d-d_{3}-2j}x_{2}^{2d_{3}-2d+3j-1}\;for\;1\leq j\leq
	    \frac{2d-d_{3}}{2} , \]
	\[ x_{3}^{2j}x_{1}^{2j-1}x_{2}^{\frac{2d+d_{3}-2}{2}-3j}, \ldots, x_{3}^{2j}x_{2}^{\frac{2d+d_{3}-4}{2}-j} 
	\;for\;1\leq j\leq \frac{2d-d_{3}-2}{2}, \]
	\[ x_{3}^{2d-d_{3}-2+2j}x_{2}^{2d_{3}-2d+2-2j}\{x_{1},x_{2}\}^{2d-d_{3}+j-2},1\leq j\leq d_{3}-d, 
	   x_{3}^{d_{3}-2+2j}\{x_{1},x_{2}\}^{d-j},1\leq j\leq d).\]
	Otherwise, if $d_3$ is odd, we have
	\[ J = (x_{1}^{d},x_{1}^{d-1}x_{2}, x_{1}^{d-j-1}x_{2}^{2j+1} \;for\; 1\leq j\leq d_{3}-d-1, \]
	\[ x_{1}^{2d-d_{3}-2j+1}x_{2}^{2d_{3}-2d+3j-2}, x_{1}^{2d-d_{3}-2j}x_{2}^{2d_{3}-2d+3j-1}\;for\;1\leq j\leq
	    \frac{2d-d_{3}-1}{2} , \]
	\[ x_{2}^{\frac{2d+d_{3}-1}{2}}, x_{3}x_{2}^{\frac{2d+d_{3}-3}{2}}, x_{3}^{2j+1}x_{1}^{2j}x_{2}^{\frac{2d+d_{3}-3}{2}-3j}, \ldots, x_{3}^{2j+1}x_{2}^{\frac{2d+d_{3}-3}{2}-j}
	,1\leq j\leq \frac{2d-d_{3}-3}{2}, \]
	\[ x_{3}^{2d-d_{3}-2+2j}x_{2}^{2d_{3}-2d+2-2j}\{x_{1},x_{2}\}^{2d-d_{3}+j-2},1\leq j\leq d_{3}-d, x_{3}^{d_{3}-2+2j}\{x_{1},x_{2}\}^{d-j},1\leq j\leq d). \]	
\end{prop}

\begin{proof}
We note first that $d\geq 3$. Indeed, if $d=2$ then the condition $2d>d_{3}+1$ implies $d_{3}=2$ which is a contradiction. We have $|J_{d}|=2$, hence $J_{d}=x_{1}^{d-1}\{x_{1},x_{2}\}$, since $J$ is strongly stable. Therefore:
\[ Shad(J_{d}) = \{x_{1}^{d-1}\{x_{1},x_{2}\}^{2} ,x_{1}^{d-1}x_{3}\{x_{1},x_{2}\}\}. \]
Assume $d_3>d+1$. Since $|J_{d+1}| - |Shad(J_{d})| = 1$ we must add a new generator to $Shad(J_{d})$ in order to obtain $J_{d+1}$. On the other hand, $J$ is strongly stable and $x_3$ is a strong Lefschetz element for $S/J$ so this new generator is $x_{1}^{d-2}x_{2}^{3}$, therefore
\[ J_{d+1} = \{x_{1}^{d-2}\{x_{1},x_{2}\}^{3}, x_{1}^{d-1}x_{3}\{x_{1},x_{2}\}\}.\]
We prove by induction on $1\leq j\leq d_{3}-d-1$ that
\[J_{d+j} = Shad(J_{d-1+j})\cup \{x_{1}^{d-j-1}x_{2}^{2j+1}\} = \{x_{1}^{d-j-1}\{x_{1},x_{2}\}^{2j+1}, \ldots, x_{1}^{d-1}x_{3}^{j}\{x_{1},x_{2}\}\}.\]
The case $j=1$ was done. Suppose the assertion is true for some $j<d_{3}-d-1$. Then, since $|J_{d+j+1}|-|Shad(J_{d+j})|=1$
it follows that we must add one generator to $Shad(J_{d+j})$ in order to obtain $J_{d+j+1}$. Since $J$ is strongly stable and $x_3$ is a strong Lefschetz element for $S/J$,
this new generator must be $x_{1}^{d-j-2}x_{2}^{2j+3}$. In particular,
\[ J_{d_{3}-1} = \{x_{1}^{2d-d_{3}}\{x_{1},x_{2}\}^{2d_{3}-2d-1}, x_{1}^{2d-d_{3}+1}x_{3}\{x_{1},x_{2}\}^{2d_{3}-2d-3}, \ldots, x_{3}^{d_{3}-d-1}x_{1}^{d-1}\{x_{1},x_{2}\}\}, \]
which is the same formula as in the case $d_3=d+1$.

We need to consider several possibilities. First, suppose $d_3 = 2d-2$. We have $|J_{d_{3}}| - |Shad(J_{d_{3}-1})| = 2$ so we must add two generators to $Shad(J_{d_{3}-1})$ in order to obtain $J_{d_{3}} = J_{2d-2}$. But since $J$ is strongly
stable and $x_3$ is a strong Lefschetz element for $S/J$, these new generators are $x_{1}x_{2}^{2d-3}$ and $x_{2}^{2d-2}$.

Suppose now $d_3<2d-2$. Since $|J_{d_{3}}| - |Shad(J_{d_{3}-1})| = 2$ we must add two generators to $Shad(J_{d_{3}-1})$
in order to obtain $J_{d_{3}}$. $J$ strongly stable and $x_3$ is a strong Lefschetz element for $S/J$ force us to choose $x_{1}^{2d-d_{3}-1}x_{2}^{2d_{3}-2d+1}, x_{1}^{2d-d_{3}-2}x_{2}^{2d_{3}-2d+2}$, so
\[ J_{d_{3}} = \{x_{1}^{2d-d_{3}-2}\{x_{1},x_{2}\}^{2d_{3}-2d}, x_{1}^{2d-d_{3}}x_{3}\{x_{1},x_{2}\}^{2d_{3}-2d-1}, x_{1}^{2d-d_{3}+1}x_{3}^{2}\{x_{1},x_{2}\}^{2d_{3}-2d-3}\}. \]
We show by induction on $1\leq j\leq \left\lfloor \frac{2d-d_{3}+1}{2} \right\rfloor$ that
\[ J_{d_{3}-1+j} = Shad(J_{d_{3}-2+j})\cup \{ x_{1}^{2d-d_{3}-2j+1}x_{2}^{2d_{3}-2d+3j-2}, x_{1}^{2d-d_{3}-2j}x_{2}^{2d_{3}-2d+3j-1}\} = \]
\[ = \{ x_{1}^{2d-d_{3}-2j}\{x_{1},x_{2}\}^{2d_{3}-2d+3j-1}, x_{3}x_{1}^{2d-d_{3}-2j+2}\{x_{1},x_{2}\}^{2d_{3}-2d+3j-4}, \ldots,\]\[ x_{3}^{j}x_{1}^{2d-d_{3}}\{x_{1},x_{2}\}^{2d_{3}-2d-1},
x_{3}^{j+1}x_{1}^{2d-d_{3}+1}\{x_{1},x_{2}\}^{2d_{3}-2d-3},\ldots, x_{3}^{d_{3}-d-1+j}x_{1}^{d-1}\{x_{1},x_{2}\}\}.\]
We already done the case $j=1$. Suppose the assertion is true for some \linebreak $j<\left\lfloor \frac{2d-d_{3}+1}{2} \right\rfloor$. We have $|J_{d_{3}+j}|-|Shad(J_{d_{3}+j-1})|=2$ so we add two generators to $Shad(J_{d_{3}+j-1})$ and they must be $x_{1}^{2d-d_{3}-2j-1}x_{2}^{2d_{3}-2d+3j+1}, x_{1}^{2d-d_{3}-2j-2}x_{2}^{2d_{3}-2d+3j+2}$
from the usual argument.
In the following, we distinguish between two possibilities: $d$ is even or $d$ is odd. If $d$ is even, we get
\[ J_{\frac{2d+d_{3}-4}{2}} = \{ x_{1}^{2}\{x_{1},x_{2}\}^{\frac{2d+d_{3}-8}{2}}, x_{3}x_{1}^{4}\{x_{1},x_{2}\}^{\frac{2d+d_{3}-14}{2}}, \ldots, \]
\[ x_{3}^{\frac{2d-d_{3}-2}{2}}x_{1}^{2d-d_{3}}\{x_{1},x_{2}\}^{2d_{3}-2d-1},\ldots
,x_{3}^{\frac{d_{3}-3}{2}}x_{1}^{d-1}\{x_{2},x_{2}\} \}.\]
We have $|J_{\frac{2d+d_{3}-2}{2}}|-|Shad(J_{\frac{2d+d_{3}-4}{2}})|=2$ so we add two generators to $Shad(J_{\frac{2d+d_{3}-4}{2}})$ and they must be $x_{1}x_{2}^{\frac{2d+d_{3}-4}{2}}, x_{2}^{\frac{2d+d_{3}-2}{2}}$, so:
\[ J_{\frac{2d+d_{3}-2}{2}} =  \{\{x_{1},x_{2}\}^{\frac{2d+d_{3}-2}{2}}, x_{3}x_{1}^{2}\{x_{1},x_{2}\}^{\frac{2d+d_{3}-8}{2}}, \ldots, \]
\[,x_{3}^{\frac{2d-d_{3}}{2}}x_{1}^{2d-d_{3}}\{x_{1},x_{2}\}^{2d_{3}-2d-1},\ldots
,x_{3}^{\frac{d_{3}-2}{2}}x_{1}^{d-1}\{x_{2},x_{2}\}\}.\]
One can easily show by induction on  $1\leq j\leq \frac{2d-d_{3}-2}{2}$, if case, that
\[ J_{\frac{2d+d_{3}-2}{2}+j} = Shad(J_{\frac{2d+d_{3}-4}{2}+j})\cup \{ x_{3}^{2j}x_{1}^{2j-1}x_{2}^{\frac{2d+d_{3}-2}{2}-3j}, \ldots, x_{3}^{2j}x_{2}^{\frac{2d+d_{3}-4}{2}-j}\} = \]
\[ = \{ \{x_{1},x_{2}\}^{\frac{2d+d_{3}-2}{2}+j}, \ldots, x_{3}^{2j}\{x_{1},x_{2}\}^{\frac{2d+d_{3}-2}{2}-j},x_{3}^{2j+1}x_{1}^{2j+2}\{x_{1},x_{2}\}^{\frac{2d+d_{3}-8}{2}-3j},
 \ldots, \]\[
x_{3}^{\frac{2d-d_{3}}{2} +j}x_{1}^{2d-d_{3}}\{x_{1},x_{2}\}^{2d_{3}-2d-1}, \ldots, 
x_{3}^{\frac{d_{3}-2}{2} + j}x_{1}^{d-1}\{x_{1},x_{2}\} \}. \]
Indeed, if $j=1$, we have $|J_{\frac{2d+d_{3}}{2}+j}| - |Shad(J_{\frac{2d+d_{3}-2}{2}+j})| = 2$ so we must add two
monomials to $Shad(J_{\frac{2d+d_{3}-2}{2}+j})$ in order to obtain $J_{\frac{2d+d_{3}}{2}+j}$. Since $J$ is strongly stable and $x_3$ is a strong Lefschetz element for $S/J$, these new monomials are $x_{3}^{2}x_{1}x_{2}^{\frac{2d+d_{3}-2}{2}-3}$ and $x_{3}^{2}x_{2}^{\frac{2d+d_{3}-2}{2}-2}$. The induction step is similar.

If $d$ is odd, we have $ J_{\frac{2d+d_{3}-3}{2}} = \{x_{1}\{x_{1},x_{2}\}^{\frac{2d+d_{3}-5}{2}}, x_{3}x_{1}^{3}\{x_{1},x_{2}\}^{\frac{2d+d_{3}-11}{2}}, \ldots, $

$x_{3}^{\frac{2d-d_{3}-1}{2}}x_{1}^{2d-d_{3}}\{x_{1},x_{2}\}^{2d_{3}-2d-1},\ldots
,x_{3}^{\frac{d_{3}-3}{2}}x_{1}^{d-1}\{x_{2},x_{2}\} \}$.

Since $|J_{\frac{2d+d_{3}-1}{2}}|-|Shad(J_{\frac{2d+d_{3}-3}{2}})| = 2$ we add two generators to $Shad(J_{\frac{2d+d_{3}-3}{2}})$ in order to obtain $J_{\frac{2d+d_{3}-1}{2}}$. Since $J$ is strongly stable and $x_3$ is a strong Lefschetz element for $S/J$, these new monomials are $x_{2}^{\frac{2d+d_{3}-1}{2}}, x_{3}x_{2}^{\frac{2d+d_{3}-3}{2}}$, therefore:
\[ J_{\frac{2d+d_{3}-1}{2}} =  \{ \{x_{1},x_{2}\}^{\frac{2d+d_{3}-1}{2}}, x_{3}\{x_{1},x_{2}\}^{\frac{2d+d_{3}-3}{2}}, x_{3}^{2}x_{1}^{3}\{x_{1},x_{2}\}^{\frac{2d+d_{3}-11}{2}}, \ldots, \]
\[x_{3}^{\frac{2d-d_{3}+1}{2}}x_{1}^{2d-d_{3}}\{x_{1},x_{2}\}^{2d_{3}-2d-1},\ldots
,x_{3}^{\frac{d_{3}-1}{2}}x_{1}^{d-1}\{x_{2},x_{2}\} \}.\]
Assume $d_3 < 2d-3$ (otherwise $\frac{2d+d_3-1}{2}=2d-2$). $|J_{\frac{2d+d_{3}+1}{2}}| - |Shad(J_{\frac{2d+d_{3}-3}{2}})| = 3$ so we must add $3$ generators to $Shad(J_{\frac{2d+d_{3}-3}{2}})$ to 
obtain $J_{\frac{2d+d_{3}+1}{2}}$. The usual argument implies that they are
$ x_{1}^{2}x_{3}^{3}x_{2}^{\frac{2d+d_{3}-9}{2}}, x_{1}x_{3}^{3}x_{2}^{\frac{2d+d_{3}-7}{2}}, x_{3}^{3}x_{2}^{\frac{2d+d_{3}-5}{2}}$ and thus
\[ J_{\frac{2d+d_{3}+1}{2}} = \{\{x_{1},x_{2}\}^{\frac{2d+d_{3}+1}{2}}, x_{3}\{x_{1},x_{2}\}^{\frac{2d+d_{3}-1}{2}}, x_{3}^{2}\{x_{1},x_{2}\}^{\frac{2d+d_{3}-3}{2}}, x_{3}^{3}\{x_{1},x_{2}\}^{\frac{2d+d_{3}-5}{2}}, \ldots, \]
\[x_{3}^{\frac{2d-d_{3}+3}{2}}x_{1}^{2d-d_{3}}\{x_{1},x_{2}\}^{2d_{3}-2d-1},\ldots
,x_{3}^{\frac{d_{3}+1}{2}}x_{1}^{d-1}\{x_{2},x_{2}\}\}.\]

One can easily prove by induction on $1\leq j\leq \frac{2d-d_{3}-3}{2}$ that
\[ J_{\frac{2d+d_{3}-1}{2}+j} = Shad(J_{\frac{2d+d_{3}-3}{2}+j})\cup \{ x_{3}^{2j+1}x_{1}^{2j}x_{2}^{\frac{2d+d_{3}-3}{2}-3j}, \ldots, x_{3}^{2j+1}x_{2}^{\frac{2d+d_{3}-3}{2}-j} \} = \]
\[ = \{ \{x_{1},x_{2}\}^{\frac{2d+d_{3}-1}{2}+j},\ldots, x_{3}^{2j+1}\{x_{1},x_{2}\}^{\frac{2d+d_{3}-3}{2}-j},
x_{3}^{2j+2}x_{1}^{2j+3}\{x_{1},x_{2}\}^{\frac{2d+d_{3}-11}{2}-3j}, \]\[
x_{3}^{\frac{2d-d_{3}+1}{2} +j}x_{1}^{2d-d_{3}}\{x_{1},x_{2}\}^{2d_{3}-2d-1}, \ldots, 
x_{3}^{\frac{d_{3}-1}{2} + j}x_{1}^{d-1}\{x_{1},x_{2}\} \}. \]

In all cases above, we get: $J_{2d-2} = \{\{x_{1},x_{2}\}^{2d-2},x_{3}\{x_{1},x_{2}\}^{2d-3},\ldots,$ \linebreak
$x_{3}^{2d-d_{3}-2}\{x_{1},x_{2}\}^{d_{3}}
,x_{3}^{2d-d_{3}-1}x_{1}^{2d-d_{3}}\{x_{1},x_{2}\}^{2d_{3}-2d-1},\ldots,x_{3}^{d-2}x_{1}^{d-1}\{x_{1},x_{2}\}\}$.  

We have $|J_{2d-1}| - |Shad(J_{2d-2})| = 2d-d_{3}$ so we must add $2d-d_{3}$ new generators to $Shad(J_{2d-2})$ to
obtain $J_{2d-1}$. We get
$ J_{2d-1} = \{ \{x_{1},x_{2}\}^{2d-1}, x_{3}\{x_{1},x_{2}\}^{2d-2}, \ldots,$ \linebreak
$ x_{3}^{2d-d_{3}}\{x_{1},x_{2}\}^{d_{3}-1}
,x_{3}^{2d-d_{3}+1}x_{1}^{2d-d_{3}+1}\{x_{1},x_{2}\}^{2d_{3}-2d-3},\ldots,x_{3}^{d-1}x_{1}^{d-1}\{x_{1},x_{2}\}\}$.

One can easily show by induction of $1\leq j\leq d_{3}-d$ that
\[ J_{2d-2+j} = Shad(J_{2d-3+j})\cup \{ x_{3}^{2d-d_{3}-2+2j}x_{2}^{2d_{3}-2d+2-2j}\{x_{1},x_{2}\}^{2d-d_{3}+j-2}\} = \]
\[ = \{ \{x_{1},x_{2}\}^{2d-2+j},x_{3}\{x_{1},x_{2}\}^{2d-3+j},\ldots,x_{3}^{2d-d_{3}-2+2j}\{x_{1},x_{2}\}^{d_{3}-j},\]
\[x_{3}^{2d-d_{3}-1+2j}x_{1}^{2d-d_{3}+j}\{x_{1},x_{2}\}^{2d_{3}-2d-1-2j},\ldots,x_{3}^{d-2+j}x_{1}^{d-1}\{x_{1},x_{2}\}\}.\]
In particular, we get: $J_{d+d_{3}-2} = \{\{x_{1},x_{2}\}^{d+d_{3}-2},x_{3}\{x_{1},x_{2}\}^{d+d_{3}-3},\ldots ,x_{3}^{d_{3}-2}\{x_{1},x_{2}\}^{d}\}$. 
We have $|J_{d+d_3-1}|-|Shad(J_{d+d_3-2})|=d$ so we must add $d$ new generators to $Shad(J_{d+d_3-2})$ in order
to obtain $J_{d+d_3-1}$. Since $J$ is strongly stable and $x_3$ is a strong Lefschetz element for $S/J$, these new generators are $x_{3}^{d_{3}}\{x_{1},x_{2}\}^{d-1}$ so
\[J_{d+d_{3}-1} = \{\{x_{1},x_{2}\}^{d+d_{3}-1},x_{3}\{x_{1},x_{2}\}^{d+d_{3}-2},\ldots ,x_{3}^{d_{3}}\{x_{1},x_{2}\}^{d-1}\}. \]
Now, one can easily prove by induction on $1\leq j\leq d-1$ that $J_{d+d_{3}-2+j}$ is the set
\[ Shad(J_{d+d_{3}-3+j}) \cup \{x_{3}^{d_{3}-2+2j}\{x_{1},x_{2}\}^{d-j}\} =
 \{ \{x_{1},x_{2}\}^{d+d_{3}-2+j},\ldots ,x_{3}^{d_{3}-2+2j}\{x_{1},x_{2}\}^{d-j}\}.\]
Finally we obtain that $J_{d_{3}+2d-2}=S_{d_{3}+2d-2}$ and thus we cannot add new minimal generators of $J$ in degree $>d_{3}+2d-2$.
\end{proof}

\begin{cor}
In the conditions of the above proposition, the number of minimal generators of $J$ is $d(d+1) - \left(\frac{2d-d_3}{2}\right)^{2} + 1$ if $d_3$ is even or $d(d+1) - \frac{(2d-d_3)^{2}-1}{4} + 1$ if $d_3$ is odd.
\end{cor}

\begin{ex}
\begin{enumerate}
	\item Let $d_{1}=d_{2}=4$ and $d_{3}=6$. We have
	\[J = (x_{1}^{4},\; x_{1}^{3}x_2,\; x_1^{2}x_2^{3},\; x_1 x_2^{5},\;x_2^{6},\; x_3^{2}x_2^{4}\{x_1,x_2\},\; 
	       x_{3}^{4}x_{2}^{2}\{x_1,x_2\}^{2},\]\[\; x_3^{6}\{x_1,x_2\}^{3},\; x_3^{8}\{x_1,x_2\}^{2},\;
	       x_3^{10}\{x_1,x_2\},\; x_3^{12}). \]
	\item Let $d_{1}=d_{2}=4$ and $d_{3}=5$. We have:
	\[J = (x_{1}^{4},\; x_{1}^{3}x_2,\; x_1^{2}x_2^{3},\; x_1 x_2^{4},\; x_2^{6},\;x_2^{5}x_{3},\; 
	       x_3^{3}x_2^{2}\{x_1,x_2\}^{2},\;\]\[\; x_3^{5}\{x_1,x_2\}^{3},\; x_3^{7}\{x_1,x_2\}^{2},\;
	       x_3^{9}\{x_1,x_2\},\; x_3^{11} ). \]
\end{enumerate}
\end{ex}

\medskip
\begin{itemize}
\item Subcase $d_1<d_2=d_3$.
\end{itemize} 
\medskip

\begin{prop}
\em Let $2\leq d_1<d_2=d_3=:d$ be positive integers. The Hilbert function of the standard graded complete intersection $A=K[x_{1},x_{2},x_{3}]/I$, where $I$ is the ideal generated by $f_{1}$, $f_{2}$, $f_{3}$, with $f_i$ homogeneous polynomials of degree $d_i$, for all $i$, with $1\leq i\leq 3$, has the form:  
\begin{enumerate}
	\item $H(A,k) = \binom{k+2}{k}$, for $k\leq d_{1}-2$.
	\item $H(A,k) = \binom{d_{1}+1}{2} + jd_{1}$, for $k = j +d_{1} - $, where $0 \leq j \leq d-d_1$.
	\item $H(A,k) = \binom{d_{1}+1}{2} + d_{1}(d-d_{1}) + \sum_{i=1}^{j}(d_{1}-2i)$, 
        for $k = j + d - 1$, where $0 \leq j \leq \left\lfloor \frac{d_{1}-1}{2} \right\rfloor$.
	\item $H(A,k) = H(A,d_{1}+2d-3-k)$, for $k\geq \left\lceil \frac{d_{1}+2d-3}{2} \right\rceil$.
\end{enumerate}
\end{prop}

\begin{proof}
It follows from \cite[Lemma 2.9(b)]{PV}.
\end{proof}

\begin{cor}
\em In the conditions of Proposition $2.3.11$, let $J=\Gin(I)$ be the generic initial ideal of $I$ with respect to the reverse lexicographic order. If we denote by $J_{k}$ the set of monomials of $J$ of degree $k$, then:
\begin{enumerate}
	\item $|J_{k}| = 0$, for $k\leq d_{1}-1$.
	\item $|J_{k}| = j(j+1)/2$, for $k=j+d_{1}-1$, where $0 \leq j \leq d-d_{1}$.	
	\item $|J_{k}| = \frac{(d-d_{1})(d-d_{1}-1)}{2} + j(d-d_{1}) + \frac{3j(j+1)}{2}$,
	      for $k=j+d-1$, where $0\leq j \leq \left\lfloor \frac{d_{1}-1}{2} \right\rfloor$.	      
	\item If $d_{1}$ is even then
	      $|J_{k}| = \frac{3d_{1}^{2}+2d_{1}+4d^{2}+4d-4dd_{1}}{8} + \frac{j(2d+d_{1})}{2} + \frac{3j(j+1)}{2}$,
	      for $k = j + \frac{2d+d_{1}-2}{2}$, where $0\leq j \leq \frac{d_{1}-2}{2}$.
	      
	      If $d_{1}$ is odd then
	      $|J_{k}| = \frac{3d_{1}^{2}+4d^{2}-4dd_{1}-3}{2} + \frac{j(2d+d_{1})}{2} + \frac{3j^{2}}{2}$
	      for $k = j + \frac{2d+d_{1}-3}{2}$, where $0\leq j \leq \frac{d_{1}-1}{2}$.		
  \item $|J_{k}| = \frac{d(d-1)}{2} + d_{1}(d_{1}-1) + j(2d_{1}+d) + \frac{j(j-1)}{2}$, for
        $k=j +d_{1}+d -2$ , where $0  \leq j \leq d-d_{1}$.                 
  \item $|J_{k}| =  \frac{	2d(2d-1) - d_{1}(d_{1}-1)}{2} + j(2d+d_{1}))$, for
        $k=j+2d-2$, where $0 \leq j \leq d_{1}-1$.        
  \item $J_{k} = S_{k}$, for $k\geq 3d-2$.
\end{enumerate}
\end{cor}

\begin{prop}
\em Let $2\leq d_1<d_2=d_3=:d$ be positive integers. Let $f_{1},f_{2},f_{3}\in K[x_{1},x_{2},x_{3}]$ be a regular sequence of homogeneous polynomials of degrees $d_{1},d_{2},d_{3}$. Let $I=(f_{1},f_{2},f_{3})$, $J=\Gin(I)$, the generic initial ideal with respect to the reverse lexicographic order and
$S/I$ has (SLP), then if $d_1$ is even we have:
	\[J=(x_{1}^{d_{1}}, x_{1}^{d_{1}-2j+1}x_{2}^{d-d_{1}-2+3j}, x_{1}^{d_{1}-2j}x_{2}^{d-d_{1}-1+3j}\;for\; 
	     1\leq j \leq  \frac{d_{1}-2}{2}, \]
	\[  x_{1}x_{2}^{\frac{d_{1}+2d-4}{2}}, x_{2}^{\frac{d_{1}+2d-4}{2}},x_{3}^{2j}x_{1}^{2j-1}x_{2}^{\frac{d_{1}+2d}{2} -3j}, \ldots, x_{3}^{2j}x_{2}^{\frac{d_{1}+2d-2}{2}-j}
	    \;for\; 1 \leq j\leq \frac{d_{1}-4}{2},\]
	\[x_{3}^{d_{1}-2+2j}x_{2}^{d-d_{1}-j+1}\{x_{1},x_{2}\}^{d_{1}-1} \;for\;1\leq j \leq d-d_{1},\]	    
  \[  x_{3}^{2d-d_{1}-2+2j}\{x_{1},x_{2}\}^{d_{1}-j}\;for\; 1\leq j\leq d_{1}).\]
	Otherwise, if $d_1$ is odd, then:
	\[J=(x_{1}^{d_{1}}, x_{1}^{d_{1}-2j+1}x_{2}^{d-d_{1}-2+3j}, x_{1}^{d_{1}-2j}x_{2}^{d-d_{1}-1+3j}\;for\; 
	     1\leq j \leq  \frac{d_{1}-1}{2},  \]
	\[ x_{2}^{\frac{d_{1}+2d-1}{2}}, x_{3}x_{2}^{\frac{d_{1}+2d-3}{2}}, x_{3}^{2j+1}x_{1}^{2j}x_{2}^{\frac{d_{1}+2d-3}{2}
	    -3j}, \ldots, x_{3}^{2j+1}x_{2}^{\frac{d_{1}+2d-3}{2}-3j}\;for\; 1 \leq j\leq \frac{d_{1}-3}{2}, \]
	\[x_{3}^{d_{1}-2+2j}x_{2}^{d-d_{1}-j+1}\{x_{1},x_{2}\}^{d_{1}-1} \;for\;1\leq j \leq d-d_{1},\]	    
  \[  x_{3}^{2d-d_{1}-2+2j}\{x_{1},x_{2}\}^{d_{1}-j}\;for\; 1\leq j\leq d_{1}).\]	
\end{prop}

\begin{proof}
We have $|J_{d_1}|=1$, hence $J_{d_1}=\{x_{1}^{d_1}\}$, since $J$ is a strongly stable. Therefore:
\[ Shad(J_{d_1}) = \{x_{1}^{d_1}\{x_{1},x_{2}\} ,x_{1}^{d_1}x_{3}\}. \]
Assume $d>d_1+1$. Using the formulas from $2.3.12$ we get $|J_{d_{1}+1}|-|Shad(J_{d_{1}})|=0$ so $J_{d_{1}+1} = Shad(J_{d_{1}})$. We show by induction on $1\leq j\leq d-d_{1}$ that 
\[ J_{d_{1}+j} = Shad(J_{d_{1}-1+j}) = \{x_{1}^{d_{1}}\{x_{1},x_{2}\}^{j}, x_{1}^{d_{1}}x_{3}\{x_{1},x_{2}\}^{j-1}, \ldots, x_{1}^{d_{1}}x_{3}^{j}\}.\]
We already prove this for $j=1$. Suppose the assertion is true for some $j<d-d_1$. Since $|J_{d_{1}+j}|=|Shad(J_{d_{1}-1+j})|$ we have $J_{d_{1}+j} = Shad(J_{d_{1}-1+j})$ thus the induction step is done.
In particular, we get:
\[ J_{d-1} = \{x_{1}^{d_{1}}\{x_{1},x_{2}\}^{d-d_{1}-1}, x_{1}^{d_{1}}x_{3}\{x_{1},x_{2}\}^{d-d_{1}-2}, \ldots, x_{1}^{d_{1}}x_{3}^{d-d_{1}-1}\}. \]
which is the same expression as in the case $d=d_1+1$.

In the following, we consider two possibilities. First, suppose $d_1=2$. We have $|J_{d}|-|Shad(J_{d-1})|=2$ so we must add two new generators to $Shad(J_{d-1})$ in order to obtain $J_{d}$. Since $J$ is strongly stable and $x_3$ is a strong Lefschetz element for $S/J$ these new generators
are $x_1 x_2^{d-1}$ and $x_2^{d}$.

Suppose now $d_1>2$. We have $|J_{d}|-|Shad(J_{d-1})|=2$ so we must add
two new generators to $Shad(J_{d-1})$ in order to obtain $J_{d}$. Since $J$ is strongly stable and $x_3$ is a strong Lefschetz element for $S/J$ these new generators
are $x_{1}^{d_{1}-1}x_{2}^{d-d_{1}+1},x_{1}^{d_{1}-2}x_{2}^{d-d_{1}+2}$. Therefore
\[ J_{d} = \{x_{1}^{d_{1}-2}\{x_{1},x_{2}\}^{d-d_{1}}, x_{1}^{d_{1}}x_{3}\{x_{1},x_{2}\}^{d-d_{1}-1}, \ldots, x_{1}^{d_{1}}x_{3}^{d-d_{1}}\}.\]
We prove by induction on $1\leq j \leq \left\lfloor \frac{d_{1}-1}{2} \right\rfloor$ that:
\[ J_{d-1+j} = Shad(J_{d-2+j})\cup \{x_{1}^{d_{1}-2j+1}x_{2}^{d-d_{1}-2+3j}, x_{1}^{d_{1}-2j}x_{2}^{d-d_{1}-1+3j}\} = \]
\[ = \{ x_{1}^{d_{1}-2j}\{x_{1},x_{2}\}^{d-d_{1}-1-3j}, x_{1}^{d_{1}-2j+2}x_{3}\{x_{1},x_{2}\}^{d-d_{1}-4-3j},\ldots,
x_{1}^{d_{1}}x_{3}^{j}\{x_{1},x_{2}\}^{d-d_{1}-1},\]
\[ x_{1}^{d_{1}}x_{3}^{j+1}\{x_{1},x_{2}\}^{d-d_{1}-2},\ldots, x_{1}^{d_{1}}x_{3}^{d-d_{1}+j-1}\}. \]
We already done the case $j=1$. Suppose the assertion is true for some $j<\left\lfloor \frac{d_{1}-1}{2} \right\rfloor$. Since $|J_{d+j}|-|Shad(J_{d-1+j})|=2$ we must add two generators to $Shad(J_{d-1+j})$ in order
to obtain $|J_{d+j}|$ and they are $x_{1}^{d_{1}-2j-1}x_{2}^{d-d_{1}+3j+1}, x_{1}^{d_{1}-2j-2}x_{2}^{d-d_{1}+3j+2}$
because $J$ is strongly stable and $x_3$ is a strong Lefschetz element for $S/J$. Therefore, the induction step is done.

In the following, we consider two possibilities: $d_1$ is even or $d_1$ is odd. First, suppose $d_1$ is even. We have
\[ J_{\frac{d_{1}+2d-4}{2}} = \{x_{1}^{2}\{x_{1},x_{2}\}^{\frac{d_{1}+2d-8}{2}}, x_{1}^{4}x_{3}\{x_{1},x_{2}\}^{\frac{d_{1}+2d-14}{2}}, \ldots, x_{1}^{d_{1}}x_{3}^{\frac{d_{1}-2}{2}} \{x_{1},x_{2}\}^{d-d_{1}-1},\] \[x_{1}^{d_{1}}x_{3}^{\frac{d_{1}}{2}} \{x_{1},x_{2}\}^{d-d_{1}-2}, \ldots,
x_{1}^{d_{1}}x_{3}^{\frac{2d-d_{1}-4}{2}}\}. \]
Since $|J_{\frac{d_{1}+2d-2}{2}}| - |Shad(J_{\frac{d_{1}+2d-4}{2}})| = 2$ we need to add two new monomials to 
$Shad(J_{\frac{d_{1}+2d-4}{2}})$ and since $J$ is strongly stable and $x_3$ is a strong Lefschetz element for $S/J$, they are 
 $x_{1}x_{2}^{\frac{d_{1}+2d-4}{2}}, x_{2}^{\frac{d_{1}+2d-4}{2}}$, thus:
\[ J_{\frac{d_{1}+2d-2}{2}} =\{\{x_{1},x_{2}\}^{\frac{d_{1}+2d-2}{2}},x_{3}x_{1}^{2}\{x_{1},x_{2}\}^{\frac{d_{1}+2d-8}{2}}, x_{1}^{4}x_{3}^{2}\{x_{1},x_{2}\}^{\frac{d_{1}+2d-14}{2}}, \ldots, \] \[x_{1}^{d_{1}}x_{3}^{\frac{d_{1}}{2}} \{x_{1},x_{2}\}^{d-d_{1}-1},x_{1}^{d_{1}}x_{3}^{\frac{d_{1}+2}{2}} \{x_{1},x_{2}\}^{d-d_{1}-2}, \ldots,
x_{1}^{d_{1}}x_{3}^{\frac{2d-d_{1}-2}{2}}\}. \]

One can easily show by induction on $1 \leq j\leq \frac{d_{1}-2}{2}$ that
\[ J_{\frac{d_{1}+2d-2}{2} + j} = Shad(J_{\frac{d_{1}+2d-4}{2} + j}) \cup \{x_{3}^{2j}x_{1}^{2j-1}x_{2}^{\frac{d_{1}+2d}{2} -3j}, \ldots, x_{3}^{2j}x_{2}^{\frac{d_{1}+2d-2}{2}-j}\} = \]
\[ = \{ \{x_{1},x_{2}\}^{\frac{d_{1}+2d-2}{2} + j}, x_{3}\{x_{1},x_{2}\}^{\frac{d_{1}+2d-4}{2} + j}, \ldots, x_{3}^{2j}\{x_{1},x_{2}\}^{\frac{d_{1}+2d-2}{2} - j}, \]\[
x_{3}^{2j+1}x_{1}^{2j+2}\{x_{1},x_{2}\}^{\frac{d_{1}+2d-8}{2} - 3j}, \ldots, x_{3}^{\frac{d_{1}}{2}}x_{1}^{d_{1}}\{x_{1},x_{2}\}^{d-d_{1}-1},\ldots,x_{1}^{d_{1}}x_{3}^{\frac{2d-d_{1}-2}{2}+j}\}.\]
The assertion was already done for $j=1$ and the induction step is similar.

If $d_1$ is odd, we get
\[ J_{\frac{d_{1}+2d-3}{2}} = \{x_{1}\{x_{1},x_{2}\}^{\frac{d_{1}+2d-5}{2}}, x_{1}^{3}x_{3}\{x_{1},x_{2}\}^{\frac{d_{1}+2d-11}{2}}, \ldots, x_{1}^{d_{1}}x_{3}^{\frac{d_{1}-1}{2}} \{x_{1},x_{2}\}^{d-d_{1}-1},\] \[x_{1}^{d_{1}}x_{3}^{\frac{d_{1}+1}{2}} \{x_{1},x_{2}\}^{d-d_{1}-2}, \ldots,
x_{1}^{d_{1}}x_{3}^{\frac{2d-d_{1}-3}{2}}\}.\]

Since $|J_{\frac{d_{1}+2d-1}{2}}| - |Shad(J_{\frac{d_{1}+2d-3}{2}})|  = 2$ we add two generators to 
$Shad(J_{\frac{d_{1}+2d-3}{2}})$ in order to obtain $J_{\frac{d_{1}+2d-1}{2}}$ and they must be
$x_{2}^{\frac{d_{1}+2d-1}{2}}, x_{3}x_{2}^{\frac{d_{1}+2d-3}{2}}$, therefore:
\[ J_{\frac{d_{1}+2d-1}{2}} =\{\{x_{1},x_{2}\}^{\frac{d_{1}+2d-1}{2}}, x_{3}\{x_{1},x_{2}\}^{\frac{d_{1}+2d-3}{2}}, x_{1}^{3}x_{3}^{2}\{x_{1},x_{2}\}^{\frac{d_{1}+2d-11}{2}}, \ldots, \]\[ x_{1}^{d_{1}}x_{3}^{\frac{d_{1}+1}{2}} \{x_{1},x_{2}\}^{d-d_{1}-1},x_{1}^{d_{1}}x_{3}^{\frac{d_{1}+3}{2}} \{x_{1},x_{2}\}^{d-d_{1}-2}, \ldots,
x_{1}^{d_{1}}x_{3}^{\frac{2d-d_{1}-1}{2}}\}.\]
Since $|J_{\frac{d_{1}+2d+1}{2}}| - |Shad(J_{\frac{d_{1}+2d-1}{2}})|  = 3$, we add $3$ new generators to 
$Shad(J_{\frac{d_{1}+2d-1}{2}})$ in order to obtain $J_{\frac{d_{1}+2d+1}{2}}$ and they are
 $x_{1}^{2}x_{3}^{3}x_{2}^{\frac{d_{1}+2d-9}{2}}, x_{1}x_{3}^{3}x_{2}^{\frac{d_{1}+2d-7}{2}}, x_{3}^{3}x_{2}^{\frac{d_{1}+2d-5}{2}}$, therefore
\[ J_{\frac{d_{1}+2d+1}{2}} = \{\{x_{1},x_{2}\}^{\frac{d_{1}+2d+1}{2}}, x_{3}\{x_{1},x_{2}\}^{\frac{d_{1}+2d-1}{2}}, x_{3}^{2}\{x_{1},x_{2}\}^{\frac{d_{1}+2d-3}{2}},x_{3}^{3}\{x_{1},x_{2}\}^{\frac{d_{1}+2d-5}{2}}, \ldots,\] 
\[ x_{1}^{d_{1}}x_{3}^{\frac{d_{1}+3}{2}} \{x_{1},x_{2}\}^{d-d_{1}-1},x_{1}^{d_{1}}x_{3}^{\frac{d_{1}+5}{2}} \{x_{1},x_{2}\}^{d-d_{1}-2}, \ldots, x_{1}^{d_{1}}x_{3}^{\frac{2d-d_{1}+1}{2}}\}. \]

One can easily prove by induction on $1\leq j\leq \frac{d_{1}-3}{2}$ that
\[ J_{\frac{d_{1}+2d-1}{2} + j} = Shad(J_{\frac{d_{1}+2d-3}{2} + j})\cup \{ x_{3}^{2j+1}x_{1}^{2j}x_{2}^{\frac{d_{1}+2d-3}{2} -3j}, \ldots, x_{3}^{2j+1}x_{2}^{\frac{d_{1}+2d-3}{2}-3j}\} = \]
\[ = \{ \{x_{1},x_{2}\}^{\frac{d_{1}+2d-1}{2} + j}, x_{3}\{x_{1},x_{2}\}^{\frac{d_{1}+2d-3}{2} + j}, \ldots, x_{3}^{2j+2}\{x_{1},x_{2}\}^{\frac{d_{1}+2d-3}{2} - j},\]
\[ x_{3}^{2j+2}x_{1}^{2j+3}\{x_{1},x_{2}\}^{\frac{d_{1}+2d-11}{2} - 3j}, \ldots, x_{3}^{\frac{d_{1}+1}{2}}x_{1}^{d_{1}}\{x_{1},x_{2}\}^{d-d_{1}-1},\ldots,x_{1}^{d_{1}}x_{3}^{\frac{2d-d_{1}-1}{2}+j}\}.\]
The assertion was already proved for $j=1$ and the induction step is similar.

In all cases, we obtain
\[ J_{d_{1}+d-2} = \{ \{x_{1},x_{2}\}^{d_{1}+d-2}, x_{3}\{x_{1},x_{2}\}^{d_{1}+d-3}, \ldots, x_{3}^{d_{1}-2}\{x_{1},x_{2}\}^{d},\] \[ x_{1}^{d_{1}}x_{3}^{d_{1}-1}\{x_{1},x_{2}\}^{d-d_{1}-1},\ldots, x_{1}^{d_{1}}x_{3}^{d-2}\}. \]
We have $|J_{d_{1}+d-1}|-|Shad(J_{d_{1}+d-2})|=d_1$ so we must add $d_1$ new monomials to $Shad(J_{d_{1}+d-2})$
in order to obtain $J_{d_{1}+d-1}$. Since $J$ is strongly stable and $x_3$ is a strong Lefschetz element for $S/J$, these new monomials are 
 $x_{3}^{d_{1}}x_{2}^{d-d_{1}}\{x_{1},x_{2}\}^{d_{1}-1}$, therefore
\[ J_{d_{1}+d-1} = ( \{x_{1},x_{2}\}^{d_{1}+d-1}, x_{3}\{x_{1},x_{2}\}^{d_{1}+d-2}, \ldots, x_{3}^{d_{1}}\{x_{1},x_{2}\}^{d-1},\] \[ x_{1}^{d_{1}}x_{3}^{d_{1}+1}\{x_{1},x_{2}\}^{d-d_{1}-2},\ldots, x_{1}^{d_{1}}x_{3}^{d-1}).\]

One can easily prove by induction on $1\leq j \leq d-d_{1}$ that
\[ J_{d_{1}+d-1+j} = Shad(J_{d_{1}+d-2+j}) \cup \{ x_{3}^{d_{1}-2+2j}x_{2}^{d-d_{1}-j+1}\{x_{1},x_{2}\}^{d_{1}-1}\} = \]
\[ = \{ \{x_{1},x_{2}\}^{d_{1}+d-2+j}, x_{3}\{x_{1},x_{2}\}^{d_{1}+d-3+j}, \ldots, x_{3}^{d_{1}-2+2j}\{x_{1},x_{2}\}^{d-j},\] \[ x_{1}^{d_{1}}x_{3}^{d_{1}-1+2j}\{x_{1},x_{2}\}^{d-d_{1}-j-1},\ldots, x_{1}^{d_{1}}x_{3}^{d-2+j}\}  \]
the case $j=1$ being already done and than, the induction step being similar. In particular,
$J_{2d-2} = \{\{x_{1},x_{2}\}^{2d-2},x_{3}\{x_{1},x_{2}\}^{2d-3},\ldots ,x_{3}^{2d-d_{1}-2}\{x_{1},x_{2}\}^{d_{1}}\}$.

We have $|J_{2d-1}| -|Shad(J_{2d-2})|=d_1$ so we must add $d_1$ new generators to $Shad(J_{2d-2})$ in order to
obtain $J_{2d-1}$. Since $J$ is strongly stable and $x_3$ is a strong Lefschetz element for $S/J$, these new monomials are $x_{3}^{2d-d_{1}}\{x_{1},x_{2}\}^{d_{1}-1}$ so
\[J_{2d-1} = \{\{x_{1},x_{2}\}^{2d-1},x_{3}\{x_{1},x_{2}\}^{2d-2},\ldots ,x_{3}^{2d-d_{1}}\{x_{1},x_{2}\}^{d_{1}-1}\}. \]
One can easily show by induction on $1\leq j\leq d_{1}$ that
\[J_{2d-2+j} = Shad(J_{2d-3+j})\cup \{x_{3}^{2d-d_{1}-2+2j}\{x_{1},x_{2}\}^{d_{1}-j}\} = \]
 \[ = \{\{x_{1},x_{2}\}^{2d-2+j},\ldots,x_{3}^{2d-d_{1}-2+2j}\{x_{1},x_{2}\}^{d_{1}-j}\}.\]
Finally, we obtain $J_{d_{1}+2d-2}=S_{d_{1}+2d-2}$ and therefore we cannot add new minimal generators of $J$ in degrees $>d_{1}+2d-2$.
\end{proof}

\begin{cor}
In the conditions of the above proposition, the number of minimal generators of $J$ is $d_1(d+1) - 
\left(\frac{d_1}{2} \right)^{2} + 1$ if $d$ is even; $d_1(d+1) - \frac{d_1^{2}-1}{4} + 1$ if $d$ is odd.
\end{cor}

\begin{ex}
If $d_{1}=4$, $d_{2}=d_{3}=6$, then $J = (x_1^{4},\; x_1^{3}x_2^{3},\; x_1^{2}x_2^{4},\;x_2^{7},\; x_3x_2^{6},\;x_3^{2}x_{1}x_{2}^{5},$ \linebreak 
$ x_3^{2}x_{2}^{6},\; x_3^{4}x_{2}^{2}\{x_1,x_2\}^{3},\;  
	     x_3^{6}x_{2}\{x_1,x_2\}^{3},\; x_3^{8}\{x_1,x_2\}^{3},\; x_3^{10}\{x_1,x_2\}^{2},\; 
	       x_3^{12}\{x_1,x_2\},\; x_3^{14})$. \vspace{5 pt}

If $d_{1}=3$ and $d_{2}=d_{3}=6$, then: $J = (x_1^{3},\; x_1^{2}x_{2}^{4},\; x_1x_2^{5},\; x_1x_2^{6},\;x_2^{7},\;x_3^{3}x_2^{3}\{x_1,x_2\}^{2},\;$
\linebreak $ x_3^{5}x_2^{2}\{x_1,x_2\}^{2},\;   
	       x_3^{7}x_2\{x_1,x_2\}^{2},\;,x_3^{9}\{x_1,x_2\}^{2},\; x_3^{11}\{x_1,x_2\}^,\; x_3^{11}  )$.       
\end{ex}

\medskip
\begin{itemize}
\item Subcase $d_1<d_2<d_3$.
\end{itemize} 
\medskip

\begin{prop}
\em Let $2\leq d_1<d_2<d_3$ be positive integers such that \linebreak $d_{1}+d_{2}>d_{3}+1$. The Hilbert function of the standard graded complete intersection $A=K[x_{1},x_{2},x_{3}]/I$, where $I$ is the ideal generated by $f_{1}$, $f_{2}$, $f_{3}$, with $f_i$ homogeneous polynomials of degree $d_i$, for all $i$, with $1\leq i\leq 3$, has the form:  
\begin{enumerate}
	\item $H(A,k) = \binom{k+2}{k}$, for $k\leq d_{1}-2$.
	\item $H(A,k) = \binom{d_{1}+1}{2} + jd_{1}$, for $k=j+d_{1}-1$, where $0\leq j \leq d_{2}-d_{1}$.		
	\item $H(A,k) = \binom{d_{1}+1}{2} + d_{1}(d_{2}-d_{1}) + \sum_{i=1}^{j}(d_{1}-i)$, 
	      for $k=j+d_{2}-1$, where $0 \leq k \leq d_{3}-d_{2}$.	      
	\item $H(A,k) = \binom{d_{1}+1}{2} + d_{1}(d_{2}-d_{1}) + \sum_{i=1}^{d_{3}-1}(d_{1}-j) + \sum_{i=1}^{j}
	               (d_{1}+d_{2}-d_{3}-2i)$, for $k=j+d_{3}-1$, where $0\leq j \leq [\frac{d_{1}+d_{2}-d_{3}-1}{2}]$.
	\item $H(A,k) = H(A,d_{1}+d_{2}+d_{3}-3-k)$ for $k>d_{3}-1 + [\frac{d_{1}+d_{2}-d_{3}-1}{2}]$.
\end{enumerate}
\end{prop}

\begin{proof}
It follows from \cite[Lemma 2.9(b)]{PV}.
\end{proof}

\begin{cor}
\em In the conditions of Proposition $2.3.16$, let $J=\Gin(I)$ be the generic initial ideal of $I$ with respect to the reverse lexicographic order. If we denote by $J_{k}$ the set of monomials of $J$ of degree $k$, then:
\begin{enumerate}	
	\item $|J_{k}| = 0$, for $k\leq d_{1}-1$.
	\item $|J_{k}| = j(j+1)/2$, for $k = j +d_{1} - 1$, where $0 \leq j \leq d_{2}-d_{1}$.	
	\item $|J_{k}| = d_{2}(d_{2}-1) + j(d_{2}-d_{1}) + j(j+1)$,
	      for $k = j +d_{2} - 1$, where $0 \leq j \leq d_{3}-d_{2}$.	      
	\item $|J_{k}| = |J_{d_{3}-1}|+j(2d_{3}-d_{1}-d_{2})+\frac{3j(j+1)}{2}$, for $k = j+d_{3}-1$, 
	      where $0 \leq j \leq [\frac{\alpha-1}{2}]$.	      
	\item If $d_{1}+d_{2}+d_{3}$ is even then 
	      $|J_{k}| = |J_{\frac{d_{1}+d_{2}+d_{3}-4}{2}}|  + \frac{(j+1)(d_{1}+d_{2}+d_{3})}{2} + \frac{3j(j+1)}{2}$,
	      for $k = j + \frac{d_{1}+d_{2}+d_{3}-2}{2}$, where $0 \leq j \leq \frac{d_{1}+d_{2}-d_{3}-2}{2}$.
	
	      If $d_{1}+d_{2}+d_{3}$ is odd then 
	      $|J_{k}| = |J_{\frac{d_{1}+d_{2}+d_{3}-3}{2}}| +  \frac{j(d_{1}+d_{2}+d_{3})}{2} + \frac{3j^{2}}{2}$, 
	      for $k = j + \frac{d_{1}+d_{2}+d_{3}-3}{2}$, where $0 \leq j \leq \frac{d_{1}+d_{2}-d_{3}-1}{2}$.
  \item $|J_{k}| = |J_{d_{1}+d_{2}-2}| + j(2d_{1}+2d_{2}-d_{3}-1)+j^{2}$, 
        for $k=j+d_{1}+d_{2}-2$, where $0 \leq d_{3}-d_{2}$.       
  \item $|J_{k}| = |J_{d_{1}+d_{3}-2}|+j(2d_{1}+d_{3}-1)+\frac{j(j+1)}{2}$, for $k = j+d_{1}+d_{3}-2$,
        where $0 \leq j \leq d_{2}-d_{1}$.        
  \item $|J_{k}| = |J_{d_{2}+d_{3}-1}|+j(d_{1}+d_{2}+d_{3})$, for $k = j + d_{1}+d_{3}-2$,
        where $0 \leq j \leq d_{1}-1$        
  \item $J_{k} = S_{k}$, for $k\geq 3d-2$.
\end{enumerate}
\end{cor}

\begin{prop}
\em Let $2\leq d_1<d_2<d_3$ be positive integers such that $d_{1}+d_{2} > d_3+1$. Let $f_{1},f_{2},f_{3}\in K[x_{1},x_{2},x_{3}]$ be a regular sequence of homogeneous polynomials of degrees $d_{1},d_{2},d_{3}$. Let
$\alpha =d_{1}+d_{2}-d_{3}$. Let $I=(f_{1},f_{2},f_{3})$, $J=\Gin(I)$, the generic initial ideal with respect to the reverse lexicographic order, and suppose $S/I$ has (SLP). If $\alpha$ is even, then:
\[J = (x_{1}^{d_{1}},x_{1}^{d_{1}-1}x_{2}^{d_{2}-d_{1}+1},x_{1}^{d_{1}-2}x_{2}^{d_{2}-d_{1}+3},\ldots, 
        x_{1}^{d_{1}+d_{2}-d_{3}}x_{2}^{2d_{3}-d_{1}-d_{2}-1}, \] 
\[  x_{1}^{d_{1}+d_{2}-d_{3}-2j}x_{2}^{2d_{3}-d_{1}-d_{2}+3j-1},
x_{1}^{d_{1}+d_{2}-d_{3}-2j+1}x_{2}^{2d_{3}-d_{1}-d_{2}+3j-2} 
 for\;j=1,\ldots,\frac{\alpha-2}{2},\]
\[x_{2}^{\frac{d_{1}+d_{2}+d_{3}-2}{2}},x_{1}x_{2}^{\frac{d_{1}+d_{2}+d_{3}-4}{2}}, 
x_{3}^{2j}x_{2}^{\frac{d_{1}+d_{2}+d_{3}}{2}-3j} \{x_{1},x_{2}\}^{2j-1} for\;j=1,\ldots,\frac{\alpha-2}{2}, \]
\[ x_{1}^{d_{1}+d_{2}-d_{3}+j-2}x_{2}^{2d_{3}-d_{1}-d_{2}-2j+2}x_{3}^{d_{1}+d_{2}-d_{3}+2j-2},..,
  x_{2}^{d_{3}-j}x_{3}^{d_{1}+d_{2}-d_{3}+2j-2}\; for\;j=1,..,d_{3}-d_{2},  \]
\[  x_{1}^{d_{1}-1}x_{3}^{d_{1}+d_{3}-d_{2}+2j-2}x_{2}^{d_{2}-d_{1}-j+1} , \ldots, x_{3}^{d_{1}+d_{3}-d_{2}+2j-2}x_{2}^{d_{2}-j}\; for\; j=1,\ldots,d_{2}-d_{1}, \]
\[ \{x_{1},x_{2}\}^{d_{1}-j}x_{3}^{d_{2}+d_{3}-d_{1}-2+2j}\; for\; j=1,\ldots,d_{1}). \]

Otherwise, if $\alpha$ is odd, then:
\[J = (x_{1}^{d_{1}},x_{1}^{d_{1}-1}x_{2}^{d_{2}-d_{1}+1},x_{1}^{d_{1}-2}x_{2}^{d_{2}-d_{1}+3},\ldots, 
        x_{1}^{d_{1}+d_{2}-d_{3}}x_{2}^{2d_{3}-d_{1}-d_{2}-1},\]
\[ x_{1}^{d_{1}+d_{2}-d_{3}-2j}x_{2}^{2d_{3}-d_{1}-d_{2}+3j-1},
x_{1}^{d_{1}+d_{2}-d_{3}-2j+1}x_{2}^{2d_{3}-d_{1}-d_{2}+3j-2} ,
 j=1,\ldots,\frac{\alpha-1}{2}, \]
\[x_{2}^{\frac{d_{1}+d_{2}+d_{3}-1}{2}},x_{3}x_{2}^{\frac{d_{1}+d_{2}+d_{3}-3}{2}} ,x_{1}^{2j}x_{3}^{2j+1}x_{2}^{\frac{d_{1}+d_{2}+d_{3}-3}{2} -3j},.., x_{3}^{2j+1}x_{2}^{\frac{d_1+d_2+d_3-3}{2} - j} ,j=1,..,\frac{\alpha-3}{2}, \]
\[ x_{1}^{d_{1}+d_{2}-d_{3}+j-2}x_{2}^{2d_{3}-d_{1}-d_{2}-2j+2}x_{3}^{d_{1}+d_{2}-d_{3}+2j-2},..,
  x_{2}^{d_{3}-j}x_{3}^{d_{1}+d_{2}-d_{3}+2j-2}\; for\;j=1,..,d_{3}-d_{2},\]
\[   x_{1}^{d_{1}-1}x_{3}^{d_{1}+d_{3}-d_{2}+2j-2}x_{2}^{d_{2}-d_{1}-j+1} , \ldots, x_{3}^{d_{1}+d_{3}-d_{2}+2j-2}x_{2}^{d_{2}-j}\; for\; j=1,\ldots,d_{2}-d_{1}, \]
\[ \{x_{1},x_{2}\}^{d_{1}-j}x_{3}^{d_{2}+d_{3}-d_{1}-2+2j}\; for\; j=1,\ldots,d_{1}). \]
\end{prop}

\begin{proof}
We have $|J_{d_1}|=1$, hence $J_{d_1}=\{x_{1}^{d_1}\}$,
since $J$ is strongly stable. Therefore:
\[ Shad(J_{d_1}) = \{x_{1}^{d_1}\{x_{1},x_{2}\} ,x_{1}^{d_1}x_{3}\}. \]
Assume $d_2>d_1+1$. Using the formulas from $2.3.17$ we get $|J_{d_1+1}|-|Shad(J_{d_1})|=0$, therefore
\[J_{d_1+1} = Shad(J_{d_1}) = \{x_{1}^{d_1}\{x_{1},x_{2}\} ,x_{1}^{d_1}x_{3}\}. \]
Using induction on $1\leq j\leq d_2-d_1-1$ we prove that 
\[ J_{d_{1}+j} = Shad(J_{d_1+j-1}) = \{ x_{1}^{d_{1}}\{x_{1},x_{2}\}^{j} ,\ldots,x_{3}^{j}x_{1}^{d_{1}}\}. \]
Indeed, this assertion was already proved for $j=1$, and if we suppose that is true for some $j<d_2-d_1-1$ we
get $|J_{d_{1}+j+1}|=|Shad(J_{d_1+j})|$ so we are done. In particular, we obtain
$ J_{d_2-1} = \{ x_{1}^{d_{1}}\{x_{1},x_{2}\}^{d_{2}-d_{1}-1}, x_{3}x_{1}^{d_{1}}\{x_{1},x_{2}\}^{d_{2}-d_{1}-2},\ldots,x_{3}^{d_{2}-d_{1}-1}x_{1}^{d_{1}} \}.$

We have $|J_{d_2}|-|Shad(J_{d_2-1})|=1$ so we must add a new generator to $Shad(J_{d_2-1})$ in order to obtain
$J_{d_2}$. But since $J$ is strongly stable and $x_3$ is a strong Lefschetz element for $S/J$, this new generator is $x_{1}^{d_{1}-1}x_{2}^{d_{2}-d_{1}+1}$ and therefore 
$ J_{d_{2}} = \{ x_{1}^{d_{1}-1}\{x_{1},x_{2}\}^{d_{2}-d_{1}+1}, x_{3}x_{1}^{d_{1}}\{x_{1},x_{2}\}^{d_{2}-d_{1}-1},\ldots,x_{3}^{d_{2}-d_{1}}x_{1}^{d_{1}}\}. $

We show by induction on $1\leq j\leq d_3-d_2$ that 
\[ J_{d_2-1+j} = Shad(J_{d_2-2+j})\cup \{x_{1}^{d_{1}-j} x_{2}^{d_{2}-d_{1}+2j-1}\} = \]
\[= \{ x_{1}^{d_{1}-j}\{x_{1},x_{2}\}^{d_{2}-d_{1}+2j-1}, \ldots, x_{3}^{j}x_{1}^{d_{1}}\{x_{1},x_{2}\}^{d_{2}-d_{1}-1},\ldots,x_{3}^{d_{2}-d_{1}+j-1}x_{1}^{d_{1}}\}.\]
The first step of induction was already done. Suppose the assertion is true for some $j<d_{3}-d_{2}$. Since
$|J_{d_2+j}|-|Shad(J_{d_2-1+j})|=1$ and $J$ is strongly stable and $x_3$ is a strong Lefschetz element for $S/J$, we add $x_{1}^{d_{1}-j-1}x_{2}^{d_{2}-d_{1}+2j+1}$ 
to $Shad(J_{d_2-1+j})$ in order to obtain $J_{d_2+j}$. Thus, we are done. In particular, we get
\[ J_{d_{3}-1} = \{x_{1}^{d_{1}+d_{2}-d_{3}}\{x_{1},x_{2}\}^{2d_{3}-d_{1}-d_{2}-1}, x_{3}x_{1}^{d_{1}+d_{2}-d_{3}+1}\{x_{1},x_{2}\}^{2d_{3}-d_{1}-d_{2}-3}, \ldots,\] \[ x_{3}^{d_{3}-d_{2}}x_{1}^{d_{1}}\{x_{1},x_{2}\}^{d_{2}-d_{1}-1},\ldots,x_{3}^{d_{3}-d_{1}-1}x_{1}^{d_{1}}\}.\]

We consider first $\alpha=2$, i.e. $d_3 = d_1 + d_2 - 2$. In this case, since 
$|J_{d_3}|-|Shad(J_{d_3 -1})|=2$ we add two new generators to $Shad(J_{d_3-1})$ in order to obtain $J_{d_3}$.
Since $J$ is strongly stable and $x_3$ is a strong Lefschetz element for $S/J$, these new generators are $x_{2}^{d_{3}}$ and $x_{1}x_{2}^{d_{3}-1}$.

Suppose now $\alpha\geq 2$. We have $|J_{d_3}|-|Shad(J_{d_3 -1})|=2$ so we must add two new generators
to $Shad(J_{d_3-1})$ to obtain $J_{d_3}$. Since $J$ is strongly stable and $x_3$ is a strong Lefschetz element for $S/J$, these new generators are $x_{1}^{d_{1}+d_{2}-d_{3}-2}x_{2}^{2d_{3}-d_{1}-d_{2}+2}$ and
$x_{1}^{d_{1}+d_{2}-d_{3}-1}x_{2}^{2d_{3}-d_{1}-d_{2}+1}$, therefore
\[J_{d_{3}}= \{x_{1}^{d_{1}+d_{2}-d_{3}-2}\{x_{1},x_{2}\}^{2d_{3}-d_{1}-d_{2}+2}, x_{3}x_{1}^{d_{1}+d_{2}-d_{3}}\{x_{1},x_{2}\}^{2d_{3}-d_{1}-d_{2}-1},\]\[ x_{3}^{2}x_{1}^{d_{1}+d_{2}-d_{3}+1}\{x_{1},x_{2}\}^{2d_{3}-d_{1}-d_{2}-3},\ldots,
x_{3}^{d_{3}-d_{2}}x_{1}^{d_{1}}\{x_{1},x_{2}\}^{d_{2}-d_{1}},\ldots,x_{3}^{d_{3}-d_{1}}x_{1}^{d_{1}}\}.\]
One can prove by induction on $1\leq j\leq \left\lfloor \frac{\alpha-1}{2} \right\rfloor$ that $J_{d_{3}-1+j}$ is the set
\[ Shad(J_{d_{3}-2+j})\cup \{ x_{1}^{d_{1}+d_{2}-d_{3}-2j}x_{2}^{2d_{3}-d_{1}-d_{2}+3j-1}, x_{1}^{d_{1}+d_{2}-d_{3}-2j+1}x_{2}^{2d_{3}-d_{1}-d_{2}+3j-2} \} = \]
\[ \{x_{1}^{d_{1}+d_{2}-d_{3}-2j}\{x_{1},x_{2}\}^{2d_{3}-d_{1}-d_{2}+3j-1},
\ldots, x_{1}^{d_{1}+d_{2}-d_{3}-2}x_{3}^{j-1}\{x_{1},x_{2}\}^{2d_{3}-d_{1}-d_{2}+2},
x_{3}^{j}J_{d_{3}-1}\}.\]
Indeed, the assertion was proved for $j=1$ and the induction step is similar. In the following, we must consider
two possibilities: $\alpha$ is even or $\alpha$ is odd. Suppose first $\alpha$ is even. We obtain that $J_{\frac{d_{1}+d_{2}+d_{3}-4}{2}}$ is the set
\[ \{ x_{1}^{2}\{x_{1},x_{2}\}^{\frac{d_{1}+d_{2}+d_{3}-8}{2}}, \ldots,
x_{1}^{d_{1}+d_{2}-d_{3}-2}x_{3}^{\frac{d_{1}+d_{2}-d_{3}-4}{2}}
\{x_{1},x_{2}\}^{2d_{3}-d_{1}-d_{2}+2}, x_{3}^{\frac{d_{1}+d_{2}-d_{3}-2}{2}}J_{d_{3}-1}\}.\] 
We have $|J_{\frac{d_{1}+d_{2}+d_{3}-2}{2}}|-|Shad(J_{\frac{d_{1}+d_{2}+d_{3}-4}{2}})|=2$ so we must add two
generators to $Shad(J_{\frac{d_{1}+d_{2}+d_{3}-4}{2}})$ and, since $J$ is strongly stable and $x_3$ is a strong Lefschetz element for $S/J$, these new generators are
$x_{2}^{\frac{d_{1}+d_{2}+d_{3}-2}{2}}$ and $x_{1}x_{2}^{\frac{d_{1}+d_{2}+d_{3}-4}{2}}$. Therefore
\[ J_{\frac{d_{1}+d_{2}+d_{3}-2}{2}} = \{ \{x_{1},x_{2}\}^{\frac{d_{1}+d_{2}+d_{3}-2}{2}},
x_{1}^{2}x_{3}\{x_{1},x_{2}\}^{\frac{d_{1}+d_{2}+d_{3}-8}{2}},\ldots,\] 
\[  x_{1}^{d_{1}+d_{2}-d_{3}-2}x_{3}^{\frac{d_{1}+d_{2}-d_{3}-2}{2}}
\{x_{1},x_{2}\}^{2d_{3}-d_{1}-d_{2}+2},x_{3}^{\frac{d_{1}+d_{2}-d_{3}}{2}}J_{d_{3}-1}\}.\]

One can easily prove by induction on $1\leq j \leq \frac{\alpha-4}{2}$ that $J_{\frac{d_{1}+d_{2}+d_{3}-2}{2} + j}$
is the set
\[ = Shad(J_{\frac{d_{1}+d_{2}+d_{3}-4}{2} + j})\cup 
\{ x_{3}^{2j}x_{1}^{2j-1} x_{2}^{\frac{d_{1}+d_{2}+d_{3}-2}{2}-3j},\ldots,x_{3}^{2j}x_{2}^{\frac{d_{1}+d_{2}+d_{3}-2}{2}-j} \} = \]
\[\{ \{x_{1},x_{2}\}^{\frac{d_{1}+d_{2}+d_{3}-2}{2}+j},\ldots,x_{3}^{2j} \{x_{1},x_{2}\}^{\frac{d_{1}+d_{2}+d_{3}-2}{2}-j},x_{1}^{2j+2}x_{3}^{2j+1}\{x_{1},x_{2}\}^{\frac{d_{1}+d_{2}+d_{3}-8}{2}-3j},\ldots,\]
\[ x_{1}^{d_{1}+d_{2}-d_{3}-2}x_{3}^{\frac{d_{1}+d_{2}-d_{3}-2}{2}+j}\{x_{1},x_{2}\}^{2d_{3}-d_{1}-d_{2}+2},
  x_{3}^{\frac{d_{1}+d_{2}-d_{3}}{2}+j}J_{d_{3}-1} \}.\]

Suppose now $\alpha$ is odd. We have that $J_{\frac{d_{1}+d_{2}+d_{3}-3}{2}}$ is the set
\[\{ x_{1}\{x_{1},x_{2}\}^{\frac{d_{1}+d_{2}+d_{3}-5}{2}}, \ldots,
x_{1}^{d_{1}+d_{2}-d_{3}-2}x_{3}^{\frac{d_{1}+d_{2}-d_{3}-3}{2}}\{x_{1},x_{2}\}^{2d_{3}-d_{1}-d_{2}+2},x_{3}^{\frac{d_{1}+d_{2}-d_{3}-1}{2}}J_{d_{3}-1}\}.\]
We have $|J_{\frac{d_{1}+d_{2}+d_{3}-1}{2}}|-| Shad(J_{\frac{d_{1}+d_{2}+d_{3}-3}{2}})|=2$ so we must add two
new generators to $Shad(J_{\frac{d_{1}+d_{2}+d_{3}-3}{2}})$ in order to obtain $J_{\frac{d_{1}+d_{2}+d_{3}-1}{2}}$.
Since $J$ is strongly stable and $x_3$ is a strong Lefschetz element for $S/J$, these new generators are $x_{2}^{\frac{d_{1}+d_{2}+d_{3}-1}{2}}$ and $x_{3}x_{2}^{\frac{d_{1}+d_{2}+d_{3}-3}{2}}$. Therefore
\[J_{\frac{d_{1}+d_{2}+d_{3}-1}{2}} =  \{ \{x_{1},x_{2}\}^{\frac{d_{1}+d_{2}+d_{3}-1}{2}},
x_{3}\{x_{1},x_{2}\}^{\frac{d_{1}+d_{2}+d_{3}-3}{2}}, x_{1}^{3}x_{3}^{2}\{x_{1},x_{2}\}^{\frac{d_{1}+d_{2}+d_{3}-11}{2}}, \ldots, \] 
\[ x_{1}^{d_{1}+d_{2}-d_{3}-2}x_{3}^{\frac{d_{1}+d_{2}-d_{3}-1}{2}}\{x_{1},x_{2}\}^{2d_{3}-d_{1}-d_{2}+2},
x_{3}^{\frac{d_{1}+d_{2}-d_{3}+1}{2}}J_{d_{3}-1} \}.\]

One can easily prove by induction on $1\leq j\leq \frac{\alpha-3}{2}$ that
\[ J_{\frac{d_{1}+d_{2}+d_{3}-1}{2} + j} = Shad(J_{\frac{d_{1}+d_{2}+d_{3}-3}{2} + j})\cup\{ 
x_{1}^{2j}x_{3}^{2j+1}x_{2}^{\frac{d_{1}+d_{2}+d_{3}-3}{2} -3j},\ldots, x_{3}^{2j+1}x_{2}^{\frac{d_{1}+d_{2}+d_{3}-3}{2} - j}\}. \]
For $j=1$, we notice that $|J_{\frac{d_{1}+d_{2}+d_{3}+1}{2}}| - |Shad(J_{\frac{d_{1}+d_{2}+d_{3}-1}{2}})|=3$ so we
must add $3$ new monomials to $Shad(J_{\frac{d_{1}+d_{2}+d_{3}-1}{2}})$ in order to obtain $J_{\frac{d_{1}+d_{2}+d_{3}+1}{2}}$. But, since $J$ is strongly stable and $x_3$ is a strong Lefschetz element for $S/J$, they are exactly 
$x_{1}^{2}x_{3}^{3}x_{2}^{\frac{d_{1}+d_{2}+d_{3}-3}{2} -3}$, $x_{1}x_{3}^{3}x_{2}^{\frac{d_{1}+d_{2}+d_{3}-3}{2}-2}$ and $x_{3}^{3}x_{2}^{\frac{d_{1}+d_{2}+d_{3}-3}{2}-1}$, as required. The induction step is similar.

In all cases, we obtain that $J_{d_{1}+d_{2}-2}$ is the set
\[ \{ \{x_{1},x_{2}\}^{d_{1}+d_{2}-2},x_{3}\{x_{1},x_{2}\}^{d_{1}+d_{2}-3},\ldots, x_{3}^{d_{1}+d_{2}-d_{3}-2}\{x_{1},x_{2}\}^{d_{3}},x_{3}^{d_{1}+d_{2}-d_{3}-1}J_{d_{3}-1}\}.\]
We have $|J_{d_{1}+d_{2}-1}|-|Shad(J_{d_{1}+d_{2}-2})| = \alpha$, so we must add $\alpha$ new generators
to obtain $J_{d_{1}+d_{2}-1}$. Since $J$ is strongly stable and $x_3$ is a strong Lefschetz element for $S/J$, they are
$x_{1}^{d_{1}+d_{2}-d_{3}-1}x_{2}^{2d_{3}-d_{1}-d_{2}}x_{3}^{d_{1}+d_{2}-d_{3}}$, $\ldots, 
x_{2}^{d_{3}-1}x_{3}^{d_{1}+d_{2}-d_{3}}$, therefore $J_{d_{1}+d_{2}-1}$ is the set
\[ \{ \{x_{1},x_{2}\}^{d_{1}+d_{2}-1},\ldots,x_{3}^{d_{1}+d_{2}-d_{3}}\{x_{1},x_{2}\}^{d_{3}-1},
x_{3}^{d_{1}+d_{2}-d_{3}+1}x_{1}^{d_{1}+d_{2}-d_{3}+1}\{x_{1},x_{2}\}^{2d_{3}-d_{1}-d_{2}-3},\] 
\[\ldots,x_{3}^{d_{1}}x_{1}^{d_{1}} \{x_{1},x_{2}\}^{d_{2}-d_{1}-1}, 
x_{1}^{d_{1}}x_{3}^{d_{1}+1} \{x_{1},x_{2}\}^{d_{2}-d_{1}-2},\ldots,x_{1}^{d_{1}}x_{3}^{d_{2}-1}\}.\]

One can easily prove by induction on $1\leq j \leq d_{3}-d_{2}$ that $J_{d_{1}+d_{2}-1+j}$ is the set
\[ Shad(J_{d_{1}+d_{2}-2+j})\cup \{
 x_{1}^{d_{1}+d_{2}-d_{3}+j-1}x_{2}^{2d_{3}-d_{1}-d_{2}-2j}x_{3}^{d_{1}+d_{2}-d_{3}+2j}, \ldots,
  x_{2}^{d_{3}-1-j}x_{3}^{d_{1}+d_{2}-d_{3}+2j}\} \]
Indeed, the case $j=1$ was already done and the induction step is similar. In particular, we get
\[J_{d_{1}+d_{3}-2} = \{\{x_{1},x_{2}\}^{d_{1}+d_{3}-2}, x_{3}\{x_{1},x_{2}\}^{d_{1}+d_{3}-3},\ldots,
                         x_{3}^{d_{1}+d_{3}-d_{2}-2}\{x_{1},x_{2}\}^{d_{2}},\]
\[x_{1}^{d_{1}}x_{3}^{d_{1}+d_{3}-d_{2}-1}\{x_{1},x_{2}\}^{d_{2}-d_{1}-1}, \ldots, x_{1}^{d_{1}}x_{3}^{d_{3}-2}\}.\]
Since $|J_{d_{1}+d_{3}-1}|-|Shad(J_{d_{1}+d_{3}-2})|=d_{1}$ we must add $d_1$ generators to $Shad(J_{d_{1}+d_{3}-2})$
in order to obtain $J_{d_{1}+d_{3}-1}$. Since $J$ is strongly stable and $x_3$ is a strong Lefschetz element for $S/J$, these new generators are 
 $x_{1}^{d_{1}-1}x_{3}^{d_{1}+d_{3}-d_{2}}x_{2}^{d_{2}-d_{1}},\ldots,x_{3}^{d_{1}+d_{3}-d_{2}}x_{2}^{d_{2}-1}$,so 
 $J_{d_{1}+d_{3}-1}$ is the set
\[\{ \{x_{1},x_{2}\}^{d_{1}+d_{3}-1},...,x_{3}^{d_{1}+d_{3}-d_{2}}\{x_{1},x_{2}\}^{d_{2}-1},
x_{1}^{d_{1}}x_{3}^{d_{1}+d_{3}-d_{2}+1}\{x_{1},x_{2}\}^{d_{2}-d_{1}-2},...,x_{1}^{d_{1}}x_{3}^{d_{3}-1}\}.\]

We prove by induction on $1\leq j\leq d_2-d_1$ that $J_{d_{1}+d_{3}-2+j}$ is the set
\[  Shad(J_{d_{1}+d_{3}-3+j})\cup\{ x_{1}^{d_{1}-1}x_{3}^{d_{1}+d_{3}-d_{2}+2j-2}x_{2}^{d_{2}-d_{1}-j+1} , \ldots, x_{3}^{d_{1}+d_{3}-d_{2}+2j-2}x_{2}^{d_{2}-j}\}. \]
Indeed, we already proved this for $j=1$ and the induction step is similar. We get
\[J_{d_{2}+d_{3}-2} = \{ \{x_{1},x_{2}\}^{d_{1}+d_{3}-2},x_{3}\{x_{1},x_{2}\}^{d_{1}+d_{3}-3},\ldots,
x_{3}^{d_{3}+d_{2}-d_{1}-2}\{x_{1},x_{2}\}^{d_{1}} \}.\] 
One can easily prove by induction on $1\leq j\leq d_1$ that 
\[J_{d_{2}+d_{3}-2+j} = Shad(J_{d_{2}+d_{3}-3+j})\cup \{
x_{3}^{d_{2}+d_{3}-d_{1}-2+2j}\{x_{1},x_{2}\}^{d_{1}-j}\} =\]\[ =  \{ \{x_{1},x_{2}\}^{d_{2}+d_{3}-2+j},\ldots, x_{3}^{d_{2}+d_{3}-d_{1}-2+2j}\{x_{1},x_{2}\}^{d_{1}-j}\}.\]
Finally, we obtain $J_{d_{1}+d_2+d_3-2}=S_{d_{1}+d_2+d_3-2}$ and therefore we cannot add new minimal generators of $J$ in degrees $>d_{1}+d_2+d_3-2$.
\end{proof}

\begin{cor}
In the above conditions of the above proposition, the number of minimal generators of $J$ is $d_1(d_2+1) - 
\left(\frac{\alpha}{2} \right)^{2} + 1$ if $\alpha$ is even or $d_1(d_2+1)-\frac{\alpha^{2}-1}{4} + 1$ if $\alpha$ is odd.
\end{cor}

\begin{ex} 

\begin{enumerate}
	\item Let $d_{1}=3$ , $d_{2}=5$ and $d_{3}=6$. Then
\[ J = (x_1^{3}, x_1^{2}x_2^{3}, x_1x_2^{5},x_2^{6},\;x_2^{4}x_3^{2}\{x_1,x_2\},\;x_2^{2}x_3^{4}\{x_1,x_2\}^{2},\;\]\[ x_2x_3^{6}\{x_1,x_2\}^{2},\;x_3^{8}\{x_1,x_2\}^{2},\; x_3^{10}\{x_1,x_2\},\; x_3^{12}).\]
  \item Let $d_{1}=4$ , $d_{2}=5$ and $d_{3}=6$. Then 
  \[J = (x_1^{4},\;x_1^{3}x_2^{2},\;x_1^{2}x_2^{4},\;x_1x_2^{5},\;x_2^{7},\;x_3x_2^{6},
  \;x_2^{3}x_3^{3}\{x_1,x_2\}^{2},\]\[
\; x_2x_3^{5}\{x_1,x_2\}^{3},\; x_3^{7}\{x_1,x_2\}^{3},\; x_3^{9}\{x_1,x_2\}^{2},\; x_3^{11}\{x_1,x_2\},\;x_3^{13}).\]
\end{enumerate}
\end{ex}

\begin{rem}
If $f_1,f_2,f_3\in S=K[x_1,x_2,x_3]$ is a regular sequence of homogeneous polynomials of given degrees $d_1,d_2,d_3$ such that $S/(f_1,f_2,f_3)$ has $(SLP)$, then the number of minimal generators of $J=Gin((f_1,f_2,f_3))$, $\mu(J) \leq d_1(d_2+1)+1$. This follows immediately from $2.2.4$, $2.2.9$, $2.3.4$, $2.3.9$, $2.3.14$ and $2.3.19$.
\end{rem}

\begin{rem}
Let $f_1,f_2,f_3\in S=K[x_1,x_2,x_3]$ be a regular sequence of homogeneous polynomials of given degrees $d_1,d_2,d_3$
such that $S/I$ has (SLP), where $I=(f_1,f_2,f_3)$. Let $J=Gin(I)$. One can compute the socle of $S/J$, as follows. Since $J$ is (strongly) stable, $(J:(x_1,x_2,x_3))=(J:x_3)$. Indeed, if $u\in (J:x_3)$, then $x_3u\in J$ and
since $J$ is stable, $x_1(x_3u)/x_3 = x_1u\in J$ and also $x_2(x_3u)/x_3 = x_2u\in J$, thus $u\in (J:(x_1,x_2,x_3))$.

On the other hand, for example, when $d_1+d_2\leq d_3+1$ and $d_1=d_2<d_3$, Proposition $2.2.3$ implies $(J:x_3) = J + T$, where \[ T=(x_{2}^{2d-2j-2}x_{3}^{d_{3}-2d+2j+1}\{x_{1},x_{2}\}^{j},\;0\leq j\leq d-2,\; x_{3}^{d_{3}+2j-3}\{x_{1},x_{2}\}^{d-j},\; 1\leq j\leq d).\] One can check that none of the minimal generators of
$T$ is in $J$. Therefore, the set of the minimal generators of $T$ form a base for $Soc(S/J)$. 

The other cases are similar, and the reader can easily compute the socle of $S/J$ for any integers $d_1,d_2,d_3\geq 2$.
\end{rem}

\section{Generic initial ideal for $(n,d)$-complete intersections.}

Let $K$ be an algebraically closed field of characteristic zero. Let $S=K[x_{1},\ldots,x_{n}]$ be the polynomial ring in $n$ variables over $K$. Let $n,d\geq 2$ be two integers. We consider \linebreak $I=(f_{1},\ldots,f_{n})\subset S$ an ideal generated by a regular sequence $f_{1},\ldots,f_{n}\in S$ of homogeneous polynomials of degree $d$. We say that $A=S/I$ is a \emph{$(n,d)$-complete intersection}. Let $J=Gin(I)$ be the generic initial of $I$, with respect to the reverse lexicographic (revlex) order (see \cite[\S 15.9]{E}, for details). 

We say that a property $(P)$ holds for a generic sequence of homogeneous polynomials 
$f_1,f_2,\ldots, f_n \in S$ of given degrees $d_1,d_2,\ldots,d_n$ if there exists a nonempty open Zariski subset $U\subset S_{d_1}\times S_{d_2}\times \cdots \times S_{d_n}$ such that for every $n$-tuple $(f_{1},f_{2},\ldots,f_{n})\in U$ the property $(P)$ holds.

For any nonnegative integer $k$, we denote by $J_{k}$ the set of monomials of $J$ of degree $k$. Conca and Sidman proved that $J_d$ is revlex if $f_{1},\ldots,f_{n}$ is a generic regular sequence, (see \cite[Theorem 1.2]{CSi}). We will prove that $J_{d}$ is a revlex set in another case, namely, when $f_{i}\in k[x_{i},\ldots,x_{n}]$. It is likely to be true that $J_{d}$ is revlex for any complete intersection, but we do not have the means to prove this assertion.

Let $I=(f_{1},\ldots,f_{n})\subset S = K[x_1,\ldots,x_n]$ be an
ideal generated by a regular sequence $f_{1},\ldots,f_{n}\in S$ of
homogeneous polynomials of degree $d$. Let $J=Gin(I)$ be the generic
initial ideal of $I$, with respect to the revlex order. It is well
known that the Hilbert series of $S/J$ is the same as the Hilbert
series of $S/I$ and moreover, $H(S/J,t) = H(S/I,t) = (1+t + \cdots +
t^{d-1})^{n}$. More precisely, we have:

\begin{prop}
\begin{enumerate}
    \item $H(S/J,k) = \binom{k+n-1}{n-1}$, for $0\leq k\leq d-1$.
    \item $H(S/J,k) = \binom{k+n-1}{n-1} - n\binom{j+n-1}{n-1}$, for $d\leq k\leq \left\lfloor \frac{n(d-1)}{2} \right\rfloor$ and $j=k-d$.
    \item $H(S/J,k) = H(S/J,n(d-1)-k)$, for $k\geq \left\lceil \frac{n(d-1)}{2}\right\rceil$.
\end{enumerate}
\end{prop}

\begin{proof}
Use induction on $n$. Denote $H_{n}(t)=(1+t + \cdots +
t^{d-1})^{n}$. The case $n=1$ is trivial. The induction step follows
from the equality $H_{n}(t) = H_{n-1}(t)(1+t + \cdots + t^{d-1})$.
\end{proof}

\begin{cor}
\begin{enumerate}
  \item $|J_{k}|=0$, for $k\leq d-1$.
  \item $|J_{k}|=n\binom{j+n-1}{n-1}$, for $d\leq k\leq \left\lfloor \frac{n(d-1)}{2} \right\rfloor$ and $j=k-d$.
  \item $|J_{k}|= \binom{\left\lceil \frac{n(d-1)}{2} \right\rceil + j + n - 1}{n-1} -
        \binom{ \left\lfloor\frac{n(d-1)}{2} \right\rfloor -j + n - 1}{n-1} +
        n\binom{\left\lfloor\frac{n(d-1)}{2} \right\rfloor - d - j - n}{n-1}$, for
        $\left\lceil \frac{n(d-1)}{2} \right\rceil \leq k \leq (n-1)(d-1)-1$, where $j=k-\left\lceil \frac{n(d-1)}{2}
        \right\rceil$
  \item $|J_{k}| = \binom{(n-1)d+j}{n-1} - \binom{n-1+d-1-j}{n-1}$, for $(n-1)(d-1)\leq k \leq n(d-1)$, where
        $j=k-(n-1)(d-1)$.
\end{enumerate}
\end{cor}

\begin{proof}
Using $|J_k| = |S_k| - H(S/J,k)$ the proof follows from $2.4.1$.
\end{proof}

Suppose $f_{i} = \sum_{k=1}^{N}b_{ik}u_{k}$ for $1\leq i \leq n$
where $u_{1},u_{2},\ldots,u_{N}\in S$ are all the monomials of
degree $d$ decreasing ordered in revlex and $N=\binom{d+n-1}{n-1}$.
We denote $u_{k}=x^{\alpha_{k}}$. For example,
$\alpha_1=(d,0,\ldots,0)$, $\alpha_2=(d-1,1,0,\ldots,0)$ etc.

We take a generic transformation of coordinates $x_{i} \mapsto
\sum_{j=1}^{n}c_{ij}x_{j}$ for $i=1,\ldots,n$. Conca and Sidman
proved in \cite{CS} that we may assume that $c_{ij}$ are algebraically
independents over $K$. More precisely, if we consider the field
extension $K\subset L=K(c_{ij}|i,j=\overline{1,n})$ and if we set
\[ F_{i} = f_{i}(\sum_{j=1}^{n}c_{1j}x_{j},\ldots,\sum_{j=1}^{n}c_{nj}x_{j}) \in L[x_{1},\ldots,x_{n}],\;i=1,\ldots,n\]
then $J = Gin(I) = in(F_{1},\ldots,F_{n})\cap S$.

We write $F_{i} = \sum_{j=1}^{n}a_{ij}u_{j} + \cdots $ the monomial
decomposition of $F_{i}$ in $L[x_{1},\ldots,x_{n}]$. With these
notations, we have the following elementary lemma:

\begin{lema}
$J_{d}$ is revlex if and only if the following condition is
fulfilled:
\[ \Delta = \left|
\begin{array}[pos]{ccc}
    a_{11} & \cdots & a_{1n} \\
    \vdots &        & \vdots \\
    a_{n1} & \cdots & a_{nn}
\end{array}
\right|\neq 0.\]
\end{lema}

\begin{proof}
Suppose $\Delta\neq 0$. Since $|J_{d}|=n$, it is enough to show that
$u_{1},\ldots,u_{n}\in J$. Let $A=(a_{ij})_{\tiny
\begin{array}{c} i=1,n\\ j=1,n
\end{array}}$. Since $\Delta=det(A)\neq 0$, $A$ is invertible and we
have
\[A^{-1}
\left( \begin{array}[pos]{c}
    F_1 \\
    \vdots \\
    F_n
\end{array} \right) = \left( \begin{array}[pos]{c}
    H_1 \\
    \vdots \\
    H_n
\end{array} \right),
 \]
 where $H_i = u_i + $ small terms in revlex order. Therefore $LM(H_i)=u_i \in J$, for all $1\leq i\leq n$,
 where $LM(H_i)$ denotes the leading monomial of $H_i$ in the revlex order.

 Conversely, since $u_{1},\ldots,u_{n}\in J_d$, we can find some polynomials $H_i\in L[x_1,\ldots,x_n]$, with
 $LM(H_i)=u_i$, $1\leq i\leq n$, as linear combination of $F_i$'s. 
 If we denote $H_i=\sum_{j=1}^N \widetilde{a}_{ij}u_j$ and $\widetilde{A} = (\widetilde{a}_{ij})_{i,j=1,\ldots,n}$, it
 follows that there exists a map $\psi:L^{n} \rightarrow L^{n}$, given by a matrix $E=(e_{ij})_{i,j=1,\ldots,n}$, such
 that $\widetilde{A} = A\cdot E$. Now, since $det(\widetilde{A})\neq 0$ it follows that $\Delta=det(A)\neq 0$, as
 required.
\end{proof}

\begin{obs}{\em
By the changing of variables $\varphi$ given by $x_{i}\mapsto
\sum_{j=1}^{n}c_{ij}x_{j}$, $x^{\alpha_{k}}$ became
\[ m_{k}:=(\sum_{j=1}^{n}c_{1j}x_{j})^{\alpha_{k1}}\cdots (\sum_{j=1}^{n}c_{nj}x_{j})^{\alpha_{kn}} =
 (\sum_{|t|=\alpha_{k1}}c_{1}^{t}x^{t}) \cdots (\sum_{|t|=\alpha_{kn}}c_{n}^{t}x^{t}) ,\]
where, for any multiindex $t=(t_{1},\ldots,t_{n})$ we denoted
$x^t=x_1^{t_1}\cdots x_n^{t_n}$ and $c_{i}^{t} =
c_{i1}^{t_{1}}\cdots c_{in}^{t_{n}}$. Let $g_{kl}$ be the
coefficient in $c_{ij}$'s of $x^{\alpha_{l}}$ in the monomial
decomposition of $m_{k}$. Using the above writing of $m_k$, we claim
that: 
\[(1)\;\; g_{kl} = \tiny \sum_{\begin{array}{c} |t_{1}| = \alpha_{k1},\ldots,|t_{n}| = \alpha_{kn} \\ t_{1}+\cdots+t_{n} = \alpha_{l} \end{array}}
\left[ \binom{\alpha_{k1}}{t_{11}}\cdots \binom{\alpha_{kn}}{t_{n1}}
\right] \left[ \binom{\alpha_{k1}-t_{11}}{t_{12}}\cdots
\binom{\alpha_{kn}-t_{n1}}{t_{n2}}\right] \cdots \] \[ \left[
\binom{\alpha_{k1}-t_{11} - \cdots t_{1n-1}}{t_{1n}}\cdots
\binom{\alpha_{kn}-t_{n1} - \cdots - t_{nn-1}}{t_{nn}} \right] \cdot
c_{1}^{t_{1}}\cdots c_{n}^{t_{n}}. \] \normalsize Indeed, the
monomial $c_{1}^{t_{1}}\cdots c_{n}^{t_{n}}$ appear in the
coefficient of $x^{\alpha_{l}}$ in the expansion of $m_k$ if and
only if $t_1+\cdots+t_n = \alpha_l$ and $|t_{1}| =
\alpha_{k1},\ldots,|t_{n}| = \alpha_{kn}$. Moreover, by Newton
binomial, the coefficient of $x_1^{t_{i1}}\cdots x_n^{t_{in}}$ in
$(\sum_{j=1}^{n}c_{ij}x_{j})^{\alpha_{k1}}$ is
$\binom{\alpha_{k1}}{t_{i1}}
\binom{\alpha_{k1}-t_{i1}}{t_{i2}}\cdots  \binom{\alpha_{k1}-t_{i1}
- \cdots t_{i,n-1}}{t_{in}}c_i^{t_i}$ for any $1\leq i\leq n$, and
thus we proved the claim.

Since $a_{il} = \sum_{k=1}^{N}b_{ik}\cdot g_{kl}$, from the
Cauchy-Binet formula we get:
\[ \Delta =   \sum_{1\leq k_{1} < k_{2} <\ldots < k_{n} \leq N} B_{k_{1},k_{2},\ldots,k_{n}}G_{k_{1},k_{2},\ldots,k_{n}}, \; where \]
\[ B_{k_{1},k_{2},\ldots,k_{n}} =
\left| \begin{array}[pos]{ccc}
    b_{1k_{1}} & \cdots & b_{1k_{n}} \\
    \vdots &  & \vdots \\
    b_{nk_{1}} & \cdots & b_{nk_{n}}
\end{array} \right| \; and \;\;
G_{k_{1},k_{2},\ldots,k_{n}} = \left| \begin{array}[pos]{ccc}
    g_{k_{1}1} & \cdots & g_{k_{n}1} \\
    \vdots &  & \vdots \\
    g_{k_{1}n} & \cdots & g_{k_{n}n}
\end{array} \right|. \]
}\end{obs}

Now, we are able to prove the main result of our paper.

\begin{teor}
If $f_{i}\in K[x_{i},\ldots,x_{n}]$ then $J_{d}$ is revlex. In
particular, if $S/I$ is a monomial complete intersection, then
$J_{d}$ is revlex.
\end{teor}

\begin{proof}
Let $k_{i}=\binom{i+d-1}{d}$, for any $i=1,\ldots,n$. Then
$u_{k_{i}} = x_{i}^{d}$. Recall our notation, $u_k=x^{\alpha_k}$. We
have $b_{11}\neq 0$, otherwise $I = (f_{1},\ldots,f_{n}) \subset
(x_{2},\ldots,x_{n})$ contradicting the fact that $I$ is an Artinian
ideal.  Using a similar argument, we get $b_{ik_{i}}\neq 0$ for all
$1 \leq i \leq n$. Thus, multiplying each $f_{i}$ with
$b_{ik_{i}}^{-1}$, we may assume $b_{ik_{i}}=1$ for all $1 \leq i
\leq n$. In other words, $f_{i}=x_{i}^{d} + f'_{i}$, where $f'_{i}$
contains monomials smaller than $x_{i}^{d}$ in the revlex order.
Also, since $f_{i}\in K[x_{i},\ldots,x_{n}]$ we have $b_{i'k_{i}}=0$
for any $i'>i$. In particular,$B_{k_{1},\ldots,k_{n}}=1$.

In the expansion of the determinant $G_{k_{1},\ldots,k_{n}}$,
appears the term $g_{k_{1}1}\cdot g_{k_{2}2} \cdots g_{k_{n}n} =
r\cdot (c_{11}^{d})(c_{21}^{d-1}c_{22})\cdots (c_{i}^{\alpha_{i}})
\cdots (c_{n}^{\alpha_{n}})$, where $r$ is a nonzero (positive)
integer. Indeed, by $(1)$, we have $g_{11}=c_{11}^d$, $g_{k_{2}2}=d
c_{21}^{d-1}c_{22}$ and, in general, $g_{k_{i}i} = $ some binomial
coefficient $\cdot c_i^{\alpha_i}$. We claim that
$m=(c_{11}^{d})(c_{21}^{d-1}c_{22})\cdots (c_{i}^{\alpha_{i}})
\cdots (c_{n}^{\alpha_{n}})$ doesn't appear again in the expansion
of $\Delta$.

Since $f_{i}\in k[x_{i},\ldots,x_{n}]$, in the monomials in
$(c_{tl})$ of $a_{ij}$ there are no $c_{tl}$'s with $t<i$. Also, all
the monomials of $f'_{i}$ contain variables $x_{t}$ with $t>i$.
Corresponding to them, in $a_{ij}$'s there are $c_{tj}$'s with
$t>i$. Thus in $a_{il}$ the only monomials in $c_{i1},\ldots,c_{in}$
of degree $d$ comes from
$\varphi(x_i^{d})=(\sum_{j=1}^{n}c_{ij}x_j)^d$, the other monomials
being multiples of some $c_{tl}$ with $t>i$. Consequently, in the
expansion of $\Delta$, the monomials of the type
$c_1^{\beta_1}\cdots c_n^{\beta_n}$, where $\beta_1,\ldots,\beta_n$
are multiindices with $|\beta_1|= \cdots = |\beta_n|=d$ comes only
from $\varphi(x_1^d),\ldots,\varphi(x_n^d)$.

On the other hand, for any $1\leq i\leq n$, $c_i^{\alpha_i}$ is unique between the
monomials in $c_{tl}$'s from $\varphi(x_n^d)$, because they are of the type $c_i^{\gamma}$, 
where $\gamma$ is a multiindex with $|\gamma|=d$. 
From these facts, it follows that the
monomial $m$ is unique in the monomial expansion of $\Delta$ and
occurs there with a nonzero coefficient. Thus $\Delta\neq 0$ and by
applying Lemma $2.4.3$ we complete the proof of the theorem.
\end{proof}

\begin{obs}{\em
In the case $n=2$ and $n=3$, $J_{d}$ is revlex for any
$(n,d)$-complete intersection. Indeed, in the case $n=2$, $J$ itself
is revlex since it is strongly stable. In the case $n=3$, since
$|J_{d}|=3$ and $J$ is strongly stable, it follows that either (a)
$J_{d}=(x_{1}^{d},x_{1}^{d-1}x_{2},x_{1}^{d-2}x_{2}^{2})$, either
(b) $J_{d}=(x_{1}^{d},x_{1}^{d-1}x_{2},x_{1}^{d-1}x_{3})$. But in
the case $(b)$, the map $(S/J)_{d-1} \stackrel{\cdot
x_{3}}{\longrightarrow} (S/J)_{d}$ is not injective, because
$x_{1}^{d-1}\neq 0$ in $(S/J)_{d-1}$ and $x_{1}^{d-1}x_{3}= 0$ in
$(S/J)_{d}$. This is a contradiction with the fact that $x_3$ is a weak
Lefschetz element on $S/J$ and therefore, $J_d$ is
revlex.}\end{obs}

\begin{lema}
(a) $a_{i1}=f_{i}(c_{11},\ldots,c_{n1})$ for all $1\leq i\leq n$.

(b) If $1\leq l\leq n$ is an integer then the sequence
$a_{1l},a_{2l},\ldots,a_{nl}$ is regular as a sequence of
          polynomials in $K[c_{ij}|\;1\leq i,j\leq n]$.
\end{lema}

\begin{proof}
Substituting $x_j=0$ for $j\neq 1$ in $F_i$ we get (a). In order to
prove (b), firstly notice that $a_{11},a_{21},\ldots,a_{n1}$ is a
regular sequence on $K[c_{11},c_{21},\ldots,c_{n1}]$, since $f_1,\ldots,f_n$ is a regular sequence on $K[x_1,\ldots,x_n]$ and $c_{11},c_{21},\ldots,c_{n1}$ are algebraically independent over $K$. 

Let $1\leq l\leq n$ be an integer. We claim that
\[(*)\; \frac{K[c_{ij}|\;1\leq i,j\leq n]}{(a_{1l},\ldots,a_{nl},c_{i1}-c_{ij}\;,\;1\leq i\leq n,\;2\leq j\leq n)} \cong
        \frac{K[c_{11},c_{21},\ldots,c_{n1}]}{(a_{11},a_{21},\ldots,a_{n1})}.\]
Indeed, by $(1)$, if we put $c_{ij}=c_{i1}$ for all $1\leq i\leq n$ and $2\leq j\leq n$ in the expansion of $g_{kl}$
we obtain $r_l \cdot g_{k1}$, where $r_l$ is a strictly positive integer, which depends only on $l$, and therefore, 
$a_{il}$ became $r_l\cdot a_{i1}$. From (*) it follows that $a_{1l},\ldots,a_{nl},c_{i1}-c_{ij}\;for\;1\leq i\leq n,\;2\leq j\leq n$ is a system of parameters for $K[c_{ij}|\;1\leq i,j\leq n]$ and thus $a_{1l},\ldots,a_{nl}$ is a regular sequence on $K[c_{ij}|\;1\leq i,j\leq n]$, so we proved (b).
\end{proof}

As we noticed in Remark $2.4.6$, for $n=3$, the conclusion of Theorem
$2.4.5$ holds for any regular sequence $f_1,f_2,f_3$ of homogeneous
polynomials of degree $d$. In the following, we give another proof
of this, without using the fact that $S/(f_1,f_2,f_3)$ has the
(WLP), i.e. $x_3$ is a weak Lefschetz element for $S/J$. Also, we
get the same conclusion for the case $n=4$ and $d=2$. However, this
approach do not works in the general case.

\begin{prop}
(a) If $f_1,f_2,f_3\in K[x_1,x_2,x_3]$ is a regular sequence of
homogeneous polynomials of degree $d\geq 2$, $I=(f_1,f_2,f_3)$
    and $J=Gin(I)$, the generic initial ideal of $I$, with respect to the reverse lexicographical order, then
    $J_d$ is a revlex set.

(b) If $f_1,f_2,f_3,f_4\in K[x_1,x_2,x_3,x_4]$ is a regular sequence
of homogeneous polynomials of degree $2$, $I=(f_1,f_2,f_3,f_4)$
    and $J=Gin(I)$, the generic initial ideal of $I$, with respect to the reverse lexicographical order, then
    $J_2$ is a revlex set.
\end{prop}

\begin{proof}
(a) Let $A=(a_{ij})_{i,j=\overline{1,3}}$. Since $Gin(f_1,f_2)$ is
strongly stable, it follows by Lemma $2.4.3$ that
    $\Delta_{3}=\left| \begin{array}[pos]{cc} a_{11} & a_{12}\\ a_{21} & a_{22} \end{array} \right| \neq 0$.
    Analogously, $\Delta_{2}=\left| \begin{array}[pos]{cc} a_{11} & a_{12}\\ a_{31} & a_{32} \end{array} \right| \neq 0$ and $\Delta_{1}=\left| \begin{array}[pos]{cc} a_{21} & a_{22}\\ a_{31} & a_{32} \end{array} \right| \neq 0$.
    We have $\Delta = a_{13}\Delta_{1} - a_{23}\Delta_{2} + a_{33}\Delta_{3}$.
    Suppose $\Delta = 0$. It follows $a_{13}\Delta_{1} = a_{23}\Delta_{2} - a_{33}\Delta_{3}$
    and therefore, since $a_{13},a_{23},a_{33}$ is a regular sequence in $K[c_{ij}| i,j=\overline{1,3}]$,
    we get $\Delta_{1} \in (a_{23},a_{33})$. The first three monomials of degree $d$ in revlex order are $x_1^d$,
    $x_1^{d-1}x_2$ and $x_1^{d-2}x_2^2$. It follows that the degree of
    $a_{i1}$, $a_{i2}$ and $a_{i3}$ in $c_{21},c_{22},c_{23}$ is $0$, $1$, respectively $2$,
    for any $1\leq i\leq 3$. Therefore, the degree of $\Delta_{1}$ in
    the variables $c_{21},c_{22},c_{23}$ is $1$, but the degree of $a_{23}$ and $a_{33}$ in $c_{21},c_{22},c_{23}$ is
    $2$, which is impossible, since $\Delta_{1} \in (a_{23},a_{33})$.

(b) Let $A=(a_{ij})_{i,j=\overline{1,4}}$. Since any three
polynomials from $f_1,f_2,f_3,f_4$ form a regular sequence, it
follows from (a) that any $3\times 3$
    minor of the matrix $\widetilde{A}=(a_{ij})_{\tiny \begin{array}{c} i=\overline{1,4} \\ j=\overline{1,3}
    \end{array} \normalsize}$ is
    nonzero. Let $\Delta_i$ be the minor obtained from $\widetilde{A}$ by erasing the $i$-row. Suppose $\Delta=0$. It
    follows that $a_{14}\Delta_1 = a_{24}\Delta_2 - a_{34}\Delta_3 + a_{44}\Delta_4$ and therefore, since
    $a_{14},a_{24},a_{34},a_{44}$ is a regular sequence in $K[c_{ij}| i,j=\overline{1,4}]$, we get
    $\Delta_1 \in (a_{24},a_{34},a_{44})$. Since the first $4$ monomials in
    revlex are $x_1^2,x_1x_2,x_2^2,x_1x_3$, we get a contradiction from the fact that
    the degree of $\Delta_1$ in the variables $c_{31},c_{32},c_{33},c_{34}$ is zero, but the degree of
    $a_{24},a_{34},a_{44}$ in $c_{31},c_{32},c_{33},c_{34}$ is $1$.
\end{proof}


\begin{obs}{\em
The hypothesis that $K$ is a field with $char(K)=0$ is essential.
Indeed, suppose $char(K)=p$ and $I = (x_1^{p},x_2^{p})\subset
K[x_1,x_2]$. Then, simply using the definition of the generic
initial ideal, we get $Gin(I)=I$ and, obviously,
$I_{p}=\{x_1^{p},x_2^{p}\}$ is not revlex.

Also, the hypothesis that $f_1,\ldots,f_n$ is a regular sequence of
homogeneous polynomials is essential. Let $I=(f_1,f_2,f_3)\subset
K[x_1,x_2,x_3]$, where $f_1=x_1^2$, $f_2=x_1x_2$ and $f_3=x_1x_3$.
In order to compute the generic initial ideal of $I$ we can take a
generic transformation of coordinates with an upper triangular
matrix, i.e. $x_1 \mapsto x_1,\; x_2 \mapsto x_2+c_{12}x_1,\;
x_3\mapsto x_3 + c_{23}x_2 + c_{13}x_1$, where $c_{ij}\in K$ for all
$i,j$ (see \cite[\S 15.9]{E}). We get
\[ F_{1}(x_1,x_2,x_3):=f_{1}(x_1,x_2+c_{12}x_1,x_3 + c_{23}x_2 + c_{13}x_1)= x_1^{2}, \]
\[ F_{2}(x_1,x_2,x_3):=f_{2}(x_1,x_2+c_{12}x_1,x_3 + c_{23}x_2 + c_{13}x_1)= c_{12}x_1^{2}+x_1x_2, \]
\[ F_{3}(x_1,x_2,x_3):=f_{3}(x_1,x_2+c_{12}x_1,x_3 + c_{23}x_2 + c_{13}x_1)= c_{13}x_1^{2}+c_{23}x_1x_2 + x_1x_3. \]
The generic initial ideal of $I$, $J=in(F_1,F_2,F_3)$ satisfies $J_2=I_{2}$, but $I_2$ is not revlex.}
\end{obs}

\section{Several examples of computation of the Gin.}

Let $I=(f_{1},\ldots,f_{n})\subset S = K[x_1,\ldots,x_n]$ be an
ideal generated by a regular sequence $f_{1},\ldots,f_{n}\in S$ of
homogeneous polynomials of degree $d$. Let $J=Gin(I)$ be the generic
initial ideal of $I$, with respect to the revlex order.

In section $2.3$, the case $n=3$ and $d\geq 2$ was treated completely,
when $S/(f_1,f_2,f_3)$ has (SLP), see Proposition $2.3.3$.
In the following, we discuss some particular cases with $n\geq 4$.

\paragraph{The case $n=4$, $d=2$.}

We assume that $S/I$ has (SLP). From Wiebe's Theorem, it follows
that $x_{4}$ is a strong Lefschetz element for $S/J$. For a positive
integer $k$, we denote $Shad(J_{k}) = \{x_1,\ldots,x_n\}J_{k}$. We
have $H(S/J,t) = (1+t)^{4} = 1+4t+6t^{2}+4t^{3}+t^{4}$.

We have $|J_{2}|=4$. From Proposition $2.4.8$, $J_2$ is revlex,
therefore
\[ J_{2}=\{x_{1}^{2},x_{1}x_{2},x_{2}^{2},x_{1}x_{3}\} = \{\{x_{1},x_{2}\}^{2},x_{1}x_{3}\} .\]
We have $|Shad(J_{2})| = 12$. On the other hand, $|J_{3}|=16$, so we
need to add $4$ new generators at $Shad(J_{2})$ to get $J_{3}$. If
we add a new monomial which is divisible by $x_{4}^{2}$, then the
map $(S/J)_{1} \stackrel{\cdot x_{4}^{2}}{\longrightarrow}
(S/J)_{3}$, will be no longer injective. Since $|(S/J)_{1}| =
|(S/J)_{3}|$, we get a contradiction with the fact that $x_{4}$ is a
strong Lefschetz element for $S/J$. But there exists only $16$
monomials in $S$ which are not multiple of $x_4^2$. Thus
\[J_{3} = \{\{x_{1},x_{2},x_{3}\}^{3}, x_{4}\{x_{1},x_{2},x_{3}\}^{2}\},\;and\; therefore\]
\[ Shad(J_{3})=\{\{x_{1},x_{2},x_{3}\}^{4}, x_{4}\{x_{1},x_{2},x_{3}\}^{3}, x_{4}^{2}\{x_{1},x_{2},x_{3}\}^{2}\}. \]
Since $|Shad(J_{3})| = 31$ and $|J_{4}| = |S_{4}|-|(S/J)_{4}| =35 -
1 = 34$ we have to add $3$ new generators at $Shad(J_{3})$ in order
to get $J_{4}$. Since $J$ is strongly stable, these new generators
are $x_{4}^{3}x_{1}$, $x_{4}^{3}x_{2}$ and $x_{4}^{3}x_{3}$. So
\[ J_{4}=\{x_{1},x_{2},x_{3},x_{4}\}^{4}\setminus \{x_{4}^{4}\}.\;We\;get\;\;
   Shad(J_{4})= \{x_{1},x_{2},x_{3},x_{4}\}^{5}\setminus \{x_{4}^{5}\}\]
   and since $J_{5}=S_{5}$ it follows that we must add $x_{4}^{5}$ at $Shad(J_{4})$ to obtain $J_{5}$. From now one,
   we cannot add any new monomial. $J$ is
   the ideal generated by all monomials added at some step $k$ to $Shad(J_{k})$, thus
   we proved the following proposition:

\begin{prop}
If $I=(f_{1},f_{2},f_{3},f_{4})$ is an ideal generated by a regular
sequence of homogeneous polynomials $f_{1},f_{2},f_{3},f_{4} \in
S=k[x_{1},x_{2},x_{3},x_{4}]$ of degree $2$ such that the algebra
$S/I$ has (SLP) then the generic initial ideal of $I$ with respect
to the revlex order is
\[J = (x_{1}^{2},\; x_{1}x_{2},\; x_{2}^{2},\; x_1x_3,\;
       x_{2}x_{3}^{3},\; x_{3}^{3}, \; x_{3}^{2}x_{4}, \; x_{3}^{2}x_{4},\;
       x_{4}^{3}x_{1},\; x_{4}^{3}x_{2},\; x_{4}^{3}x_{3},\; x_{4}^{5}). \]
In particular, this assertion holds for a generic sequence of
homogeneous polynomials $f_{1},f_{2},f_{3},f_{4}\in S$ or if
$f_{i}\in k[x_{i},\ldots,x_{4}]$, $1\leq i\leq 4$.
\end{prop}

\paragraph{The case $n=5$, $d=2$.}

In the following, we suppose that $S/I$ has (SLP), so $x_{5}$ is a
strong Lefschetz element for $S/J$. Also, we suppose that $J_{2}$ is
revlex. We have $H(S/J,t) = (1+t)^{5} = 1+5t+10t^{2}+10t^{3}+5t^{4}
+ t^{5}$. We have $|J_{2}|=5$. Since $J_{2}$ is revlex from the
assumption, we have $J_{2} =
\{\{x_1,x_2\}^{2},x_{3}\{x_{1},x_{2}\}\}$. So \[ Shad(J_{2}) = \{
\{x_{1},x_{2}\}^{3}, \{x_{1},x_{2}\}^{2}\{x_3,x_4,x_5\},
x_3\{x_1,x_2\}\{x_3,x_4,x_5\}\}.\] We have $|Shad(J_{2})|=19$. On
the other hand $|J_{3}| = |S_{3}| - |(S/J)_{3}| = 35 - 10 = 25$, so
we must add $6$ new generators, from a list of $16$ monomials, at
$Shad(J_{2})$ to get $J_{3}$.

Since $x_{5}$ is a strong Lefschetz element for $S/J$ it follows
that we cannot add any monomial of the form $x_{5}\cdot m$, where
$m$ is nonzero in $(S/J)_{2}$ because, in that case, the map
$(S/J)_{2} \stackrel{\cdot x_{5}}{\rightarrow} (S/J)_{3}$ will be no
longer injective. But there are $|(S/J)_{2}| = 10$ such monomials
$m$. Therefore, we must add the remaining $6$ monomials,
$x_{3}^{3},x_{3}^{2}x_{4}, x_{1}x_{4}^{2},x_{2}x_{4}^{2},
x_{3}x_{4}^{2}, x_{4}^{3}$. Thus
\[ J_{3} = \{\{x_{1},x_{2},x_{3},x_{4}\}^{3} , x_{5}(\{x_{1},x_{2},x_{3}\}^{2}\setminus\{x_{3}^{2}\})\}.\;Therefore:\]
\[ Shad(J_{3}) =\{ \{x_{1},x_{2},x_{3},x_{4}\}^{4} , x_{5}\{x_{1},x_{2},x_{3},x_{4}\}^{3} ,
x_{5}^{2}( \{x_{1},x_{2},x_{3}\}^{2}\setminus\{x_{3}^{2}\})\}. \] We
have $|Shad(J_{3})|=60$ and $|J_{4}| = |S_{4}| - |(S/J)_{4}| = 70 -
5 = 65$. So we need to add $5$ new generators at $Shad(J_{3})$ to
get $J_{4}$. If we add a monomial which is divisible by $x_{5}^{3}$
we obtain a contradiction from the fact that the map $(S/J)_{1}
\stackrel{\cdot x_{5}^{3}}{\rightarrow} (S/J)_{4}$ is no longer
injective. Therefore, we must add: $x_{3}^{2}x_{5}^{2},
x_{1}x_{4}x_{5}^{2}, x_{2}x_{4}x_{5}^{2}, x_{3}x_{4}x_{5}^{2},
x_{4}^{2}x_{5}^{2}$, and so
\[ J_{4} = \{ \{x_{1},x_{2},x_{3},x_{4}\}^{4} , x_{5}\{x_{1},x_{2},x_{3},x_{4}\}^{3} ,
x_{5}^{2} \{x_{1},x_{2},x_{3},x_{4}\}^{2} \}. \]\[ So
\;\;Shad(J_{4}) = \{ \{x_{1},x_{2},x_{3},x_{4}\}^{5} , \cdots ,
x_{5}^{3} \{x_{1},x_{2},x_{3},x_{4}\}^{2} \}.\] 

We have $|J_{5}| - |Shad(J_{4})| = 4$, so we must add $4$ new generators at $Shad(J_{4})$ to
get $J_{5}$. Since $J$ is strongly stable, these new generators are:
$x_{5}^{4}x_{1},x_{5}^{4}x_{2},x_{5}^{4}x_{3},x_{5}^{4}x_{4}$. Therefore $J_{5} = \{ \{x_{1},\ldots,x_{5}\}^{5}\setminus\{x_{5}^{5}\}\}$. Finally, we must add $x_5^6$ to $Shad(J_5)$ in
order to obtain $J_6$. We proved the following proposition, with the help of 
\cite[Theorem 1.2]{CS} and Theorem $2.4.5$.

\begin{prop}
If $I=(f_{1},f_{2},\ldots,f_{5})\subset K[x_1,\ldots,x_5]$ is an
ideal generated by a generic (regular) sequence of homogeneous
polynomials of degree $2$ or if $f_{1},f_{2},\ldots,f_{5}$ is a
regular sequence of homogeneous polynomials of degree $2$ with
$f_{i}\in K[x_{i},\ldots,x_{5}]$ for $i=1,\ldots,5$ then $J=Gin(I)$
the generic initial ideal of $I$ with respect to the revlex order
is:
\[ J = (x_{1}^{2},\; x_{1}x_{2},\; x_{2}^{2},\;x_{1}x_{3},\; x_{2}x_{3}, \;
       x_{3}^{3},\; x_{3}^{2}x_{4},\;  x_{1}x_{4}^{2},\; x_{2}x_{4}^{2},\;  x_{3}x_{4}^{2},\;  x_{4}^{3},\;\]\[
       x_{3}^{2}x_{5}^{2},\;x_{1}x_{4}x_{5}^{2},\; x_{2}x_{4}x_{5}^{2},\;x_{3}x_{4}x_{5}^{2},\;  x_{4}^{2}x_{5}^{2},\;
       x_{5}^{4}x_{1},\; x_{5}^{4}x_{2},\; x_{5}^{4}x_{3},\; x_{5}^{4}x_{4},\;  x_{5}^{6} ) \]
\end{prop}

\paragraph{The case $n=4$, $d=3$.}

We suppose that $S/I$ has (SLP), so $x_{4}$ is a strong Lefschetz
element for $S/J$. Also, we suppose that $J_{3}$ is revlex.

We have $H(S/J,t) = (1+t+t^{2})^{4} = (1+2t+3t^{2}+2t{3}+t^{4})^{2}
= $
\[ = 1 +4t + 10t^{2} +16t^{3} + 19t^{4} + 16t^{5} + 10t^{6} + 4t^{7} +t^{8}.\]
Since $|J_{3}|=4$ and $J_{3}$ is revlex, it follows that $J_{3} =  \{x_{1},x_{2}\}^{3}$.
Therefore, we have \linebreak $Shad(J_{3}) = \{ \{x_{1},x_{2}\}^{4}, \{x_{1},x_{2}\}^{3}\{x_{3},x_{4}\}\}$. 
Since $|J_{4}| - |Shad(J_{3})|= 4$, we must add $4$ new generators to $Shad(J_{3})$ to
obtain $J_{4}$. Since $x_{4}$ is a strong Lefschetz element for
$S/J$ we cannot add any monomial of the form $x_{4}\cdot m$, where
$m\neq 0$ in $J_{3}$. Therefore, since $J$ is strongly stable, we
have to choose $3$ monomials from the list
$x_{3}^{2}\{x_{1},x_{2}\}^{2},x_{3}^{3}\{x_{1},x_{2}\},x_{3}^{4}$. There
are two different chooses: either we add (I) $x_{3}^{2}\{x_{1},x_{2}\}^{2}$, either (II) $x_{3}^{2}x_{1}\{x_{1},x_{2},x_{3}\}$.

In the case (I), we get $J_{4}=\{ \{x_{1},x_{2}\}^{4},
\{x_{1},x_{2}\}^{3}\{x_{3},x_{4}\}, x_{3}^{2}\{x_{1},x_{2}\}^{2}\}$, so
\[Shad(J_{4}) = \{ \{x_{1},x_{2}\}^{5}, \{x_{1},x_{2}\}^{4}\{x_{3},x_{4}\}, \{x_{1},x_{2}\}^{3}\{x_{3},x_{4}\}^{2}, x_{3}^{2}\{x_{3},x_{4}\}\{x_{1},x_{2}\}^{2}\}. \]
Since $|J_{5}|-|Shad(J_{4})|=40-34=6$, we need to add $6$ new generators
at $Shad(J_{4})$ to get $J_{5}$. Since $x_{4}$ is a strong Lefschetz
element for $S/J$ we cannot add any monomial of the form
$x_{4}^{2}m$, where $m$ is a nonzero monomial in $J_{3}$. So, we
must add: $x_{3}^{4}\{x_{1},x_{2},x_{3}\},
x_{4}x_{3}^{3}\{x_{1},x_{2},x_{3}\}$. Thus
$J_{5} = \{ \{ x_{1},x_{2},x_{3}\}^{5}, x_{4}\{ x_{1},x_{2},x_{3}\}^{4}, x_{4}^{2}\{x_{1},x_{2}\}^{3}\}$. 

In the case (II), we have $J_{4}=\{ \{x_{1},x_{2}\}^{4},
\{x_{1},x_{2}\}^{3}\{x_{3},x_{4}\},
x_{1}x_{3}^{2}\{x_{1},x_{2},x_{3}\} \}$, so $Shad(J_{4})$ is the set
$\{ \{x_{1},x_{2}\}^{5}, \{x_{1},x_{2}\}^{4}\{x_{3},x_{4}\}, \{x_{1},x_{2}\}^{3}\{x_{3},x_{4}\}^{2}, x_{3}^{2}x_{1}\{x_{3},x_{4}\}\{x_{1},x_{2}\}, x_{3}^{3}x_{1}\{x_{3},x_{4}\}\}$.
Since $|J_{5}|-|Shad(J_{4})|=40-34=6$, we must add $6$ new generators at $Shad(J_{4})$ to get $J_{5}$. 
Since $x_{4}$ is a strong-Lefschetz element for $S/J$, we cannot add any monomial of
the form $x_{4}^{2}m$, where $m\neq 0$ in $J_{3}$. So, we must add:
$x_{3}^{3}x_{2}^{2}, x_{3}^{4}x_{2}, x_{3}^{5},
x_{4}x_{3}^{2}x_{2}^{2}, x_{4}x_{3}^{3}x_{2}, x_{4}x_{3}^{4}$. Thus
\[ J_{5} = \{ \{ x_{1},x_{2},x_{3}\}^{5}, x_{4}\{ x_{1},x_{2},x_{3}\}^{4}, x_{4}^{2}\{x_{1},x_{2}\}^{3}\}, \]
the same as in the case (I). Thus, in both cases (I) and (II), we get:
\[Shad(J_{5}) = \{ \{ x_{1},x_{2},x_{3}\}^{6}, x_{4}\{ x_{1},x_{2},x_{3}\}^{5}, x_{4}^{2}\{x_{1},x_{2},x_{3}\}^{4} , x_{4}^{3}\{x_{1},x_{2}\}^{3}\}. \]
Since $|Shad(J_{5})| = |S_{6}|-16$ and $|J_{6}| = |S_{6}| - 10$, we must add $6$ new generators to
$Shad(J_{5})$ in order to obtain $J_6$. Since $x_{4}$ is a strong-Lefschetz element for $S/J$, these new
generators are not divisible by $x_{4}^{4}$. So, we add
$x_{4}^{3}x_{3}\{x_{1},x_{2}\}^{2},x_{4}^{3}x_{3}^{2}\{x_{1},x_{2}\},
x_{4}^{3}x_{3}^{3}$ and thus,
\[ J_{6} = \{ \{x_{1},x_{2},x_{3}\}^{6}, x_{4}\{ x_{1},x_{2},x_{3}\}^{5}, x_{4}^{2}\{x_{1},x_{2},x_{3}\}^{4},x_{4}^{3}\{x_{1},x_{2},x_{3}\}^{3}\}. \;So\]
\[ Shad(J_{6}) = \{ \{x_{1},x_{2},x_{3}\}^{7}, x_{4}\{ x_{1},x_{2},x_{3}\}^{6}, \ldots, x_{4}^{4}\{x_{1},x_{2},x_{3}\}^{3}\}. \]
$|S_{7}|-|Shad(J_{6})| = 6+4 = 10$ and $|S_{7}| - |J_{7}|=4$, so we
must add $6$ new generators at $Shad(J_{6})$ to get $J_{7}$. Using
the same argument, these new generators must be
$x_{4}^{5}\{x_{1},x_{2},x_{3}\}^{2}$ and therefore
$J_{7} = \{ \{x_{1},x_{2},x_{3}\}^{7}, x_{4}\{ x_{1},x_{2},x_{3}\}^{6}, \ldots, x_{4}^{5}\{x_{1},x_{2},x_{3}\}^{2} \}$. We get
\[ Shad(J_{7}) = \{ \{x_{1},x_{2},x_{3}\}^{8}, x_{4}\{ x_{1},x_{2},x_{3}\}^{7}, \ldots, x_{4}^{6}\{x_{1},x_{2},x_{3}\}^{2}\}. \]
Since $|S_{8}| - |Shad(J_{7})| = 4$ and $|S_{8}|-|J_{8}| = 1$, we
must add $3$ new generators at $Shad(J_{7})$ in order to get $J_8$.
Since $x_{4}$ is strong-Lefschetz, these new generators are $x_{4}^{7}\{x_{1},x_{2},x_{3}\}$, so
$J_8 = \{x_{1},x_{2},x_{3},x_{4}\}^{8}\setminus\{x_{4}^{8}\}$. Finally, we must add $x_{4}^{9}$ to
$Shad(J_8)$ in order to obtain $J_9$. We proved the following proposition, with the help of 
\cite[Theorem 1.2]{CS} and Theorem $2.4.5$.

\begin{prop}
If $I=(f_{1},f_{2},f_{3},f_{4})\subset K[x_1,x_2,x_3,x_4]$ is an
ideal generated by a generic (regular) sequence of homogeneous
polynomials of degree $3$ or if $f_{1},f_{2},f_{3},f_{4}$ is a
regular sequence of homogeneous polynomials of degree $3$ with
$f_{i}\in k[x_{i},\ldots,x_{4}]$, for $i=1,\ldots,4$, then
$J=Gin(I)$ the generic initial ideal of $I$ with respect to the
revlex order has one of the following forms:
\[(I)\;\;\;\;\;\;\;\;\;\;\;\;\;\;\;\;\;\;\;\; J = (\{x_{1},x_{2}\}^{3},\;x_{3}^{2}\{x_{1},x_{2}\}^{2},\; x_{3}^{4}\{x_{1},x_{2},x_{3}\},\;
                  x_{4}x_{3}^{3}\{x_{1},x_{2},x_{3}\}, \]
          \[x_{4}^{3}x_{3}\{x_{1},x_{2}\}^{2},x_{4}^{3}x_{3}^{2}\{x_{1},x_{2}\}, x_{4}^{3}x_{3}^{3},\;
            x_{4}^{5}\{x_{1},x_{2},x_{3}\}^{2},\; x_{4}^{7}\{x_{1},x_{2},x_{3}\},\; x_{4}^{9} ) \]
\[(II)\; J = (\{x_{1},x_{2}\}^{3},\;x_{3}^{2}x_{1}\{x_{1},x_{2},x_{3}\},\;x_{3}^{3}x_{2}^{2},\; x_{3}^{4}x_{2},\;
                   x_{3}^{5},\; x_{4}x_{3}^{2}x_{2}^{2},\; x_{4}x_{3}^{3}x_{2},\; x_{4}x_{3}^{4},  \]
          \[x_{4}^{3}x_{3}\{x_{1},x_{2}\}^{2},x_{4}^{3}x_{3}^{2}\{x_{1},x_{2}\}, x_{4}^{3}x_{3}^{3},\;
            x_{4}^{5}\{x_{1},x_{2},x_{3}\}^{2},\;x_{4}^{7}\{x_{1},x_{2},x_{3}\}, \; x_{4}^{9}) \]
\end{prop}

\begin{obs}
It seems Conca-Herzog-Hibi noticed in \cite{CHH}, page $838$, that, if $f_{1},f_{2},f_{3},f_{4}$
is a generic sequence of homogeneous polynomials of degree $3$ then the generic initial ideal $J$
has the form (I), and $J=Gin(x_1^3,x_2^3,x_3^3,x_4^3)$ has the form (II).
\end{obs}

\vspace{2mm} \noindent {\footnotesize
\begin{minipage}[b]{10cm}
 Mircea Cimpoeas, Junior Researcher\\
 Institute of Mathematics of the Romanian Academy\\
 Bucharest, Romania\\
 E-mail: mircea.cimpoeas@imar.ro
\end{minipage}}

\begin{thebibliography}{99}
  \bibitem{saf}Sarfraz Ahmad, Imran Anwar, "An upper bound for the regularity of ideals of Borel type", Preprint, 2006.
  \bibitem{A}  D.\ Anick "Thin algebras of embedding dimension three", J.Algebra 100 (1986), 235-259.
  \bibitem{ah} Annetta Aramova, J\"urgen Herzog "p-Borel principal ideals", Illinois J.Math.41,no 1.(1997),103-121.
  \bibitem{BH} W.\ Bruns, J.\ Herzog "Cohen-Macaulay rings", Second Edition, Cambridge, 1998.
  \bibitem{BS} D.\ Bayer, M. Stillman "A criterion for detecting m-regularity", Invent. Math 87 (1987) 1-11.
  \bibitem{CS} G.\ Caviglia, E. Sbarra, "Characteristic-free bounds for the Castelnuovo Mumford regularity", Compos.
               Math. 141(2005), no.6, 1365-1373.
  \bibitem{mir}Mircea Cimpoea\c s "A generalisation of Pardue's formula", Bull. Math. Soc. Sci. Math. Roumanie (N.S.)
               49(97), no. 4, 2006, 315-334.
  \bibitem{mir2}Mircea Cimpoea\c s "Monomial ideals with linear upper bound regularity", Preprint 2006.
  \bibitem{mir3}Mircea Cimpoea\c s "Generic initial ideal for complete intersections of embedding dimension three with
               strong Lefschetz property", Bull. Math. Soc. Sci. Math. Roumanie (N.S.) 50(98), no. 1, 2007, 3-31.
  \bibitem{mir4}Mircea Cimpoea\c s "A stable property of Borel type ideals", to appear in Communications in algebra. 
  \bibitem{mir5}Mircea Cimpoea\c s "A note on the generic initial ideal for complete intersections", Bull. Math. Soc.
                Sci. Math. Roumanie, tome 50 (98), no. 2, 2007, p. 119-130.
  \bibitem{mir6}Mircea Cimpoea\c s "Regularity for certain classes of monomial ideals", Analele Stiintifice ale
               Univ.Ovidius, seria matematica, Vol. XV (2007) fascicola 1.
  \bibitem{mir7}Mircea Cimpoea\c s "Some remarks on Borel type ideals", Preprint 2007.
  \bibitem{CH} A.\ Conca, J.\ Herzog "Castelnuovo-Mumford regularity of products of ideals",
               Collect. Math.  54  (2003),  no. 2, 137-152.
  \bibitem{CHH}A.\ Conca, J.\ Herzog, T.\ Hibi "Rigid resolutions and big Betti numbers", Comment. Math. Helv. 79
              (2004), 826-839.  
  \bibitem{CSi}A.\ Conca, J.\ Sidman "Generic initial ideals of points and curves", Preprint 2005,
               arXiv:math.AC/0402418.
  \bibitem{ert}D.\ Eisenbud, A.\ Reeves, B.Totaro "Initial ideals, veronese subrings and rates of algebras", Adv.Math.  
               109 (1994), 168-187.
  \bibitem{E}  D.\ Eisenbud "Commutative algebra with a view toward algebraic geometry", Springer 1994.
  \bibitem{F}  R.\ Fr\"oberg "An inequality for Hilbert series of graded algebras", Mathematica Scandinavica 56, 1985,
               117-144.
  \bibitem{hmnw}T.\ Harima, J.\ Migliore, U.\ Nagel, J.\ Watanabe "The weak and strong Lefschetz properties for artinian
               $K$-algebras", J. Algebra 262 (2003), no. 1, 99--126.
  \bibitem{HW} T.\ Harima, J.\ Watanabe, "The finite free extension of $K$-algebras with the strong Lefschetz property",
               Rend. Sem. Mat. Univ. Padova, 110 (2003), 119-146.
  \bibitem{HP} J.\ Herzog, D.\ Popescu "The strong Lefschetz property and simple extensions", Preprint.  
  \bibitem{hpv}J\"urgen Herzog, Dorin Popescu, Marius Vladoiu "On the Ext-Modules of ideals of Borel type",
               Contemporary Math. 331 (2003), 171-186.
  \bibitem{hp} J\"urgen Herzog, Dorin Popescu "On the regularity of p-Borel ideals", Proceed.of AMS, Volume 129,
               no.9, 2563-2570.    
  \bibitem{hp-f}J\"urgen Herzog, Dorin Popescu "Finite filtrations of modules and shellable multicomplexes",
               Manuscripta Math., 121(2006), 385-410.
  
  \bibitem{Mor}G.\ Moreno, "These", Ecole Politechnique, 1991.
  \bibitem{Par}K.\ Pardue "Generic polynomials", Preprint 1999.
  \bibitem{par}Keith Pardue, "Non standard Borel fixed ideals", Dissertation, Brandeis University, 1994.  
  \bibitem{p}  D.\ Popescu "Extremal Betti numbers and regularity of Borel type ideals", Bull. Math. Soc. Sc. Math.
               Roum. 48(96), no 1, (2005), 65-72.
  \bibitem{P}  D.\ Popescu "The strong Lefschetz property and certain complete intersection extensions", Bull.Math
               Soc.Sc.Math.Roumanie, 48(96),no 4,(2005),421-431.
  \bibitem{PV} D.\ Popescu, M.\ Vladoiu "Strong Lefschetz property on algebras of embedding dimension three",
               Bull. Math. Soc. Sci. Math. Roumanie (N.S.) 49(97) (2006), no. 1, 75--86.
  \bibitem{S}  R.\ P.\ Stanley "Combinatorics and Commutative Algebra", Birkh\"auser, 1983.
  \bibitem{V}  R.\ Villareal "Monomial Algebras", Monographs and Textbooks in Pure and Applied Mathematics, vol. 238,
               Marcel Dekker Inc., New York, 2001
  \bibitem{W}  A.\ Wiebe "The Lefschetz property for componentwise linear ideals and Gotzmann ideals",
               Comm. Algebra 32 (2004), no. 12, 4601--4611. 
\end{thebibliography}
\end{document}